\documentclass[leqno]{amsart}
\usepackage{amsmath,amssymb,amsthm,amsfonts, amscd}
 \usepackage{amsaddr}
\usepackage{color}
\newtheorem{RemarkA}{Remark}
\newtheorem{theoremA}{Theorem}
\numberwithin{equation}{section}

\newtheorem{proposition}{Proposition}[section]

\newtheorem{theorem}{Theorem}[section]

\newtheorem{corollary}{Corollary}[section]

\newtheorem{Remark}{Remark}[section]

\makeatletter
\renewcommand{\email}[2][]{%
  \ifx\emails\@empty\relax\else{\g@addto@macro\emails{,\space}}\fi%
  \@ifnotempty{#1}{\g@addto@macro\emails{\textrm{(#1)}\space}}%
  \g@addto@macro\emails{#2}%
}
\makeatother
\title{On gamma factors of Rankin--Selberg integrals for $\mathrm{U}_{2\ell} \times \mathrm{Res}_{E \slash F} \mathrm{GL}_n$}
\author{Kazuki Morimoto}
\thanks{The research of the second author was supported in part by
JSPS KAKENHI Grant Number 17K14166, 21K03164}
\date{\today}   
\email{morimoto@math.kobe-u.ac.jp}                                        
\subjclass[2010]{Primary: 11S40; Secondary: 11F70}
\begin{document}
\begin{abstract}
In this paper, we prove the fundamental properties of gamma factors defined by Rankin-Selberg integrals of Shimura type for 
pairs of generic representations $(\pi, \tau)$ of $\mathrm{U}_{2\ell}(F)$ and  $\mathrm{GL}_n(E)$
for a local field $F$ of characteristic zero and a quadratic extension $E$ of $F$.
We also prove similar results for pairs of generic representations $(\pi, \tau_1 \otimes \tau_2)$
of $\mathrm{GL}_{2\ell}(F)$ and $\mathrm{GL}_n(F) \times \mathrm{GL}_n(F)$.
As a corollary, we prove that the gamma factors arising from Langlands--Shahidi method
and our gamma factors coincide.
\end{abstract}
\maketitle
%%%%%%%%%%%%%%%%%%%%%%%%%%%%%%%%%%%%%%%%%%%%%%%
%
%
%
%
%
%
%
%
%
%
%
%
%
%
%%%%%%%%%%%%%%%%%%%%%%%%%%%%%%%%%%%%%%%%%%%%%%%
\section{Introduction}
Let $F$ be a local field of characteristic zero and $E$ be a \'{e}tale quadratic extension of $F$, i.e. $E$ is a quadratic extension of $F$ or $E= F \times F$.
Let $(\pi, \tau)$ be a pair of irreducible generic representations of quasi-split even unitary group $\mathrm{U}_{2\ell}(F)$ and $\mathrm{GL}_n(E)$.
When $E=F \times F$,  $\mathrm{U}_{2\ell}(F) \simeq \mathrm{GL}_{2\ell}(F)$ and 
$\mathrm{GL}_n(E) \simeq \mathrm{GL}_n(F) \times \mathrm{GL}_n(F)$, and thus $\pi$ is an irreducible representation of $\mathrm{GL}_{2\ell}(F)$ and
$\tau=\tau_1 \otimes \tau_2$ with irreducible representations $\tau_i$
of $\mathrm{GL}_n(F)$. 
In this paper, we study a family of Rankin--Selberg integrals of Shimura type for $(\pi, \tau)$ introduced and studied in Ben-Artzi and Soudry~\cite{BAS2,BAS}. When $\ell=n$, these local integrals were also studied by Watanabe~\cite{Wat}.
Note that Ginzburg, Rallis and Soudry~\cite{GRS2} studied global Rankin--Selberg integrals of Shimura type when $\ell < n$.

Let us set some notation to state our main results.
We fix an additive character $\psi_F$ of $F$ and define a character $\psi$ of $E$ by 
$x \mapsto \psi_F\left( x+\bar{x} \right)$ with the action $x \mapsto \bar{x}$ of non-trivial element of $\mathrm{Gal}(E \slash F)$.
Let $\omega_{E \slash F}$ denote the (possibly trivial) quadratic character of $F^\times$
corresponding to $E \slash F$.
Then we fix a character $\Upsilon$ of $E^\times$ such that $\Upsilon |_{F^\times} = \omega_{E \slash F}$.
Suppose that $\pi$ is an  irreducible $\chi$-generic representation of $\mathrm{U}_{2\ell}(F)$ with a non-degenerate character $\chi$
and $\tau$ is an  irreducible generic representation of $\mathrm{GL}_n(F)$.
Let $Q_n(F)$ be a maximal parabolic subgroup of $\mathrm{U}_{2n}(F)$ whose Levi part is isomorphic to $\mathrm{GL}_n(E)$.
For $s \in \mathbb{C}$, we consider the parabolic induction $\mathrm{Ind}^{\mathrm{U}_{2n}(F)}_{Q_n(F)}(\tau |\det|^{s-\frac{1}{2}})$.
Let us denote by $M^\ast(\tau, s) : \mathrm{Ind}^{\mathrm{U}_{2n}(F)}_{Q_n(F)}(\tau |\det|^{s-\frac{1}{2}}) \rightarrow 
\mathrm{Ind}^{\mathrm{U}_{2n}(F)}_{Q_n(F)}(\tau^\ast |\det|^{-s+\frac{1}{2}})$
the normalized intertwining operator
where $\tau^\ast$ denotes the conjugate of the contragredient $\tau^\vee$ of $\tau$.

The Rankin-Selberg integral $\mathcal{L}(W, f_s, \phi)$ is defined for a Whittaker function $W$ for $\pi$, the Weil representation attached to $\Upsilon$ and $\psi$,
a Schwartz function $\phi \in \mathcal{S}(E^{\mathrm{min}(\ell, n)})$ and $f_s \in \mathrm{Ind}^{\mathrm{U}_{2n}(F)}_{Q_n}(\tau |\det|^{s-\frac{1}{2}})$.
Then $\mathcal{L}(W, f_s, \phi)$ and $\mathcal{L}(W, M^\ast(\tau,s)f_s, \phi)$  satisfy the same certain equivalence property.
These integrals are related by a proportionality factor $\gamma(s, \pi \times \tau, \Upsilon, \psi)$ such that 
\[
\mathcal{L}(W, M^\ast(\tau,s)f_s, \phi) = \gamma(s, \pi \times \tau, \Upsilon, \psi)\mathcal{L}(W, f_s, \phi).
\]
Define 
\[
\Gamma(s, \pi \times \tau, \Upsilon, \psi) :=\omega_\pi(-1)^n \omega_\tau(-1)^\ell  \gamma(s, \pi \times \tau, \Upsilon, \psi)
\]
where $\omega_\pi$ and $\omega_{\tau}$ denote the central characters of $\pi$ and $\tau$, respectively.
When $E = F \times F$, we may write $\Upsilon= \mu \otimes \mu^{-1}$ with a character $\mu$ of $F^\times$.
In this case, for $\tau= \tau_1 \otimes \tau_2$, we often write our $\gamma$-factor and the normalized $\gamma$-factor
by $ \gamma(s, \pi \times (\tau_1 \otimes \tau_2), \mu, \psi)$ and $\Gamma(s, \pi \times (\tau_1 \otimes \tau_2), \mu, \psi)$, respectively.
%Similarly, for irreducible representations $\pi$ and $\tau_1 \otimes \tau_2$ of $\mathrm{GL}_{2\ell}(F)$ and $\mathrm{GL}_n(F) \times \mathrm{GL}_n(F)$
%and a character $\mu$ of $F^\times$, we can define normalized gamma factor $\Gamma(s, \pi \times (\tau_1 \otimes \tau_2), \mu, \psi)$.
Our aim in this paper is to study basic properties of $\gamma(s, \pi \times \tau, \Upsilon, \psi)$
and $ \gamma(s, \pi \times (\tau_1 \otimes \tau_2), \mu, \psi)$.
%and $\Gamma(s, \pi \times (\tau_1 \otimes \tau_2), \mu, \psi)$. 

Similar local zeta integrals for orthogonal groups were introduced by Ginzburg~\cite{Gi} and studied by Soudry~\cite{S1, S2, S3} and Kaplan~\cite{Ka2}. 
Ginzburg-Rallis-Soudry~\cite{GRS98} introduced and studied similar local zeta integrals for symplectic groups and metaplectic groups.
These local zeta integrals and our zeta locals are important ingredients to study local and global descent (see \cite{GRS2} 
for global descent and Ginzburg-Raliis-Soudry~\cite{GRS99} and Soudry-Tanay~\cite{ST} for local descent).

Kaplan~\cite[Theorem~1]{Ka} proved a list of fundamental properties
for $\gamma$-factors defined by similar local zeta integrals in \cite{GRS2} for orthogonal groups, symplectic groups,
metaplectic groups using results and ideas in \cite{S1,S2,S3} and \cite{Ka2}.
Recently, Cheng-Wang~\cite{CW} proved a similar result for odd unitary groups.
%Lapid and Rallis~\cite{LR} proved a full list of fundamental properties 
%for $\gamma$-factors defined by the doubling method.
%Moreover, 
In this paper, we prove a similar list of fundamental properties as \cite[Theorem~1]{Ka} for our $\gamma$-factors.
See Theorem~\ref{main thm} for a precise list.
\begin{theorem}
\label{main thm intro}
The $\gamma$-factor $\Gamma(s, \pi \times \tau, \Upsilon, \psi)$ satisfies a list of fundamental properties
including multiplicativity, unramified cases, minimal cases, Archimedean property and a global identity.
Further, these properties define them uniquely. 
\end{theorem}
Recall that most of these properties were stated in Shahidi~\cite{Sh2} for $\gamma$-factors defined by local coefficients.
Also, he showed that those properties define his $\gamma$-factors uniquely.
Further, recall that the full list of fundamental properties of $\gamma$-factors was stated in Lapid--Rallis~\cite{LR} for
 the local zeta integral given by the doubling method.
 For other instances studying basic properties of local $\gamma$-factors, see the introduction of \cite{Ka}.

The main part of the proof of Theorem~\ref{main thm intro} is to prove multiplicativity properties. 
There are proved here following the idea used in Soudry~\cite{S1,S2,S3} and Kaplan~\cite{Ka,Ka2}.

As stated in Theorem~\ref{main thm intro}, the above list defines the $\gamma$-factors uniquely.
Hence, as a corollary of this theorem, we obtain the following consequence.
\begin{corollary}
\label{1}
Our $\gamma$-factors $\Gamma(s, \pi \times \tau, \Upsilon, \psi)$ and $\Gamma(s, \pi \times (\tau_1 \otimes \tau_2), \mu, \psi)$ are identical with Shahidi's $\gamma$-factors
$\gamma^{Sh}(s, (\pi \otimes \Upsilon) \times \tau, \psi )$ and $\gamma^{Sh}(s, (\pi \otimes \mu) \times \tau_1,\psi) \gamma^{Sh}(s, (\pi  \otimes \mu) \times \tau_2, \psi)$, respectively. Here, $\gamma^{Sh}(s, (\pi \otimes \Upsilon) \times \tau, \psi )$ and $\gamma^{Sh}(s, (\pi \otimes \mu) \times \tau_i,\psi)$
denote Shahidi's $\gamma$-factors \cite{Sh2}.
%This implies that in the split case, $\Gamma(s, \pi \times (\tau_1, \tau_2), \mu, \psi)$ is equal to $\gamma(s, \pi \times (\tau_1 \otimes \mu), \psi) 
%\gamma(s, \pi \times (\tau_2 \otimes \mu^\ast), \psi)$
%Moreover, $\Gamma(s, \pi \times (\tau_1, \tau_2), \mu, \psi) = \gamma(s, \pi \times (\tau_1 \otimes \mu^{-1}), \psi)  \gamma(s, \pi \times (\tau_2 \otimes \mu), \psi)$
%where the right-hand side is the gamma factor defined by Jacquet, Piatetski-Shapiro and Shalika~\cite{JPSS}.
\end{corollary}
We have a similar equality of $\gamma$-factors for metaplectic groups defined by Rankin-Selberg integrals of Shimura type 
and $\gamma$-factors defined by Langlands-Shahidi method.
This equality is used in a proof of a refined form of the formal degree conjecture by Ichino--Lapid--Mao~\cite{ILM}.
The above equality is also used in the author's work~\cite{Mo}
for a proof of  a refined form of the formal degree conjecture for even unitary groups.
Note that Beuzart-Plessis~\cite{BP} proved the formal degree conjecture for unitary groups in a different method.
%We also note another application of our result.
%For an irreducible generic representation of $\mathrm{GL}_n$, Jacquet conjectured a precise family which determines
%by the coincidence of Rankin-Selberg gamma factors. This conjecture was recently proved by Chai~\cite{Cha} and  independently 
%by Jacquet-Liu~\cite{JL}. More recently, Zhang~\cite{Zha,Zha2} give another proof of Jacquet's conjecture
%assuming the coincidence of gamma factors given in Corollary~\ref{1}.
%Then our result complete his proof of Jacquet's conjecture.

This paper is organized as follows.
In Section~\ref{section:notation}, we set basic notation and recall the definition of the local zeta integrals.
In Section~\ref{section:Main Theorem}, we state our main result, namely the list of fundamental properties of $\gamma$-factors.
Also, we prove fundamental properties except for the multiplicativity.
In Section~\ref{proof:s}, we prove the multiplicativity of $\gamma$-factors.
In \ref{section:app}, we prove uniqueness of Fourier-Jacobi modules for certain representations of Whittaker type over non-archimedean local fields.
\subsection*{Acknowledgements}
This paper is part of the author's project to prove a conjecture of Lapid and Mao and the formal degree conjecture 
in the unitary group case.
The author expresses his deep gratitude to Erez Lapid for suggesting him to the project.
The author is very grateful to the anonymous referee for reading the paper carefully, pointing out several inaccuracies and making many useful suggestions.
%For special orthogonal groups this has been achieved, essentially, in[51--53,34], 
%but a few details were missing.
%Once proved, multiplicativity, along with the computations of the γ-factors for unramified data, can be used in a standard global argument 
%(asin[50,Section5],[53,Section0],[39,Section9]and[13,Sections10,11]),
%to derive most parts of the theorem. The computations for unramified data were conducted in all cases (including the metaplectic),
%in[14,16,51,18,32,33].
%The bulk of the work here is to prove multiplicativity properties for the symplectic and metaplectic groups. 
%To this end we will use the tools and ideas of [51--53,18,19,34].

%%%%%%%%%%%%%%%%%%%%%%%%%%%%%%%%%%%%%%%%%%%%%%%
%
%
%
%
%
%
%
%
%
%
%
%
%
%
%%%%%%%%%%%%%%%%%%%%%%%%%%%%%%%%%%%%%%%%%%%%%%%
\section{Notation}
\label{section:notation}
%\begin{align*}
%U_r &= \text{unipotent radical of $G_r$ whose Levi part is $\mathrm{GL}_r$}
%\\
%Z_{n_1, n_2}&= \left\{\begin{pmatrix}1_{n_1}&z\\ &1_{n_2} \end{pmatrix} \right\}
%\\
%\Omega_{a, b} &= \begin{pmatrix}&1_{a}\\ 1_b& \end{pmatrix}
%\end{align*}
%We follow the notation in \cite{BAS}.
\subsection{Groups and elements}
Let $F$ be a local field of characteristic zero.
Let $E$ be an \'etale quadratic extension of $F$, so that $E$ is either a quadratic extension of $F$ or $F \times F$.
We denote by $E \ni x \mapsto \bar{x}$ the action of  non-trivial element of $\mathrm{Gal}(E \slash F)$.
We denote by $N_{E \slash F}$ and $\mathrm{Tr}_{E \slash F}$
the norm map and the trace map from $E$ to $F$.
Let $\omega_{E \slash F}$ denote the (possibly trivial) quadratic character of $F^\times$
corresponding to $E \slash F$.
Fix a character $\Upsilon$ of $E^\times$ such that $\Upsilon |_{F^\times} = \omega_{E \slash F}$.

Then for a positive integer $k$, we consider the quasi-split even unitary group $G_k$,
which is realized by
\[
G_k  := \{ g \in \mathrm{GL}_{2k}(E) : {}^{t} \bar{g} J_{2k} g = J_{2k} \}
\]
where ${}^{t} \bar{g}$ denotes the transpose of $\bar{g}$,
\[
J_{2k} = \begin{pmatrix}&w_k\\ -w_k \end{pmatrix},
\]
and $w_k$ is the $k \times k$ permutation matrix, which has $1$ on the main skew diagonal.
For $k_0 < k$, we regard $G_{k_0}$ as a subgroup of $G_k$ by the embedding
\[
g \mapsto \begin{pmatrix} I_{k-k_0}&&\\ &g&\\ &&I_{k-k_0} \end{pmatrix}, \quad g \in G_{k_0}.
\]
When $E = F \times F$, $\mathrm{GL}_i(E)$ is isomorphic to $\mathrm{GL}_i(F) \times \mathrm{GL}_i(F)$ and 
$G_k$ is isomorphic to $\mathrm{GL}_{2k}(F)$ (see \cite[Section~1.3]{GRS2} for a concrete isomorphism). 

Let $U_{G_k}$ be the group of upper triangular unipotent matrices of $G_k$.
Let $T_k$ be the subgroup of $G_k$ consisting of diagonal matrices.
Then $T_k U_{G_k}$ is a Borel subgroup of $G_k$.
We denote by $U_{\mathrm{GL}_k}$ and $Z_k$ the group of upper triangular unipotent matrices of $\mathrm{GL}_k(E)$
and $\mathrm{GL}_k(F)$, respectively.
Let $Q_r = Q^{k}_r$ be the standard parabolic subgroup of $G_k$ whose
Levi part $M_r = M^k_r$ is isomorphic to $\mathrm{GL}_{r}(E) \times G_{k-r}$.
We denote the unipotent radical of $Q_r$ by $U_r = U_r^k$.
%When $E = F \times F$, $Q_r$ denotes the standard parabolic subgroup of $\mathrm{GL}_{2k}(F)$
%whose Levi part $M_r$ is isomorphic to $\mathrm{GL}_{r}(F) \times \mathrm{GL}_{2(k-r)}(F) \times \mathrm{GL}_r(F)$.
For a standard parabolic subgroup $P$ of one of the groups $G_k$, $\mathrm{GL}_i(E)$ or $\mathrm{GL}_i(F)$, 
we denote by $\overline{P}$ the opposite parabolic containing the Levi part of $P$. Then we have $\overline{P}=\{ {}^t p : p \in P\}$.
We write unipotent radical of $\bar{P}$ by $\overline{U}$. Then we have
$\bar{U}=\left\{ {}^t u : u \in U\right\}$ with unipotent radical $U$ of $P$.
We often identify $F^\times$ with the center of $\mathrm{GL}_i(F)$ by $a \mapsto a \cdot I_i$.

For any $g \in \mathrm{GL}_{r}(E)$, we put $g^\ast = w_ {r} {}^{t} \bar{g}^{-1} w_{r}$.
Let $\psi_F$ be a non-trivial character of $F$ and we extend it to a character of $E$
by $E \ni x  \mapsto \psi_F \left(\mathrm{Tr}_{E \slash F}(x)\right)$.
Diagonal or block-diagonal matrices will be denoted $\mathrm{diag}(\cdots)$.
For elements $x$ and $y$ of a group and a subgroup $Y$,
denote ${}^{x} y = x^{-1} y x$ and ${}^{x}Y = \{ {}^{x}y : y \in Y \}$.
%For $a \in \mathrm{GL}_r$, put $a^\ast = w_r {}^{t}a^{-1} w_r$.
%%%%%%%%%%%%%%%%%%%%%%%%%%%%%%%%%%%%%%%%%%%%%%%
%
%
%
%
%
%
%
%
%
%
%
%
%
%
%%%%%%%%%%%%%%%%%%%%%%%%%%%%%%%%%%%%%%%%%%%%%%%
\subsection{Representations}
\label{represen:s}
Throughout this paper, by a representation, we always means a smooth, admissible and finitely generated complex representation.
Further, when $F$ is an archimedean field, it is assumed to be on a Fr\'{e}chet space and of moderate growth.

Let $H$ be one of the groups of $G_k, \mathrm{GL}_k(E)$ or $\mathrm{GL}_k(F)$.
Let $V$ be a maximal unipotent subgroup of $H$.
Let $(\pi, V_\pi)$ be a (possibly reducible) representation of $H$.
A linear functional $\lambda$ on $V_\pi$ is called a Whittaker functional with respect to a generic character $\chi$ of $V$
if $\lambda(\pi(u)v) = \chi(u) \lambda(v)$ for all $u \in V$ and $v \in V_\pi$. 
If the field is archimedean, we further require that $\lambda$ is continuous. 
We say that $\pi$ is $\chi$-generic, if the complex vector space of these functionals is non-zero.
Further, if this space is one-dimensional, we say that $\pi$ is of $\chi$-Whittaker type.
In this case, we define the space $\mathcal{W}(\pi, \chi)$ by  the space of complex-valued functions $W_v$ on $H$, 
defined for some $v \in V_\pi$ by $W_v(g) = \lambda(\pi(g)v)$. 
If $\pi$ is of $\chi$-Whittaker type, the representation of $H$ by right-translations 
on $\mathcal{W}(\pi, \chi)$  is a quotient of $\pi$. We call $\mathcal{W}(\pi, \chi)$ a $\chi$-Whittaker space of $\pi$. 
If $\pi$ is irreducible, $\mathcal{W}(\pi, \chi)$ is the $\chi$-Whittaker model of $\pi$.
When a choice of a character $\chi$ is clear, we simply omit $\chi$ from the notation, e.g. generic instead of $\chi$-generic.
\begin{center}
\emph{Throughout this paper, we assume that all representations are of Whittaker type.}
\end{center}
If a representation $\pi$ has a central character, we denote it by $\omega_\pi$.
The contragredient representation of $\pi$ is denoted by $\pi^\vee$.
%
%Let $V_{Q_n}^{\mathrm{U}_{2n}} \left( \mathcal{W}(\tau, \psi), s \right)$ be the space of $\mathrm{U}_{2k}$-smooth
%left $U_n$-invariant functions $W : \mathrm{U}_{2n} \rightarrow \mathbb{C}$
%such that for all $g \in \mathrm{GL}_{2n}$, the function $m \mapsto \delta_{Q_n}(m)^{-\frac{1}{2}} |\det m|^s W(mg)$
%on $M_n$ belongs to $\mathcal{W}(\pi, \psi_{N_M})$.
%

First, let us consider the non-split case, namely the case where $E$ is a quadratic extension of $F$.
Let $H_k := E^{2k} \oplus F$ be the Heisenberg group in $4k+1$-variable over $F$
with the multiplication
\[
(u; z) \cdot (v; z^\prime) = \left(u+v; z+z^\prime+\frac{1}{2} \mathrm{Tr}_{E \slash F}(uJ_{2k} {}^{t} \bar{v}) \right),
\quad u, v \in E^{2k}, z, z^\prime \in F
\]
(see \cite[Section~5]{BAS2}).
Then the center of $H_k$ is the subgroup consisting of elements $(0; t)$ with $t \in F$.
%We may embed $H_k$ into $G_{k+1}$ by
%\[
%(u; r) \mapsto \begin{pmatrix}1&v&2r+v J_{2k} {}^{t}\bar{v}\\ &1_{2k}&v^\prime\\ &&1 \end{pmatrix}
%\]
%where $v^\prime = J_{2k} {}^{t} \bar{v}$. 
The group $G_k$ acts on $H_k$ on the right by
\[
(v; z) \cdot g = (vg; z), \quad v \in E^{2k}, z \in F, g \in G_k.
\]
We embed $H_k$ into $G_{k+1}$ by
\[
(v ; z) \mapsto \begin{pmatrix}1&v&z+\frac{1}{2}v J_{2k} {}^{t}\bar{v} \\ &I_{2k}&J_{2k}{}^{t}\bar{v}\\ &&1\end{pmatrix}
\]
where $v \in E^{2k}$ and $z \in F$.
%This action is translated to conjugation via the above embedding:
%\[
%(vg; r) \mapsto \begin{pmatrix}1&&\\ &g&\\ &&1 \end{pmatrix}^{-1}
% \begin{pmatrix}1&v&r+\frac{1}{2} v J_{2k} {}^{t}\bar{v}\\ &1_{2k}&v^\prime\\ &&1 \end{pmatrix}
% \begin{pmatrix}1&&\\ &g&\\ &&1 \end{pmatrix}.
%\]
Let us recall the formula the Weil representation $\omega_{\psi, \Upsilon}$ of $G_k \ltimes H_k$
attached to $\psi$ and $\Upsilon$ acting on the Schwartz space $\mathcal{S}(E^k)$ on $E^k$.
It will be convenient to write an element of $H_k$
as $(x, y; t)$ with $x, y \in E^k$ and $t \in F$.

Let $\phi \in \mathcal{S}(E^k)$.
Then the action of $H_k$ is given by the following formulas: for $\xi  \in E^k$, 
\[
\omega_{\psi, \Upsilon}((x, 0; 0))\phi(\xi) = \phi(\xi+x), \quad x \in E^k
\]
\[
\omega_{\psi, \Upsilon}((0, y; 0))\phi(\xi) = \psi \left( \xi w_k {}^{t}\bar{y} \right) \phi(\xi), \quad y \in E^k
\]
\[
\omega_{\psi, \Upsilon}((0, 0; z))\phi(\xi) = \psi \left(\frac{z}{2} \right) \phi(\xi), \quad  t \in F.
\]
On the other hand, the action of $G_k$ is given by the following formulas: for $\xi \in E^k$, 
\[
\omega_{\psi, \Upsilon} \left( \begin{pmatrix} a&\\ &a^\ast \end{pmatrix} \right)\phi(\xi)
= \Upsilon(\det (a)) |\det (a)|^{\frac{1}{2}} \phi(\xi a),
\]
\[
\omega_{\psi, \Upsilon} \left( \begin{pmatrix} 1&X\\ &1\end{pmatrix} \right)\phi(\xi)
= \psi \left(\frac{1}{2} \left(\xi w_k {}^{t}(\overline{\xi X}) \right) \right) \phi(\xi),
\]
\[
\omega_{\psi, \Upsilon} \left(J_{2k} \right)\phi(\xi)
= \hat{\phi}(\xi):= \int_{E^k} \phi(x) \psi(\bar{x} {}^{t} \xi) \, dx
\]
where $a \in \mathrm{GL}_k(E)$, $X \in \mathrm{Mat}_{k \times k}(E)$ such that ${}^{t}(w_k \overline{X}) = w_k X$
and the measure on $E^k$ is the self-dual Haar measure with respect to $\psi$.

In the split case, namely $E= F \times F$, we have a similar formula for the Weil representation of $\mathrm{GL}_{2k} \ltimes H_k$.
In this case, the Heisenberg group $ H_k :=F^{2k} \oplus F^{2k} \oplus F$ is defined with the multiplication
\[
((u_1, u_2); z) \cdot ((v_1, v_2); z^\prime) = ((u_1+v_1, u_2+v_2); z+z^\prime+\frac{1}{2}(u_1w_{2k}{}^{t} v_2-v_1w_{2k}{}^{t} u_2))
\]
for $u_i, v_i \in F^{2k}$, $z, z^\prime \in F$ (see \cite[Section~5]{BAS2}). 
%We embed $\mathcal{H}_{4n}$ into $\mathrm{GL}_{2k+2}(F)$ by
%\[
%((u_1, u_2); t) \mapsto \begin{pmatrix} 1&u_1&t+\frac{1}{2} u_1J_{2k} {}^{t}u_2\\ &1_{2k}&J_{2k} {}^{t} u_2 \\ &&1\end{pmatrix}
%\]
The action of $\mathrm{GL}_{2k}(F)$ on $H_k$ on the right is given by
\[
((u_1, u_2); z) \cdot g = ((u_1 g, u_2 {}^\ast g); z)
\]
where for $g \in \mathrm{GL}_{2a}(F)$, we define ${}^\ast g = J_{2k} {}^{t} g^{-1} J_{2k}^{-1}$.
For later use, we extend $\ast$ to $\mathrm{GL}_{2k}(E) = \mathrm{GL}_{2k}(F) \times \mathrm{GL}_{2k}(F)$ by $(g_1, g_2)^\ast = (g_1^\ast, g_2^\ast)$.
%Note that this action corresponds to the conjugation by $\mathrm{diag}(1, g, 1)$ as in the non-split case.

In this case,  we may define the Weil representation $\omega_{\psi, \Upsilon}$ of $\mathrm{GL}_{2k}(F) \ltimes H_k$
attached to $\psi$ and $\Upsilon$. Since $E^\times = F^\times \times F^\times$, 
$\Upsilon=\mu \otimes \mu^\prime$ with characters $\mu, \mu^\prime$ of $F^\times$.
By the condition $\Upsilon |_{F^\times} =1$,  we should have $ \mu^\prime=\mu^{-1}$ and $\Upsilon=\mu \otimes \mu^{-1}$. 
Hence, in this case, we often simply write this Weil representation 
by $\omega_{\psi, \mu}$.
Then it is realized on the space $\mathcal{S}(F^k \times F^k)$ by the following formulas.
Indeed, for $\phi \in \mathcal{S}(F^k \times F^k)$, $z \in F$ and $x_i, y_i, \xi_i \in F^k$,
 the action of $H_k$ is given by the following formulas:
\[
\omega_{\psi, \mu} (((x_1, 0), (x_2, 0)); z) \phi(\xi_1, \xi_2)
= \psi_F \left(z \right) \phi(\xi_1+x_1, \xi_2+x_2)
\]
\[
\omega_{\psi, \mu} (((0, y_1), (0, y_2)); 0) \phi(\xi_1, \xi_2)
= \psi_F \left( \xi_1 w_k {}^{t} y_2 + \xi_2 w_k {}^{t} y_1 \right) \phi(\xi_1, \xi_2).
\]
Further, the action of $\mathrm{GL}_{2k}(F)$ is given as follows:
\[
\omega_{\psi, \mu} \left( \begin{pmatrix} a&\\ &b \end{pmatrix}\right)\phi(\xi_1, \xi_2)
=\mu(\det(ab)) \left| \frac{\det(a)}{\det(b)} \right|^{\frac{1}{2}} \phi(\xi_1 a, \xi_2 b^\ast), \quad a, b \in \mathrm{GL}_k(F);
\]
\[
\omega_{\psi, \mu} \left( \begin{pmatrix} 1_k&S\\ &1_k \end{pmatrix}\right)\phi(\xi_1, \xi_2)
= \psi_F \left(  \xi_1 S w_k {}^{t} \xi_2  \right) \phi(\xi_1, \xi_2), \quad S \in \mathrm{Mat}_{k \times k}(F)
\]
\[
\omega_{\psi, \mu} \left(J_{2k} \right)\phi(\xi)
= \hat{\phi}(\xi_1, \xi_2):= \int_{(F^k)^2} \phi(x_1, x_2) \psi_F \left(x_2 {}^{t} \xi_1 +x_1 {}^{t} \xi_2 \right) \, dx_1 \, dx_2,
\]
where $\xi_i, x_i \in F^k$ and $b^\ast = w_k {}^{t} b^{-1} w_k$.
We embed $H_k$ into $\mathrm{GL}_{2k}(F)$ by
\[
((x_1, x_2), (y_1, y_2); z) \mapsto \begin{pmatrix} 1&x_1&y_1&z+\frac{1}{2}(x_1w{}^{t}y_2-x_2 w{}^{t}y_1) \\ &I_k&&w{}^{t}y_2\\ &&I_k&-w{}^{t}x_2
\\&&&1 \end{pmatrix}
\]
with $x_i, y_i \in F^{k}$ and $z \in F$.
%%%%%%%%%%%%%%%%%%%%%%%%%%%%%%%%%%%%%%%%%%%%%%%
%
%
%
%
%
%
%
%
%
%
%
%
%
%
%%%%%%%%%%%%%%%%%%%%%%%%%%%%%%%%%%%%%%%%%%%%%%%
\subsection{Local zeta integrals} 
\label{s:def local zeta}
Let us recall the definitions of the local zeta integrals of Shimura type defined in \cite{BAS2,BAS}
and their associated gamma factors. 
Let $\ell$ and $n$ be positive integers.
%Let $N_{\ell, n}$ (resp $N_{\ell, n}^\circ$) be the unipotent radical of the standard parabolic subgroup of $G_n$
%whose Levi part is $\mathrm{GL}_1(E)^{n-\ell} \times G_\ell$ (resp. $\mathrm{GL}_1(E)^{n-\ell-1} \times G_{\ell+1}$).
%Then the map
%\[
%v \mapsto v_{\mathcal{H}} := \left((v_{n-\ell, n-\ell+j})_{j=1, \dots, 2\ell}, \frac{1}{2} \mathrm{Tr}_{E \slash F}(v_{n-\ell, n+\ell+1}) \right)
%\]
%gives an isomorphism from $N_{\ell, n}^\circ \backslash N_{\ell, n}$ to the Heisenberg group $\mathcal{H}_\ell$.
Let $\pi$ and $\tau$ be a pair of admissible representations of $G_\ell$ and $\mathrm{GL}_{n}(E)$, respectively.
Llet $\psi_{U_{G_k}} $ be a non-degenerate character of $U_{G_k}$ defined by
\[
\psi_{U_{G_k}}(n) = \psi \left( \sum_{i=1}^{k-1}n_{i, i+1} -\frac{1}{2} n_{k, k+1} \right),
\qquad n \in U_{G_k}.
\]
When $E= F \times F$, from the concrete isomorphism \cite[Section~1.3]{GRS2},
we have 
\begin{equation}
\label{gen char G_k}
\psi_{U_{G_k}}(n) =  \psi_F \left( z_{1,2}+ \cdots +z_{k-1, k} - z_{k, k+1}- \cdots -z_{2k-1, 2k} \right)
\end{equation}
for $n \in U_{G_k}$.
\begin{Remark}
Suppose that $E$ is a quadratic extension of $F$.
Recall that $T_k$ acts on the set of non-degenerate characters of $U_{G_k}$ by a conjugation.
Then any non-degenerate character is conjugate to a character defined by 
\[
\psi_{U_{G_k}, a}(n) = \psi \left( \sum_{i=1}^{k-1}n_{i, i+1} -\frac{a}{2} n_{k, k+1} \right)
\]
with $a \in F^\times$. 
%Moreover, $\psi_{U_{G_k}, a}$ and $\psi_{U_{G_k}}$ are in the same orbit
%if and only if $a \in N_{E \slash F}(E^\times)$. 
On the other hand, $\psi_{U_{G_k}, a}$ and $(\psi_a)_{U_{G_k}}$ are in the same orbit under this action
where $\psi_a$ is the additive character of $F$ defined by $\psi_a(x) = \psi(ax)$.
Since we do not impose any assumption on $\psi$, it suffices to study generic representations 
with respect to the non-degenerate character $\psi_{U_{G_k}}$.
We note that when $E = F \times F$, there is only one orbit.
\end{Remark}
We define a non-degenerate character of $\psi_{k}$ of  $U_{\mathrm{GL}_k}$  by
\[
\psi_{k}(u) = \psi \left(\sum_{i=1}^{k-1} u_{i, i+1} \right).
\]
When $E = F \times F$, we define a non-degenerate character $\psi_k$ of $Z_k$ by 
\begin{equation}
\label{gen char GL}
\psi_{k}(u) = \psi_F \left(\sum_{i=1}^{k-1} u_{i, i+1} \right)
\end{equation}
Then as in the case of $\psi_{U_{G_k}}$, we have 
\[
\psi_{k}((u, u^\prime)) = \psi_k(u)\psi_k^{-1}(u^\prime).
\]
We often simply write $\psi= \psi_k$ if the index $k$ is clear.
For a representation $\tau$ of $\mathrm{GL}_n(E)$, we denote by $\mathcal{W}(\tau, \psi)$ of $\tau$ the Whittaker space of $\tau$ 
with respect to the character $\psi$.
Let $\tau^\ast$ (resp. $\bar{\tau}$) be a representation of $\mathrm{GL}_n(E)$
define by $\tau^\ast(a) = \tau(a^\ast)$ (resp. $\bar{\tau}(a) = \tau(\bar{a})$).
Then we note that $\tau^\ast \simeq \bar{\tau}^\vee$.

Regard $\tau$ as a representation of $M_k^k$ since $M_k^k \simeq \mathrm{GL}_k(E)$.
Recall that when $E= F \times F$, $\tau = \tau_1 \otimes \tau_2$ with representations $\tau_i$ of $\mathrm{GL}_k(F)$.
Moreover, in this case, the embedding $\mathrm{GL}_k(F) \times \mathrm{GL}_k(F)$ into $M_k^k$ is given by $(g_1, g_2) \mapsto \mathrm{diag}(g_1, g_2^\ast)$.

For $s \in \mathbb{C}$, we form the parabolic induction
$\mathrm{Ind}_{Q_k}^{G_k} \left(\mathcal{W}(\tau, \psi) |\det|^{s-\frac{1}{2}} \right)$ of $G_k$ .
We realize it on a space $V_{Q_k}^{G_k} \left( \mathcal{W}(\tau, \psi), s \right)$ of smooth,
complex-valued functions $f_s(-, -)$ on $G_k \times \mathrm{GL}_k(E)$
such that 
\[
f_s(mug, a) = \delta_{Q_k}^{\frac{1}{2}}(m) |\det m|^{s-\frac{1}{2}} f_s(g, am)
\]
and the function $a \mapsto f_s(g, a)$ belongs to $\mathcal{W}(\tau, \psi)$.

For $W \in \mathcal{W}(\pi, \psi_{U_{G_\ell}}^{-1})$, $f_{s} \in V_{Q_n}^{G_n} \left( \mathcal{W}(\tau, \psi), s \right)$
and $\phi \in \mathcal{S}(E^{\mathrm{min}(\ell, n)})$,
we define the local zeta integral $\mathcal{L}(W, f_{s}, \phi)$ by 
\begin{equation}
\label{def local zeta Bessel}
%\mathcal{L}(W, f_{s}, \phi)
%= 
\left\{
\begin{array}{ll}
 \int_{U_{G_\ell} \backslash G_\ell} \int_{Y_\ell \backslash H_\ell} \int_{R_{\ell, n}}
W(g) f_{s}(\beta_{\ell,n} rh g)  \omega_{\psi^{-1}, \Upsilon^{-1}}(h g) \phi(\xi_\ell) \, dr \, dh \, dg & \text{ if $\ell < n$}\\
&\\
 \int_{U_{G_\ell} \backslash G_\ell} 
W(g) f_{s}(g)  \omega_{\psi^{-1}, \Upsilon^{-1}}(g) \phi(\xi_\ell) \, dg & \text{ if $\ell = n$}\\
&\\
 \int_{U_{G_n} \backslash G_n} \int_{R^{\ell, n}} \int_{X_n}
W({}^{\omega_{\ell-n, n}}(rxg)) f_{s}(g)  \omega_{\psi^{-1}, \Upsilon^{-1}}(g) \phi(x) \, dy \, dg & \text{ if $\ell > n$}
\end{array}
\right.
\end{equation}
where $\xi_\ell = (0, \dots, 0, 1) \in E^{\ell}$, $X_k = \{ (x, 0, 0) \in H_k : x \in E^k \}$, $Y_k = \{ (0, y, 0) \in H_k : y \in E^k \}$,
\[
R_{\ell, n} = \left\{\begin{pmatrix}I_{n-\ell-1} &0&r_1&0&r_2&r_3\\ &1&&&&r_2^\prime \\ &&I_\ell&&&0\\ &&&I_\ell&&r_1^\prime \\ &&&&1&0\\ &&&&&I_{n-\ell-1} \end{pmatrix} 
\in G_n \right\} 
\]
\[
R^{\ell, n} = \left\{ \begin{pmatrix}I_{\ell-n-1}&0&r_1&&&\\ &1&&&&\\  &&I_n&&&\\  &&&I_n&&r_1^\prime\\  &&&&1&0\\ &&&&&I_{\ell-n-1} \end{pmatrix}  \in G_\ell \right\}
\]
\[
\omega_{\ell-n, n} = \begin{pmatrix} &I_{\ell-n}&&\\ I_n&&&\\ &&&I_n\\ &&I_{\ell-n}&\end{pmatrix} \in G_\ell,
\qquad
\beta_{\ell,n} = \begin{pmatrix} &I_\ell&&\\ &&&I_{n-\ell}\\ -I_{n-\ell}&&&\\ &&I_\ell& \end{pmatrix}.
\]
When $\ell < n$, in \cite[(2.16)]{BAS}, they define a local zeta integral by
\[
 \int_{U_{G_\ell} \backslash G_\ell} \int_{X_{\ell, n}} \int_{R_{\ell, n}}
W(g) f_{s}(\beta_{\ell,n} rh g)  \omega_{\psi^{-1}, \Upsilon^{-1}}(h g) \phi(\xi_\ell) \, dr \, dh \, dg
\]
where
\[
X_{\ell, n} = \left\{y(x, c):= \begin{pmatrix}I_{n-\ell}&x&0&c\\ &I_\ell&0&0\\ &&I_\ell&-J_\ell {}^{t}\bar{x}\\ &&&I_{n-\ell} \end{pmatrix} \in G_{n}:
x \in E^\ell, c \in E \right\}.
\]
Since we have 
\[
Y_\ell \backslash H_{\ell} \simeq X_{\ell, n},
\]
this integral is equal to our local zeta integral.

Recall that our local zeta integral converges absolutely for $\mathrm{Re}(s) \gg 0$ and has a meromorphic continuation
to $\mathbb{C}$ when $\ell \geq n$ by \cite[Proposition~6.4]{BAS2}.
Moreover, when $\ell < n$, the absolute convergence is proved in \cite[Proposition~3.1, Proposition~3.2]{BAS}, and 
as noted in the beginning of \cite[Section~5]{BAS}, we have the meromorphic continuation when $F$ is non-archimedean.
In the case where $\ell <n$ and $F$ is archimedean, the meromorphic continuation follows from the proof of 
the multiplicativity of $\gamma$-factors with the argument in \cite[Section~6]{S2} (See Remark~\ref{conv rem}).
Moreover, in the archimedean case,  the local zeta integral is continuous as a trilinear form.
Indeed, when $s$ is in the domain of absolute convergence, a similar continuity was proved by Soudry~\cite{S2} in the case of odd orthogonal groups,
and our case is proved in a similar way. 
For any $s$, the case of $l \geq n$ is proved in a similar way as \cite{S2} (see also \cite[Lemma~3.5]{GRS98}). 
For $\ell < n$, as in \cite{S2}, we may show the continuity using the multiplicativity proved in this paper.

Put  $m_0=\mathrm{max}(\ell, n)$ and $m= \mathrm{min}(\ell, n)$.
Defne
\[
U_{m_0} := \left\{ \begin{pmatrix}z&v_1&v_2\\ &I_{2(m+1)}&v_1^\prime\\ &&z^\ast \end{pmatrix} \in G_{m_0}, z \in U_{GL_{m_0-m-1}} \right\}
\]
and we write by $\psi_{m_0}$ the trivial extension of $\psi_{m_0-m-1}(u)$ to $U_{m_0}$, namely we define
\begin{equation}
\label{def Ul char}
\psi_{m_0}(u)= \psi \left(\sum_{i=1}^{m_0-m-1} u_{i, i+1} \right)
\end{equation}
when $E$ is a quadratic extension and 
\begin{equation}
\label{def Ul char GL}
\psi_{m_0}(u) = \psi_F \left(\left( \sum_{i=1}^{m_0-m-1} u_{i, i+1} \right)-u_{2m_0-1-i, 2m_0-i} \right)
\end{equation}
when $E=F \times F$.
Then $G_m$ is the stabilizer of the character $\psi_{m_0}$ of $U_{m_0}$, and thus
these integrals satisfy the following equivalence properties
\begin{equation}
\label{3.6 inert}
\left\{
\begin{array}{ll}
\mathcal{L}(g \cdot W,  (uhg) \cdot f_{s}, (hg) \cdot \phi)
=\psi_{m_0}^{-1}(u) \mathcal{L}(W, f_{s}, \phi)
& \text{ $\ell <n$,}\\
&\\
\mathcal{L}(g \cdot W, g \cdot  f_{s}, g \cdot \phi)
=\mathcal{L}(W, f_{s}, \phi)
& \text{ $\ell =n$,}\\
&\\
\mathcal{L}((uhg) \cdot W,  g \cdot f_{s},(hg) \cdot \phi)
=\psi_{m_0}^{-1}(u) \mathcal{L}(W, f_{s}, \phi)
& \text{ $\ell > n$,}\\
\end{array}
\right.
\end{equation}
where  $g \in G_m$, $h \in H_{m_0}$ and $u \in U_{m_0}$.
%The meromorphic continuation of this integral is proved in the same argument as 
%the proof of \cite[Lemma~3.5]{GRS98} when 

The intertwining operator $M(\tau, s) : V_{Q_n}^{G_n} \left( \mathcal{W}(\tau, \psi), s \right) \rightarrow 
V_{Q_n}^{G_n} \left( \mathcal{W}(\tau^\ast, \psi), 1-s \right)$ is defined by 
\[
M(\tau, s) f_s(h,a) = \int_{U_n} f_s(\omega_n^{-1} u h, d_n a^\ast)
\]
where $\omega_n = \left( \begin{smallmatrix}&I_n\\ -I_n& \end{smallmatrix} \right)$
and $d_i = \mathrm{diag}(-1, 1, \dots, (-1)^i) \in \mathrm{GL}_i(F)$.
Let us denote by $M^\ast(\tau, s)$ the intertwining operator normalized by Shahidi's $\gamma$-factor \cite{Sh2},
which is defined to satisfy the following identity
\begin{multline}
\label{int op}
\int_{U_n^n} f_{s}(d_n \omega_n u, I_n) \psi^{-1} \left( \frac{1}{2}u_{n, n+1} \right) \, du
\\
= \int_{U_n^n} M^\ast(\tau, s) f_{s}(d_n \omega_n u, I_n) \psi^{-1} \left( \frac{1}{2}u_{n, n+1} \right) \, du
\end{multline}
when $E$ is a quadratic extension and 
\begin{multline}
\label{3.7 split}
\int_{U_n} f_s(d_{2n} \omega_n u, I_n, I_n) \psi_F^{-1} \left(u_{n, n+1} \right)\, du
\\
= \int_{U_n} M^\ast(\tau, s) f_s(d_{2n} \omega_n^\circ u, I_n, I_n) \psi_F^{-1} \left(u_{n, n+1} \right)\, du
\end{multline}
when $E = F \times F$.
Here,  the measure $du$ is defined as a product measure, where
each root subgroup is isomorphic to $E$ or $F$, and 
we take self-dual measure on them with respect to $\psi$.
Then we define 
\[
\mathcal{L}^\ast(W, f_{s}, \phi) = \mathcal{L}(W, M^\ast(\tau, s)f_{s}, \phi).
\]
%%%%%%%%%%%%%%%%%%%%%%%%%%%%%%%%%%%%%%%%%%%%%%%
%
%
%
%
%
%
%
%
%
%
%
%
%
%
%%%%%%%%%%%%%%%%%%%%%%%%%%%%%%%%%%%%%%%%%%%%%%%
%%%%%%%%%%%%%%%%%%%%%%%%%%%%%%%%%%%%%%%%%%%%%%%
%
%
%
%
%
%
%
%
%
%
%
%
%
%
%%%%%%%%%%%%%%%%%%%%%%%%%%%%%%%%%%%%%%%%%%%%%%%
\section{Main Theorem}
\label{section:Main Theorem}
Suppose that $\pi$ and $\tau$ are representations of $G_\ell$ and $\mathrm{GL}_k(E)$ with central characters $\omega_\pi$ and $\omega_\tau$, respectively.
Moreover, when $F$ is archimedean, we assume that $\pi$ and $\tau$ are irreducible.

The space of trilinear forms satisfying \eqref{3.6 inert} is one-dimensional 
outside of a discrete subset of $s$.
If the representations are irreducible, this result follows from \cite[Corollary~16.3]{GGP} when $F$ is non-archimedean 
and \cite[Theorem~A]{LiSu} when $F$ is archimedean.
In the case of non-archimedean local fields, 
one-dimensionality holds also for the more general class of representations of Whittaker type.
This is proved in \cite[Theorem~5.3]{BAS} when $\ell < n$ and Theorem~\ref{unique} in the appendix when $\ell \geq n$.

The meromorphic continuations of the integrals $\mathcal{L}(W, f_s, \phi)$ 
and $\mathcal{L}^\ast(W,f_s, \phi)$ regarded as trilinear forms, satisfy the same equivariance properties.
Thus they are proportional. Furthermore, the form $\mathcal{L}(W, f_s, \phi)$  is not identically zero by \cite[Proposition~6.5]{BAS2} and  \cite[Proposition~4.1]{BAS}.
Therefore there exists a proportionality factor $\gamma(s, \pi \times \tau, \Upsilon, \psi)$ satisfying
%Since $\mathcal{L}(W)$ are trilinear forms satisfying the same equivalence properties,
%and thus these are proportional for all but discrete subset of $s$.
%Then we have a proportionality factor $\gamma(s, \pi \times \tau, \Upsilon, \psi)$ (resp. $\gamma(s, \pi \times (\tau_1, \tau_2), \Upsilon, \psi)$) 
%such that 
\begin{equation}
\label{def of our gamma}
\gamma(s, \pi \times \tau, \Upsilon, \psi)\mathcal{L}(W, f_{s}, \phi) = \mathcal{L}^\ast(W, f_{s}, \phi) 
\end{equation}
%\[
%\left(\text{resp.  } \gamma(s, \pi \times (\tau_1, \tau_2), \mu, \psi)\mathcal{L}(W, f_{s}, \phi) = \mathcal{L}^\ast(W, f_{s}, \phi)  \right)
%\]
for all $W, f_s$ and $\phi$.
Moreover, we define normalized $\gamma$-factors by
\[
\Gamma(s, \pi \times \tau, \Upsilon, \psi) :=\omega_\pi(-1)^n \omega_\tau(-1)^\ell \gamma(s, \pi \times \tau, \Upsilon, \psi)
\]
In the split case, let us we write $\Upsilon = \mu \otimes \mu^{-1}$ with a character $\mu$ of $F^\times$.
When we specify $\mu$, we often write 
\[
\gamma(s, \pi \times \tau, \Upsilon, \psi)
=\gamma(s, \pi \times \tau, \mu, \psi)
\]
and
\[
\Gamma(s, \pi \times \tau, \Upsilon, \psi)=\Gamma(s, \pi \times \tau, \mu, \psi).
\]
%and
%\[
%\Gamma(s, \pi \times (\tau_1, \tau_2), \mu, \psi) :=\omega_\pi(-1)^n \omega_{\tau_1}(-1)^\ell \omega_{\tau_2}(-1)^\ell \gamma(s, \pi \times (\tau_1, \tau_2), \mu, \psi).
%\]
Now, let us state our main result of this paper.
\begin{theorem}
\label{main thm}
Let $\pi$ and $\tau$ be irreducible generic representations
of $G_\ell$ and $\mathrm{GL}_k(E)$, respectively. 
Then the $\gamma$-factors $\Gamma(s, \pi \times \tau, \Upsilon, \psi)$ satisfy the following properties.
\begin{enumerate}
\item \label{e:1} $(\text{Unramified twisting})$ $\Gamma(s, \pi \times \tau |\det|^{s_0}, \Upsilon, \psi) = \Gamma(s+s_0, \pi \times \tau, \Upsilon, \psi) $.
% and $ \Gamma(s, \pi \times (\tau_1 |\det|^{s_0}, \tau_2 |\det|^{s_0}), \mu, \psi) = \Gamma(s+s_0, \pi \times (\tau_1, \tau_2), \mu, \psi) )$.
%%%%%%%%%%
\item \label{e:2}
\begin{enumerate}
\item $(\text{Multiplicativity : non-split case})$ 
Let $\sigma_1 \otimes \cdots \otimes \sigma_m \otimes \pi^\prime$ be an irreducible generic representation
of $\mathrm{GL}_{r_1}(E) \times \cdots \times \mathrm{GL}_{r_m}(E) \times G_{(l-r)}$,
$r = r_1 + \cdots + r_m$ $(0 \leq r \leq l)$ and let $\tau_1 \otimes \cdots \otimes \tau_a$
be an irreducible generic representation of $\mathrm{GL}_{n_1}(E) \times \cdots \times \mathrm{GL}_{n_a}(E)$ $( 1 \leq a \leq n)$.
Assume that $\pi$ (resp. $\tau$) is the irreducible generic quotient of a representation parabolically induced from 
$\sigma_1 \otimes \cdots \otimes \sigma_m \otimes \pi^\prime$ (resp. $\tau_1 \otimes \cdots \otimes \tau_a$).
Then 
\begin{align*}
&\Gamma(s, \pi \times \tau, \Upsilon, \psi) \\
= &\Gamma(s, \pi^\prime \times \tau, \Upsilon, \psi)
 \prod_{i=1}^m \gamma(s, (\sigma_i \otimes \Upsilon) \times \tau, \psi)\gamma(s, (\sigma_i \otimes \Upsilon)^\ast \times \tau, \psi)
\end{align*}
and
\[
\Gamma(s, \pi \times \tau, \Upsilon, \psi) = \prod_{i=1}^a \Gamma(s, \pi \times \tau_i, \Upsilon, \psi) .
\]
Here, $ \gamma(s, \sigma_i \times \tau, \psi)$ is the $\mathrm{GL}_{r_i} \times \mathrm{GL}_n$ 
$\gamma$-factor of Jacquet, Piatetski-Shapiro and Shalika~\cite{JPSS} in the non-archimedean
case and Jacquet and Shalika \cite{JS2} in the archimedean case.
\item $(\text{Multiplicativity : split case})$ 
Let $\sigma_{1,1} \otimes \cdots \otimes \sigma_{1, m} \otimes \pi^\prime \otimes \sigma_{2,1} \otimes \cdots \otimes \sigma_{2,m}$ 
be an irreducible generic representation of $\mathrm{GL}_{r_1}(F) \times \cdots \times \mathrm{GL}_{r_m}(F) \times G_{\ell-r}
\times \mathrm{GL}_{r_1}(F) \times \cdots \times \mathrm{GL}_{r_m}(F)$,
$r = r_1 + \cdots + r_m$ $(0 \leq r \leq l)$ and let $\tau_{1,i} \otimes \cdots \otimes \tau_{a,i}$
be an irreducible generic representation of $\mathrm{GL}_{n_1}(F) \times \cdots \times \mathrm{GL}_{n_a}(F)$ $( 1 \leq a \leq n)$.
Assume that $\pi$ (resp. $\tau_{i}$) is the irreducible generic quotient of a representation parabolically induced from 
$\sigma_{1,1} \otimes \cdots \otimes \sigma_{1, m} \otimes \pi^\prime \otimes \sigma_{2,1}^\ast \otimes \cdots \otimes \sigma_{2,m}^\ast$ (resp. $\tau_{1,i} \otimes \cdots \otimes \tau_{a,i}$).
Then 
\begin{multline*}
\Gamma(s, \pi \times \tau, \Upsilon, \psi) 
\\
= \Gamma(s, \pi^\prime \times \tau, \Upsilon, \psi)
\prod_{i=1}^2 \prod_{j=1}^m \gamma(s, (\sigma_{i,j} \otimes \mu) \times \tau_1, \psi)  \gamma(s, (\sigma_{i, j} \otimes \mu) \times \tau_2, \psi)
\\
\times
 \gamma(s, (\sigma_{i, j} \otimes \mu)^\ast \times \tau_1, \psi) \gamma(s, (\sigma_{i, j} \otimes \mu)^\ast \times \tau_2, \psi)
\end{multline*}
and
\[
\Gamma(s, \pi \times (\tau_1 \otimes \tau_2), \Upsilon, \psi) = \prod_{i=1}^a \Gamma(s, \pi \times (\tau_{i,1} \otimes\tau_{i, 2}), \Upsilon, \psi) 
\]
with $\Upsilon = \mu \otimes \mu^{-1}$.
\end{enumerate}
%%%%%%%%%%
\item \label{e:3} $(\text{Unramified case})$
When all data are unramified,
\[
\Gamma(s, \pi \times \tau, \Upsilon, \psi) = \frac{L(1-s, \pi^\vee \times (\tau^\vee \otimes \Upsilon))}{L(s, \pi \times (\tau \otimes \Upsilon^{-1})) } 
\]
in the non-split case, and 
\[
\Gamma(s, \pi \times (\tau_1 \otimes \tau_2), \mu, \psi) = \frac{L(1-s, \pi^\vee \times (\tau_1^\vee \otimes \mu)) L(1-s, \pi \times (\tau_2^\vee \otimes \mu^{-1}))}
{L(s, \pi \times (\tau_1 \otimes \mu^{-1}))L(s, \pi^\vee \times (\tau_2 \otimes \mu)) } 
\]
in the split case. See \cite{BAS2, BAS} for the definition of local $L$-factors.
\item \label{e:4} $(\text{Functional equation})$  $\Gamma(s, \pi \times \tau, \Upsilon, \psi) \Gamma(1-s, \pi^\vee \times \tau^\vee, \Upsilon^{-1}, \psi^{-1}) =1$.
%and  $\Gamma(s, \pi \times (\tau_1, \tau_2), \mu, \psi) \Gamma(1-s, \pi^\vee \times (\tau_1^\vee, \tau_2^\vee), \mu^{-1}, \psi^{-1}) =1$.
%\item  $(\text{Self-duality})$ $\Gamma(s, \pi^\vee \times \tau, \Upsilon^{-1}, \psi^{-1}) = \Gamma(s, \pi \times \tau, \Upsilon, \psi)$
%and  $\Gamma(s, \pi^\vee \times \tau, \psi^{-1}) = \Gamma(s, \pi \times \tau, \psi)$
\item  \label{e:5} $(\text{Dependence on $\psi$})$  For any $b = \bar{a} a$ with $a \in E^\times$, let $\psi_b$ be the character 
given by $\psi_b(x) = \psi(bx)$.
Then 
\[
\Gamma(s, \pi \times \tau, \Upsilon, \psi_b) =\omega_\tau(b)^{2\ell}  |b|^{2n\ell \left(s-\frac{1}{2} \right)}\Gamma(s, \pi \times \tau, \Upsilon, \psi).
\]
%For any $b \in F^\times$, we have 
%\[
%\Gamma(s, \pi \times (\tau_1, \tau_2), \mu, \psi_b) =\omega_{\tau_1} \omega_{\tau_2}(b)^{2\ell}  |b|^{2n\ell \left(s-\frac{1}{2} \right)} \Gamma(s, \pi \times (\tau_1, \tau_2), \mu, \psi).
%\]
\item \label{e:6} $(\text{Minimal case})$ For $\ell=0$, $\Gamma(s, \pi \times \tau, \Upsilon, \psi) = 1$.
%$\Gamma(s, \pi \times (\tau_1, \tau_2), \mu, \psi) = 1$.
\item \label{e:7} $(\text{Archimedean property})$
For any archimedean field $F$,
\[
\Gamma(s, \pi \times \tau, \Upsilon, \psi) = 
\gamma^{\rm Artin}(s, \pi \times (\tau \otimes \Upsilon^{-1}), \psi)
\]
in the non-split case and 
\[
\Gamma(s, \pi \times (\tau_1\otimes \tau_2), \mu, \psi) = 
\gamma^{\rm Artin}(s, \pi \times (\tau_1 \otimes \mu^{-1}), \psi) \gamma^{\rm Artin}(s, \pi^\vee \times (\tau_2 \otimes \mu), \psi)
\]
in the split case, where $\gamma^{\rm Artin}(\cdots)$ denotes the Artin $\gamma$-factor defined via local Langlands correspondence.
\item \label{e:8} $(\text{Global property})$ Let $L \slash K$ be a quadratic extension of number fields
with the adeles $\mathbb{A}_L \slash \mathbb{A}_K$.
Let $\Pi$ and $\Xi$ be irreducible cuspidal automorphic globally generic representations
of $\mathrm{U}_{2k}(\mathbb{A}_K)$ and $\mathrm{GL}_{n}(\mathbb{A}_L)$.
Let $\xi$ be a character of $\mathbb{A}_L^\times \slash L^\times$ whose restriction to $\mathbb{A}^\times_K$
is the quadratic character corresponding to $L \slash K$.
Let $S$ be a finite set of places of $K$ such that for $v \not \in S$, all data are unramified.
Then
\[
L^{S}(s, \Pi \times (\Xi \otimes \xi^{-1})) 
= \left( \prod_{v \in S}\Gamma(s, \Pi_v \times \Xi_v, \xi_v, \psi_v) \right)
L^{S}(1-s, \Pi^\vee \times (\Xi^\vee \otimes \xi)) .
\]
Here, $L^{S}(s, \Pi \times (\Xi \otimes \xi^{-1}))$ is the partial $L$-function with respect to $S$.
\end{enumerate}
\end{theorem}
Before proceeding with a proof of Theorem~\ref{main thm}, we would like to give a proof of Corollary~\ref{1}.
We note the following proof of Corollary~\ref{1} in the split case is slightly different from the usual argument. 
Indeed, the property \ref{e:2} is not full multiplicativity but partial multiplicativity, and thus we will modify the argument as follows.
\renewcommand*{\proofname}{Proof of Corollary~\ref{1}}
\begin{proof}
When $F$ is archimedean, Corollary~\ref{1} follows from properties \ref{e:2} and \ref{e:6} and \cite{Sh85}.
Hereafter, we suppose that $F$ is non-archimedean. First, let us consider the case where $E$ is a quadratic extension. 
The proof of this case is well-known, however for the convenience of the reader,
let us briefly recall it. Indeed, by the multiplicativity of $\gamma$-factors (i.e., property \ref{e:2}), 
it suffices to show the required equality when $\pi$ and $\tau$ are irreducible unitary supercuspidal representations.
In this case, by \cite[Proposition~5.1]{Sh2}, there exists
\begin{enumerate}
\item[$\bullet$] a totally complex number field $K$,
\item[$\bullet$] a quadratic extension field $L$ of $K$,
\item[$\bullet$] a finite place $v_0$ of $K$,
\item[$\bullet$] a generic character $\chi$ of $U_{G_\ell}(\mathbb{A}_K)$,
\item[$\bullet$] an irreducible $\chi$-generic cuspidal automorphic representation $\Pi = \otimes \Pi_v$ of $G_\ell(\mathbb{A}_K)$,
\end{enumerate}
where $\mathbb{A}_K$ is the ring of adeles of $K$, such that
\begin{enumerate}
\item[$\bullet$] $K_{v_0} \simeq F$ and $L_{v_0} \simeq E$,
\item[$\bullet$] $\chi_{v_0} = \psi_{U_{G_\ell}}$,
\item[$\bullet$] $\Pi_{v_0} \simeq \pi$,
\item[$\bullet$] $\Pi_v$ is unramified for every finite place $v \ne v_0$.
\end{enumerate}
Similarly, we may take an irreducible cuspidal automorphic representation $\Sigma$ (resp. $\xi$) of $\mathrm{GL}_n(\mathbb{A}_L)$
such that $\Sigma_{v_0} \simeq \tau$ (resp. $\xi_{v_0} = \Upsilon$) and $\Sigma_{v}$ (resp $\xi_v$) is unramified at every finite place $v \ne v_0$.
Then by properties \ref{e:3} and \ref{e:8}, we have 
\[
L^{S}(s, \Pi \times (\Sigma \otimes \xi^{-1})) 
= \left( \prod_{v \in S}\Gamma(s, \Pi_v \times \Sigma_v, \xi_v, \chi_v) \right)
L^{S}(1-s, \Pi^\vee \times (\Sigma^\vee \otimes \xi)) 
\]
with a finite set $S$ of places such that for $v \not \in S$, all data are unramified.
For a finite place $w  \in S \setminus \{ v_0\} $, we find that our $\gamma$-factor coincides with Shahidi's $\gamma$-factor
by the property \ref{e:2} since $\Pi_w$ and $\Sigma_w$ are principal series representations.
Moreover, as we noted, our $\gamma$-factor coincides with Shahidi's $\gamma$-factor at the place $\infty$.

By \cite{Sh2}, the above functional equation holds when replacing the above product of $\gamma$-factors by 
$\prod_{v \in S} \gamma^{Sh}(s, (\Pi_v \otimes \xi_v) \times \Sigma_v, \chi_v)$, and thus these products of $\gamma$-factors
coincide. Therefore, we obtain the required equality of our $\gamma$-factors and Shaidi's $\gamma$-factors at the place $v_0$.
\\

Let us consider the split case.
As above, we may suppose that $\tau$ is an irreducible unitary supercuspidal representation.
By Zelevinsky~\cite{Zel}, $\pi$ is an irreducible subquotient of a representation parabolically induced from 
$\sigma_1 \otimes (\sigma_1|\cdot|) \otimes \cdots \otimes \sigma_1|\cdot|^{b_1-1} \otimes \cdots \otimes 
\sigma_t \otimes (\sigma_t|\cdot|) \otimes \cdots \otimes \sigma_t|\cdot|^{b_t-1}$
where $b_i \in \mathbb{N}$ and $\sigma_i$ are irreducible supercuspidal representations of $\mathrm{GL}_{n_i}$ such that $b_1n_1+\cdots b_tn_t=2\ell$.
We note that in Theorem~\ref{main thm},  the multiplicativity is stated only for a generic quotient of an induced representation, however
it holds also for a generic subquotient since an induced representation has unique generic subquotient and Whittaker functional on 
an induced representation gives a Whittaker functional on this unique generic subquotient.
Then by the multiplicativity, we may suppose that $b_i=1$ and that $n_i \ne n_j$ for $i \ne j$.
In this case, the induced representation is irreducible and thus it is equal to $\pi$ (\cite{BZ77}).

Let us take $a_i \in \mathbb{R}$ such that  $\sigma_i^0 := \sigma_i|\cdot|^{a_i}$ is unitary.
For $z_i \in \mathbb{C}$, we put $\sigma_i^0[z_i] = \sigma_i^0|\cdot|^{z_i}$. When $z_i \in \sqrt{-1} \mathbb{R}$, $\sigma_i^0[z_i]$ is unitary.
Then as in the non-split case, for $z_i \in \sqrt{-1} \mathbb{R}$,
there exists irreducible cuspidal generic automorphic representations 
$\Sigma_i[z_i]$ of quasi-split unitary group $\mathrm{U}_{n_i}(\mathbb{A}_K)$ of degree $n_i$
such that $\Sigma_{i}[z_i]_{v_0} = \sigma$ and $\Sigma_{i}[z_i]_{v}$ is unramified at every finite place $v \ne v_0$.

By the global descent method \cite[Theorem~3.1]{GRS2}, there is an irreducible cuspidal generic automorphic representation $\Sigma[z_1, \dots, z_t]$ 
of $G_\ell(\mathbb{A}_K)$ such that its base change lift to $\mathrm{GL}_{2\ell}(\mathbb{A}_L)$ is 
$\mathrm{BC}(\Sigma_1[z_1]) \boxplus \cdots \boxplus \mathrm{BC}(\Sigma_t[z_t])$ with the base change lift $\mathrm{BC}(\Sigma_i[z_i])$ of $\Sigma_i[z_i]$.
Then $\Sigma[z_1, \dots, z_t]_{v_0}$ is isomorphic to the parabolic induction of
$\sigma_{1}^0[z_1] \otimes \cdots \otimes \sigma_{t}^0[z_t]$ which is irreducible by \cite{BZ77}. 
Moreover, for these representations, 
the required equality holds as in the non-split case.

For any $z_i \in \mathbb{C}$, let us consider the parabolic induction of
$\sigma_{1}^0|\cdot|^{z_1} \otimes \cdots \otimes \sigma_{t}^0|\cdot|^{z_t}$ which is irreducible by  \cite{BZ77} and thus generic.
We denote this representation by $\pi(z_1, \dots, z_t)$.
Then from \cite[Proposition~6.4]{BAS2} and \cite[Proposition~3.2]{BAS} and their proofs, there exists $k_{z_1, \dots, z_t}>0$ such that 
our local zeta integral for $\pi(z_1, \dots, z_t) \times \tau$ converges absolutely for $\mathrm{Re}(s) >k_{z_1, \dots, z_t}$ 
and uniformly for $s$ in a compact set in $\mathrm{Re}(s) >k_{z_1, \dots, z_t}$. 
Moreover, from their proofs, we may take $k_{z_i, \dots, z_t}$ uniformly when $z_i$ is in a given compact set.
Then we see that $\Gamma(s, \pi(z_1, \dots, z_t) \times \tau, \Upsilon, \psi)$
is holomorphic on the domain $\mathrm{Re}(s) \ll 0, |\mathrm{Re}(z_i)| < |a_i|+1$.
Indeed, our $\gamma$-factor is defined as a quotient of our local integrals from the definition \eqref{def of our gamma}. 
We may choose data so that our local zeta integral is equal to 1 by \cite[Proposition~4.9, Proposition~6.5]{BAS2} and  \cite[Proposition~4.1]{BAS},
and thus the left hand side of \eqref{def of our gamma} is equal to our $\gamma$-factor.
Then from the above convergences, the right-hand side of \eqref{def of our gamma} becomes holomorphic and thus $\Gamma(s, \pi(z_1, \dots, z_t) \times \tau, \Upsilon, \psi)$
also becomes holomorphic as required.

Similarly, we see that 
$\gamma^{Sh}(s, (\pi(z_1, \dots, z_t) \otimes \mu) \times \tau_1,\psi) \gamma^{Sh}(s, (\pi(z_1, \dots, z_t)  \otimes \mu) \times \tau_2, \psi)$ is holomorphic. 
As we showed above, when $z_i \in \sqrt{-1} \mathbb{R}$, these $\gamma$-factors coincide.
Hence, we have for any $z_i$ in the domain $\mathrm{Re}(s) \ll 0, |\mathrm{Re}(z_i)| < |a_i|+1$, 
\begin{multline*}
\Gamma(s, \pi(z_1, \dots, z_t) \times \tau, \Upsilon, \psi)
\\
=\gamma^{Sh}(s, (\pi(z_1, \dots, z_t) \otimes \mu) \times \tau_1,\psi) \gamma^{Sh}(s, (\pi(z_1, \dots, z_t)  \otimes \mu) \times \tau_2, \psi).
\end{multline*}
In particular, our required equality holds by specializing $z_i=a_i$.
\end{proof}
\renewcommand*{\proofname}{Proof}
In the next sections, we shall show only the property \ref{e:2}.
The property \ref{e:1} is immediate from the definition.
The property \ref{e:3} is proved in \cite{BAS2, BAS}.
The property \ref{e:4} is proved in the same argument as other classical groups or metaplectic groups (for example, see \cite[p.413]{Ka}).
The property \ref{e:6} is trivial.
The property \ref{e:7} follows from the property \ref{e:2} by the Casselman's subrepresentation theorem.
The property \ref{e:8} for unnormalized $\gamma$-factors $\gamma(s, \Pi_v \times \Xi_v, \xi_v, \psi_v)$ 
easily follows from the properties of the Eisenstein series in a standard way as \cite[Theorem~A (c)]{Gi}.
Moreover, since $\prod_v \omega_{\Pi_v}(-1) = \prod_{v} \omega_{\Xi_v}(-1)=1$, 
the property \ref{e:8} holds for our  normalized gamma factors.

%%%%%%%%%%%%%%%%%%%%%%%%%%%%%%%%%%%%%%%%%%%%%%%
%
%
%
%
%
%
%
%
%
%
%
%
%
%
%%%%%%%%%%%%%%%%%%%%%%%%%%%%%%%%%%%%%%%%%%%%%%%
Let us give a proof of the property \ref{e:5}.
Let us consider the non-split case. For all $\ell, n$, we set
\begin{align*}
t_a &= \mathrm{diag}(a^{\ell} \bar{a}^{\ell-1}, a^{\ell-1}\bar{a}^{\ell-2}, \dots, a, \bar{a}^{-1}, \dots, \bar{a}^{-l} a^{-\ell+1}) \in G_\ell, \\
y_a &= \mathrm{diag}(a^{\ell} \bar{a}^{\ell-1},a^{\ell-1}\bar{a}^{\ell-2}, \dots, a^{\ell-n+1} \bar{a}^{\ell-n}) \in \mathrm{GL}_{n}(E),\\
z_a &= \mathrm{diag}(a^{-1} I_{m+1}, \bar{a} I_{m+1}) \in G_{m+1}, \quad m= \mathrm{min}(\ell, n).
\end{align*}
For $W \in \mathcal{W}(\pi, \psi^{-1})$, the functions $g \mapsto W(t_a g)$ span the space $\mathcal{W}(\pi, \psi_b^{-1})$.
Similarly, the functions $g \mapsto W(y_a g)$ $(W \in \mathcal{W}(\tau, \psi))$ span the space $\mathcal{W}(\tau, \psi_b)$.
Also, $\omega_{\psi_b, \Upsilon}(hg) = \omega_{\psi, \Upsilon}({}^{z_a}hg)$ $(h \in H_m, g \in G_n)$.
Let us denote by $\mathcal{L}_b(W, f_s, \phi)$  the integral adapting to $\psi_b$ up to measure constants 
which cancel out later in the computation we see that, when $\ell \leq n$,
\begin{multline*}
\mathcal{L}_b(W, f_s, \phi) 
\\
=\omega_{\tau}(b)^{2\ell-n} |b|^{(2n\ell-n^2) \left(s-\frac{1}{2} \right)}
\mathcal{L}(t_a \cdot W, ({}^{\gamma_{\ell, n}} y_a) \cdot f_s, \omega_{\psi, \Upsilon}(t_a \cdot \mathrm{diag}(a I_{\ell}, \bar{a}^{-1} I_\ell))\phi)
\end{multline*}
and when $\ell > n$, replace ${}^{\gamma_{\ell, n}} y_a$ by $y_a$ in the above identity.
%For a measure $dx$ of some group or quotient space, defined with respect to $\psi$, 
%denote by $d_bx$ the measure adapted to $\psi_b$.
%For the fixed measureon $F$, $d_bx=|b|^{1 \slash 2} dx$.
%In general the measures $dx$ and $d_bx$ are related by 
%a constant and we set $d_bx=c(b,dx)dx$.
Moreover, let us denote by $M_b(\tau, s)$ the intertwining operator of $V_{Q_n}^{G_n}(\mathcal{W}(\tau, \psi_b),s)$ 
with respect to the measure associated to $\psi_b$. Also, we write its normalized intertwining operator by $M_b^\ast(\tau, s)$.
Then by a direct computation as in \cite[5.1]{Ka}, we can show that 
\[
M_b^\ast(\tau, s) f_s(h, d)
=\delta_{Q_n}(y_b)^{\frac{1}{2}} \omega_{\tau}(b)^{2\ell} |\det y_a|^{s -\frac{1}{2}} |\det t_a^\prime|^{2s-1}
M^\ast(\tau, s)f_s( y_a h, d)
\]
where $t_a^\prime = \mathrm{diag}(a^{n}\bar{a}^{n-1}, a^{n-1} \bar{a}^{n-2}, \dots, a) \in \mathrm{GL}_n(E)$.
From these relations, it is easy to see that the required dependence follows in a similar argument as \cite[Section~5.1]{Ka}.

Let us consider the split case. For all $l, n$, we set
\begin{align*}
t_b &= \mathrm{diag}(b^{\ell-1}, b^{\ell-2}, \dots, 1, b^{-1}, \dots, b^{-\ell+1}, b^{-\ell}) \in \mathrm{GL}_{2\ell}(F), \\
y_{b, 1}& = \mathrm{diag}(b^{\ell-1},b^{\ell-2}, \dots, b^{-n+\ell})\\
y_{b,2} &= \mathrm{diag}(b^{-\ell+n-1}, b^{-\ell+n-2}, \dots, b^{-\ell})\\
y_b &= \mathrm{diag}(y_{b,1}, y_{b,2}) \in \mathrm{GL}_{2n}(F) \\
z_b &= \mathrm{diag}(I_{m+1}, b I_{m+1}) \in \mathrm{GL}_{2m+2}(F), \quad m= \mathrm{min}(\ell, n),\\
t_{b, 1}^\prime &=\mathrm{diag}(b^{n-1}, b^{n-2}, \dots, 1)\\
t_{b, 2}^\prime &=\mathrm{diag}(1, b^{-1}, \dots, b^{-n+1} ) \in \mathrm{GL}_n(F)\\
t_b^\prime &= \mathrm{diag}(t_{b,1}^\prime, t_{b, 2}^\prime) \in \mathrm{GL}_{2n}(F).
\end{align*}
For $W \in \mathcal{W}(\pi, \psi^{-1})$, the functions $g \mapsto W(t_b g)$ span the space $\mathcal{W}(\pi, \psi_b^{-1})$.
Similarly, the functions $g \mapsto W(y_{b, 1} g)$ $(W \in \mathcal{W}(\tau_1, \psi))$ 
(resp.  $g \mapsto W(y_{b, 2} g)$ $(W \in \mathcal{W}(\tau_2^{\ast}, \psi^{-1}))$ ) span the space $\mathcal{W}(\tau_1, \psi_b)$ (resp. $\mathcal{W}(\tau_2^{\ast}, \psi_b^{-1})$).
Then the integral $\mathcal{L}_b(W, f_s, \phi)$ defined via $\psi_b$ is equal to, up to measure constants,  when $\ell \leq n$
\begin{multline*}
\mathcal{L}_b(W, f_s, \phi) =\omega_{\tau_1}(b)^{2\ell-n}\omega_{\tau_2}(b)^{2\ell-n}  |b|^{(2n\ell-n^2) \left(2s-1 \right)}
\\
\cdot \mathcal{L}(t_b \cdot W,  ({}^{\gamma_{\ell, n}} y_b) \cdot f_s, \omega_{\psi, \mu}(t_b \cdot \mathrm{diag}(I_{\ell},  b^{-1}I_\ell))\phi)
\end{multline*}
and when $\ell > n$, replace ${}^{\gamma_{\ell, n}} y_b$ by $y_b$ in the above identity.
Moreover, we have
\begin{multline*}
M_b^\ast((\tau_1 \otimes \tau_2), s) f_s(h, d)
=\omega_{\tau_1}(b)^{2\ell}\omega_{\tau_2}(b)^{2\ell} |\det y_{b,1} y_{b,2}^{-1} |^{s -\frac{1}{2}} 
\\ \cdot |\det t_{b,1}^\prime (t_{b,2}^\prime)^{-1} |^{2s-1}
M^\ast((\tau_1 \otimes \tau_2), s)f_s( y_b h, d)
\end{multline*}
and required identity easily follows as in non-split case.
%%%%%%%%%%%%%%%%%%%%%%%%%%%%%%%%%%%%%%%%%%%%%%%
%
%
%
%
%
%
%
%
%
%
%
%
%
%
%%%%%%%%%%%%%%%%%%%%%%%%%%%%%%%%%%%%%%%%%%%%%%%
\section{Proof of the property \ref{e:2}}
\label{proof:s}
Our aim in this section is to prove the multiplicativity property of $\gamma$-factors.
Indeed, it follows from the following propositions.
\begin{proposition}
\label{M1}
Let $\pi$ be an irreducible representation of $G_\ell$ and $\tau_1, \tau_2$
be representations of $\mathrm{GL}_{n_i}(E)$ with $n_1+n_2=n$.
Suppose that $\tau_i$ are irreducible if $F$ is archimedean and that 
$\tau_i$ are quotients of representations parabolically induced from irreducible supercuspidal representations
if $F$ is non-archimedean.
Let $\Sigma_{\tau_1, \tau_2} = \mathrm{Ind}_{P_{n_1, n_2}}^{\mathrm{GL}_n(E)}(\tau_1 \otimes \tau_2)$ if $F$ is non-archimedean 
and let $\Sigma_{\tau_1, \tau_2}$ be the unique generic irreducible subquotient of $\mathrm{Ind}_{P_{n_1, n_2}}^{\mathrm{GL}_n(E)}(\tau_1 \otimes \tau_2)$ if $F$ is archimedean.
Then we have
\[
\gamma \left(s, \pi \times \Sigma_{\tau_1, \tau_2}, \Upsilon, \psi \right)
= \gamma(s, \pi \times \tau_1, \Upsilon, \psi)  \gamma(s, \pi \times \tau_2, \Upsilon, \psi) .
\]
\end{proposition}
%\begin{proposition}
%\label{M1 split}
%Let $\pi$ be an irreducible representation of $\mathrm{GL}_{2\ell}(F)$ and $\tau_{i, j}$
%be representations of $\mathrm{GL}_{n_{i}}(F)$ with $n_1+n_2=n$.
%Put $\tau_i = \mathrm{Ind}_{P_{n_1, n_2}}^{\mathrm{GL}_{n_{i}}}(\tau_{i,1} \otimes \tau_{i,2})$.
%Suppose that $\tau_{i, j}$ are irreducible if $F$ is archimedean and that 
%$\tau_{i,j}$ are quotients of representations parabolically induced from supercuspidal representations
%f $F$ is non-archimedean.
%Then we have
%\begin{multline*}
%\gamma \left(s, \pi \times \left(\mathrm{Ind}_{P_{n_{1}, n_{2}}}^{\mathrm{GL}_{n}(F)}(\tau_{1,1} \otimes \tau_{1,2}) \otimes
%\mathrm{Ind}_{P_{n_{1}, n_{2}}}^{\mathrm{GL}_{n}(F)}(\tau_{2,1} \otimes \tau_{2,2}) \right), \mu, \psi \right)
%\\
%= \gamma(s, \pi \times (\tau_{1,1} \otimes \tau_{2,1}), \mu, \psi)  \gamma(s, \pi \times (\tau_{1,2} \otimes \tau_{2,2}), \mu, \psi) .
%\end{multline*}
%\end{proposition}
\begin{proposition}
\label{M2}
Let $\tau$ be an irreducible representation of $\mathrm{GL}_n(E)$.
Let $\pi^\prime$ be a representation of $G_{\ell_0}$ and $\sigma$
be a representation of $\mathrm{GL}_{r}(E)$ with $\ell_0+r=\ell$.
Suppose that $\pi^\prime$ and $\sigma$ are irreducible if $F$ is archimedean and that 
$\pi^\prime$ and $\sigma$ are quotients of representations parabolically induced from supercuspidal representations
if $F$ is non-archimedean.
Let $\Pi_{\sigma, \pi^\prime}=\mathrm{Ind}_{Q_r}^{G_\ell}(\sigma \otimes \pi^\prime)$ if $F$ is non-archimedean and 
let $\Pi_{\sigma, \pi^\prime}$ be the unique generic irreducible subquotient of $\mathrm{Ind}_{Q_r}^{G_\ell}(\sigma \otimes \pi^\prime)$ if $F$ is archimedean.
Then we have
\begin{multline*}
\gamma \left(s, \Pi_{\sigma, \pi^\prime} \times \tau, \Upsilon, \psi \right)
\\
=\omega_\sigma(-1)^n \omega_\tau(-1)^r \gamma(s, (\sigma \otimes \Upsilon)  \times \tau, \psi)   \gamma(s, \pi^\prime \times \tau, \Upsilon, \psi) 
\gamma(s, (\sigma \otimes \Upsilon)^\ast \times \tau, \psi) 
\end{multline*}
when $E$ is a quadratic extension, and 
\begin{multline*}
\gamma \left(s, \Pi_{\sigma, \pi^\prime} \times (\tau_1 \otimes \tau_2), \mu, \psi \right)\\
=\omega_{\sigma_1}(-1)^n \omega_{\sigma_2}(-1)^n  \omega_{\tau_1}(-1)^r \omega_{\tau_2}(-1)^r  \gamma(s, (\sigma_1 \otimes \mu) \times \tau_1, \psi) 
\gamma(s, (\sigma_1 \otimes \mu) \times \tau_2, \psi) \\
  \gamma(s, \pi^\prime \times (\tau_1 \otimes \tau_2), \mu, \psi) \gamma(s, (\sigma_2 \otimes \mu)^\ast \times \tau_1, \psi)  \gamma(s, (\sigma_2 \otimes \mu)^\ast \times \tau_2, \psi) 
\end{multline*}
%where $ \Pi_{\sigma_1, \sigma_2 \pi^\prime}$ is the unique generic irreducible subquotient of $\mathrm{Ind}_{Q_{r}}^{\mathrm{GL}_{2\ell}(F)}(\sigma_1 \otimes \pi^\prime \otimes \sigma_2)$
when $E = F \times F$ and $\sigma = \sigma_1 \otimes \sigma_2$.
\end{proposition}
%\begin{proposition}
%\label{M2 split}
%Let $\tau_1$ and $\tau_2$ be irreducible representations of $\mathrm{GL}_n(F)$.
%Let $\pi^\prime$ be a representation of $\mathrm{GL}_{2\ell_0}(F)$ and $\sigma_1, \sigma_2$
%be representations of $\mathrm{GL}_{r}(F)$ with $\ell_0+r=\ell$.
%Suppose that $\pi^\prime$ and $\sigma_i$ are irreducible if $F$ is archimedean and that 
%$\pi^\prime$ and $\sigma_i$ are quotients of representations parabolically induced from supercuspidal representations
%if $F$ is non-archimedean.
%Then we have
%\begin{multline*}
%\gamma \left(s, \mathrm{Ind}_{Q_{m}}^{\mathrm{GL}_{2\ell}(F)}(\sigma_1 \otimes \pi^\prime \otimes \sigma_2) \times (\tau_1 \otimes \tau_2), \mu, \psi \right)\\
%=\omega_{\sigma_1}(-1)^n \omega_{\sigma_2}(-1)^n  \omega_{\tau_1}(-1)^r \omega_{\tau_2}(-1)^r  \gamma(s, (\sigma_1 \otimes \mu) \times \tau_1, \psi) 
%\gamma(s, (\sigma_1 \otimes \mu) \times \tau_2, \psi) \\
%  \gamma(s, \pi^\prime \times (\tau_1 \otimes \tau_2), \mu, \psi) \gamma(s, (\sigma_2 \otimes \mu)^\ast \times \tau_1, \psi)  \gamma(s, (\sigma_2 \otimes \mu)^\ast \times \tau_2, \psi) .
%\end{multline*}
%\end{proposition}

We note that Proposition~\ref{M1} follows from Proposition~\ref{M2}
if $F$ is archimedean as in \cite{Sh2}.
Indeed, by Casselman's subrepresentation theorem, we may suppose that $\pi$ is subrepresentation
of $\mathrm{Ind}_{Q_{\ell}}^{G_\ell}(\pi_0^\prime)$ for some representation $\pi^\prime_0$ of $\mathrm{GL}_\ell(E)$. Here, we note that 
Levi component of $Q_\ell$ does not contain unitary group.
Hence, by Proposition~\ref{M2}, we may reduce Proposition~\ref{M1} to the multiplicativity 
of the $\gamma$-factors for $\mathrm{GL}_\ell \times \mathrm{GL}_{n_i}$, and we obtain required multiplicativity.

In the non-split case, Proposition~\ref{M1} and Proposition~\ref{M2} are proved by the same argument as the proof of 
\cite[Section~7]{Ka}, indeed these are proved by word for word. 
Hence, we omit its proof and we shall prove only these propositions in the split case.
See Remark~\ref{qs rem1} in Section~\ref{sect 4.1} and Remark~\ref{qs rem2} in Section~\ref{9.2} for main equalities in the non-split case.
Further, we note that any convergence of integrals studied in the following sections is proved in the standard argument similarly as \cite{S3} or \cite[Section~4.10]{Ka2}.
Hence, in the following sections, we omit them and only give formal computations (cf. Remark~\ref{verify conv} in Section~\ref{s:ell<n2}).

\subsection{Proof of Proposition~\ref{M1} in the split case}
\label{Proof of Proposition split}
From now on, we simply write $\mathrm{GL}_i = \mathrm{GL}_i(F)$.
As we remarked above, we give formal proof.
Indeed, in order to verify the following argument, we need to twist $\tau_{i,j}$ by a parameter $\zeta \in \mathbb{C}$.
As we noted, in the following argument, we omit $\zeta$.

As we remarked in the previous section, we assume that $F$ is non-archimedean throughout this section.
Let $V_{Q_n}^{\mathrm{GL}_{2n}}(\mathcal{W}(\tau_1, \psi) \otimes \mathcal{W}(\tau_2^\ast, \psi^{-1}), s)$, $M((\tau_1 \otimes \tau_2), s)$
and $M^\ast((\tau_1 \otimes \tau_2), s) $ be 
as in Section~\ref{s:def local zeta}.
We put 
\[
m(b_1, b_2) = \begin{pmatrix} b_1&\\ &b_2 \end{pmatrix}
\]
for $b_1 \in \mathrm{GL}_{i}$ and $b_2 \in \mathrm{GL}_j$.
For $0 <r_0 \leq r$, let $Q_{r_0} < \mathrm{GL}_{2r}$ be the standard maximal parabolic subgroup whose Levi part is 
\[
\left\{\begin{pmatrix}g_1&&\\ &h&\\ &&g_2 \end{pmatrix} : g_i \in \mathrm{GL}_{r_0}, \, h \in \mathrm{GL}_{2(r-r_0)} \right\}.
\]
Then its unipotent radical $U_{r_0}$ is 
\[
\left\{ \begin{pmatrix} I_{r_0}&\ast&\ast\\ &I_{2(r-r_0)}&\ast\\ &&I_{r_0} \end{pmatrix} \right\}.
\]
Let us define
\[
Z_{n_2, n_1} = \left\{ \begin{pmatrix} I_{n_2}&z_1&&\\ &I_{n_1}&&\\ &&I_{n_1}&z_2\\ &&&I_{n_2}\end{pmatrix}\right\}
\quad
\text{and}
\quad
\omega_{n_1, n_2} = \begin{pmatrix} &I_{n_1}&&\\ I_{n_2}&&&\\ &&&I_{n_2}\\ &&I_{n_1}&\end{pmatrix}.
\]
%Let $\tau = (\tau_1, \tau_2)$ be a pair of irreducible representations of $\mathrm{GL}_{n}$.
%Let $\varepsilon_{i,j}$ be an irreducible representation of $\mathrm{GL}_{n_j}$ for $1 \leq i, j \leq 2$.
%We write $\varepsilon_i = (\varepsilon_{i1}, \varepsilon_{i2})$.
%Suppose that 
%\[
%\tau_i = \mathrm{Ind}_{P_{n_1, n_2}}^{\mathrm{GL}_n}(\varepsilon_{i1} \otimes \varepsilon_{i2} )
%\]
%Then our aim is to prove the following identity
%\[
%\gamma(s, \pi \times \tau, \psi) = \gamma(s, \pi \times \varepsilon_1, \psi)\gamma(s, \pi \times \varepsilon_2, \psi).
%\]
For $s_1, s_2 \in \mathbb{C}$, we form the indued representation
\begin{multline*}
\mathrm{Ind}_{Q_{n_1}}^{\mathrm{GL}_{2n}} (\mathcal{W}(\tau_{11}, \psi)|\det|^{s_1-\frac{1}{2}} 
\\
\otimes 
\mathrm{Ind}_{Q_{n_2}}^{\mathrm{GL}_{2n_2}} \left(\mathcal{W}(\tau_{12}, \psi)|\det|^{s_2-\frac{1}{2}} \otimes \mathcal{W}(\tau_{22}^\ast, \psi^{-1})|\det|^{-s_2+\frac{1}{2}} \right)
\\
\otimes \mathcal{W}(\tau_{21}^\ast, \psi^{-1})|\det|^{-s_1+\frac{1}{2}}
)
\end{multline*}
which is realized on a space of functions on $\mathrm{GL}_{2n} \times \mathrm{GL}_{n_1} \times \mathrm{GL}_{n_1}
 \times \mathrm{GL}_{2n_2} \times \mathrm{GL}_{n_2} \times \mathrm{GL}_{n_2}$,
% \[
% f_{\tau_{i,j}, s}()
% \]
and we denote this space  by $V^{\mathrm{GL}_{2n}}_{Q_{n_1}}((\tau_{11}, \tau_{12}),(\tau_{21}^\ast, \tau_{22}^\ast), (s_1, s_2))$.
We may write this induced representation as
\[
 \mathrm{Ind}_{Q_n}^{\mathrm{GL}_{2n}}(V^{\mathrm{GL}_{n}}_{P_{n_1, n_2}}((\tau_{11}, \tau_{12}), (s_1, s_2)) \otimes
 V^{\mathrm{GL}_{n}}_{P_{n_2, n_1}}((\tau_{22}^\ast, \tau_{21}^\ast), (-s_2, -s_1))).
\]
Here, $V^{\mathrm{GL}_{n}}_{P_{n_1, n_2}}((\tau_{11}, \tau_{12}), (s_1, s_2))$ is the space of functions
 $f_{s_1, s_2}$ on $\mathrm{GL}_n \times \mathrm{GL}_{n_1} \times \mathrm{GL}_{n_2}$
such that 
\begin{multline*}
f_{s_1, s_2}( m(b_1, b_2)g, a_1, a_2) 
\\
= \delta_{P_{n_1, n_2}}(m(b_1, b_2))^\frac{1}{2} |\det b_1|^{s_1-\frac{1}{2}}
|\det b_2|^{s_2-\frac{1}{2}}f_{\tau_1, s_1, s_2}( g, a_1b_1, a_2b_2)
\\
=
|\det b_1|^{s_1+\frac{n_2-1}{2}}
|\det b_2|^{s_2-\frac{n_1+1}{2}}f_{\tau_1, s_1, s_2}( g, a_1b_1, a_2b_2).
\end{multline*}
For any $\varphi_{s} \in V^{\mathrm{GL}_{2n}}_{Q_{n_1}}((\tau_{11}, \tau_{12}),(\tau_{21}^\ast, \tau_{22}^\ast), (s, s))$, we define an element 
$f_{\varphi_s} \in V_{Q_n}^{\mathrm{GL}_{2n}}(\mathcal{W}(\tau_1, \psi) \otimes \mathcal{W}(\tau_2^\ast, \psi^{-1}), s)$ via
\begin{multline}
\label{8.2 split}
f_{\varphi_s}(g, a_1, a_2) 
= |\det a_1^{-1}a_2|^{s+\frac{n-1}{2}} 
\\
\times
\int_{Z_{n_2, n_1}} 
\varphi_s(\omega_{n_1, n_2}z m(a_1, a_2) g, I_{n_1}, I_{n_1}, I_{2n_2}, I_{n_2}, I_{n_2}) \psi^{-1}(z)\, dz
\end{multline}
where $\psi$ denotes the non-generate character of $U_{\mathrm{GL}_n}$.
We often simply write $\mathcal{L}(W, f_{\varphi_s}, \phi) = \mathcal{L}(W, \varphi_s, \phi)$.

For $\theta_s \in V^{\mathrm{GL}_{n}}_{P_{n_1, n_2}}((\tau_{11}, \tau_{22}^\ast), (s_1, s_2))$, 
we define the intertwining operator $M(\tau_{11} \otimes \tau_{22}^\ast, (s_1, s_2))$ by
\[
(M(\tau_{11} \otimes \tau_{22}^\ast, (s_1, s_2))\theta_s)(g, a_2, a_1)
=\int_{Z_{n_2, n_1}}\theta_s(\Omega_{n_1, n_2} zg, a_1, a_2) \, dz
\]
where $a_i \in \mathrm{GL}_{n_i}$, $g \in \mathrm{GL}_n$ and
\[
\Omega_{n_1, n_2} = \begin{pmatrix} &I_{n_1}\\ I_{n_2}&\end{pmatrix}.
\]
Then 
$M(\tau_{11} \otimes \tau_{22}^\ast, (s_1, s_2))\theta_s \in V^{\mathrm{GL}_{n}}_{P_{n_2, n_1}}((\tau_{22}^\ast, \tau_{11}), (s_2, s_1))$.
Further, the normalized intertwining operator $M^\ast(\tau_{11} \otimes \tau_{22}^\ast, (s_1, s_2))$ is defined by the functional equation
\begin{multline}
\label{3.7 sp}
\int_{Z_{n_2, n_1}}\theta_s(\Omega_{n_1, n_2} zg, I_{n_1}, I_{n_2}) \psi^{-1}(z) \, dz
\\
= \int_{Z_{n_2, n_1}} (M(\tau_{11} \otimes \tau_{22}^\ast, (s_1, s_2))^\ast \theta_s)(\Omega_{n_1, n_2} zg, I_{n_2}, I_{n_1}) \psi^{-1}(z) \, dz.
\end{multline}
In the above integral, we take the product measure of self dual Haar measure with respect to $\psi$
as the measure of $Z_{n_2, n_1}$. 
Similarly, we define intertwining operators $M(\tau_{12}^\ast \otimes \tau_{21}, (s_1, s_2))$
and $M^\ast(\tau_{12}^\ast \otimes \tau_{21}, (s_1, s_2))$,
whose images are in $V^{\mathrm{GL}_{n}}_{P_{n_1, n_2}}((\tau_{21}, \tau_{12}^\ast), (s_2, s_1))$.

Given $\varphi_s \in V^{\mathrm{GL}_{2n}}_{Q_{n_1}}((\tau_{11}, \tau_{12}),(\tau_{21}^\ast, \tau_{22}^\ast), (s, s))$,
we can repeatedly apply the intertwining operators and we get the following sections :
\begin{align*}
\varphi_s^\prime =& M_0^\ast((\tau_{12} \otimes \tau_{22}^\ast), s) \varphi_s 
\in V^{\mathrm{GL}_{2n}}_{Q_{n_1}}((\tau_{11}, \tau_{22}^\ast, \tau_{12}, \tau_{21}^\ast), (s, 1-s)),\\
 \varphi_s^{\prime \prime }=&  M^\ast((\tau_{11} \otimes \tau_{22}^\ast), (s, 1-s))  M^\ast((\tau_{21}^\ast \otimes \tau_{12}), (s, 1-s)) \varphi_s^{\prime \prime} 
 \\
 &\in V^{\mathrm{GL}_{2n}}_{Q_{n_1}}((\tau_{22}^\ast, \tau_{11}, \tau_{21}^\ast, \tau_{12}), (1-s, s)),\\
 \varphi_s^{\prime \prime \prime} =& M_0^\ast((\tau_{21} \otimes \tau_{11}^\ast), s) \varphi_s^{\prime \prime} \in 
 V^{\mathrm{GL}_{2n}}_{Q_{n_2}}((\tau_{12}^\ast, \tau_{11}^\ast), (\tau_{21}, \tau_{22}), (1-s, 1-s)).
\end{align*}
According to the multiplicativity of intertwining operators,
\[
M^\ast((\tau_1 \otimes \tau_2), s) f_{\varphi_s} = f_{ \varphi_s^{\prime \prime \prime}}.
\]
Then in a similar argument as the proof \cite[(7.1), Section~8]{Ka}, the following proposition yields the required multiplicativitiy.
\begin{proposition}
\label{8.5}
\[
\gamma(s, \pi \times (\tau_{12} \otimes \tau_{22}), \mu, \psi) \mathcal{L}(W, \varphi_s, \phi)
= \mathcal{L}(W, \varphi_s^\prime, \phi).
\]
\end{proposition}
We will prove this proposition in the following sections.
 %%%%%%%%%%%%%%%%%%%%%%%%%%%%%%%%%%%%%%%%%%%%%%%
%
%
%
%
%
%
%
%
%
%
%
%
%
%
%%%%%%%%%%%%%%%%%%%%%%%%%%%%%%%%%%%%%%%%%%%%%%%
\subsubsection{$\ell < n_2$}
\label{s:ell<n2}
Plugging \eqref{8.2 split} into $ \mathcal{L}(W, f_{\varphi_s}, \phi)$, it is equal to
\begin{multline}
\label{8.9}
 \int_{U_{\mathrm{GL}_{2\ell}} \backslash \mathrm{GL}_{2\ell}} \int_{Y_\ell \backslash H_\ell} \int_{R_{\ell, n}} \int_{Z_{n_2, n_1}}
W(g)\varphi_s(\omega_{n_1, n_2}z \beta_{\ell,n} rh g, I_{n_1}, I_{n_1}, I_{2n_2}, I_{n_2}, I_{n_2}) 
\\ \omega_{\psi^{-1}, \mu}(hg)\phi(\xi_\ell)
\psi^{-1}(z)\, dz
\, dr \, dh \, dg. 
\end{multline}

Write $R_{\ell, n} = R_{\ell, n_2} \ltimes (R_{\ell, n} \cap U_{n_1})$.
Here, recall that $U_{n_1}$ is the unipotent radical of the standard parabolic subgroup $Q_{n_1}$ of $\mathrm{GL}_{2n}$
whose Levi part is $\mathrm{GL}_{n_1} \times  \mathrm{GL}_{2n_2} \times \mathrm{GL}_{n_1}$.
Then ${}^{\beta_{\ell, n}} Z_{n_2, n_1}$ normalizes $R_{\ell, n}$ and $R_{\ell, n} \cap U_{n_1}$
and also $U_{n_1} = (R_{\ell, n} \cap U_{n_1}) \rtimes {}^{\beta_{\ell, n}} Z_{n_2, n_1}$.
Recall that 
\[
\beta_{\ell, n_2} = \begin{pmatrix}&I_{\ell}&&\\ &&&I_{n_2-\ell}\\ -I_{n_2-\ell}&&&\\ &&I_{\ell}& \end{pmatrix}.
\]
Then we may write $\omega_{n_1, n_2} \beta_{\ell, n} =\beta_{\ell, n_2} w^\prime$ with
\[
w^\prime = \begin{pmatrix} &&I_{n_1}\\ &I_{2n_2}&\\ -I_{n_1}&& \end{pmatrix}.
\]
Note that $w^\prime$ commutes with $g \in \mathrm{GL}_{2\ell}$, $r \in R_{\ell, n_2}$
and $H_\ell$.
Since $\mathrm{GL}_{2\ell}$ and $H_\ell$ normalize $U_{n_1}$ and $\psi^{-1}(z)$, we obtain
\begin{multline}
\label{8.10}
\int_{U_{n_1}} \int_{U_{\mathrm{GL}_{2\ell}} \backslash \mathrm{GL}_{2\ell}} \int_{Y_\ell \backslash H_\ell} \int_{R_{\ell, n_2}}
W(g) \omega_{\psi^{-1}, \mu}(hg)\phi(\xi_\ell) 
\\
\varphi_s(w^\prime u, I_{n_1}, I_{n_1}, \beta_{\ell, n_2}rhg, I_{n_2}, I_{n_2}) \psi(u)
\,dr \, dh \, dg\, du.
\end{multline}
The above manipulations are valid in the domain of absolute convergence of \eqref{8.9}
with the implicit parameter $\zeta$ as indicated in the beginning of Section~\ref{Proof of Proposition split}.
Then we apply the functional equation for $\pi \times (\tau_{21}, \tau_{22})$ to the inner $dr \, dh \, dg$-integration.
Formally, multiplying the last integral by $\gamma(s, \pi \times (\tau_{12}, \tau_{22}), \mu, \psi)$, we get
\begin{multline}
\label{8.11}
\int_{U_{n_1}} \int_{U_{\mathrm{GL}_{2\ell}} \backslash \mathrm{GL}_{2\ell}} \int_{Y_\ell \backslash H_\ell} \int_{R_{\ell, n_2}}
W(g) \omega_{\psi^{-1}, \mu}(hg)\phi(\xi_\ell, \xi_\ell) 
\\
\varphi_s^\prime(w^\prime u, I_{n_1}, I_{n_1}, \beta_{\ell, n_2}rhg, I_{n_2}, I_{n_2}) \psi(u)
\,dr \, dh \, dg\, du.
\end{multline}
\begin{Remark}
\label{verify conv}
The above argument is verified in a similar way as \cite{S1,S3}, \cite[4.10]{Ka2} and \cite[p.431]{Ka} as follows.
In the domain of absolute convergence of \eqref{8.9} (with the implicit parameter $\zeta$), \eqref{8.10} equals to \eqref{8.9}.
This gives analytic continuation of \eqref{8.10}, and \eqref{8.10} satisfies \eqref{3.6 inert} for all $s$.
Since \eqref{8.11} also satisfies \eqref{3.6 inert}, by the uniqueness \cite[Theorem~5.3]{BAS}, 
when $s$ is in the domain of absolute convergence of \eqref{8.11},
 \eqref{8.11} equals  \eqref{8.10} up to constant multiple depending on $s$.
We may take $W, f_{\tau,s}$ and $\phi$ so that $\mathcal{L}(W, f_{s}, \phi) =1$ by \cite[Proposition~4.1]{BAS} for any $s$.
%and in particular, $\mathcal{L}(W, f_s, \phi)$ is holomorphic for all $s$ in this domain.
For such data, we find that the above constant equals required $\gamma$-factor in the domain of the convergence.
Then $\gamma(s, \pi \times (\tau_{12}, \tau_{22}), \psi) \times \text{\eqref{8.10}} = \text{\eqref{8.11}}$,
and this also gives analytic continuation of \eqref{8.11}. Therefore, this identity holds for any $s$.
On the other hand, we find that 
\eqref{8.11} equals to $\mathcal{L}(W, f_{\varphi_{s}^\prime}, \phi)$ in the domain of absolute convergence of $\mathcal{L}(W, f_{\varphi_{s}^\prime}, \phi)$.
This gives analytic continuation of $\mathcal{L}(W, f_{\varphi_{s}^\prime}, \phi)$. Therefore, Proposition~\ref{8.5} holds for any $s$.
Finally, we note that when we consider a twist by parameter $\zeta$, the resulting representation may be reducible.
Then in order to study $\gamma$-factors, we need the uniqueness given in Theorem~\ref{GL unique} for reducible representations in the non-archimedean case.
On the other hand, in the archimedean case, we may choose $\zeta$ so that the resulting representation is irreducible. This also works in the non-archimedean case but it is not needed.
\end{Remark}
\subsubsection{$n_2 < \ell < n$}
Since ${}^{\beta_{\ell, n}} z$ normalizes $R_{\ell, n}$, we can change $z \beta_{\ell, n} r \mapsto \beta_{\ell, n} r {}^{\beta_{\ell, n}} z$.
Write 
\begin{multline*}
z = \begin{pmatrix} I_{n_2}&z_{11}&&&\\ &I_{\ell-n_2}&&&\\ &&1&&\\ &&&I_{n-\ell-1}&\\ &&&&I_n  \end{pmatrix}
\begin{pmatrix} I_{n_2}&&z_{21}&&\\ &I_{\ell-n_2}&&&\\ &&1&&\\ &&&I_{n-\ell-1}&\\ &&&&I_n  \end{pmatrix}
\\
\begin{pmatrix} I_{n_2}&&&z_{31}&\\ &I_{\ell-n_2}&&&\\ &&1&&\\ &&&I_{n-\ell-1}&\\ &&&&I_n  \end{pmatrix}
\begin{pmatrix} I_{n}&&&&\\ &I_{n-\ell-1}&&&\\ &&1&&\\ &&&I_{\ell-n_2}&z_{12}\\ &&&&I_{n_2}  \end{pmatrix}
\\
\begin{pmatrix} I_{n}&&&&\\ &I_{n-\ell-1}&&&\\ &&1&&z_{22}\\ &&&I_{\ell-n_2}&\\ &&&&I_{n_2}  \end{pmatrix}
\begin{pmatrix} I_{n}&&&&\\ &I_{n-\ell-1}&&&z_{32}\\ &&1&&\\ &&&I_{\ell-n_2}&\\ &&&&I_{n_2}  \end{pmatrix}
\\
=z_{11}z_{21}z_{31}z_{12}z_{22}z_{32}
\end{multline*}
and decompose $dz = dz_{11} \, dz_{21}\, dz_{31} \, dz_{12} \, dz_{22} \, dz_{32}$.
Because $\ell > n_2$, $\psi^{-1}(z) = \psi^{-1}(z_{11}+z_{12})$.
Moreover, we have ${}^{\beta_{\ell, n}} (z_{21} z_{22}) \in Y_{\ell}$
and $\omega_{\psi^{-1}, \mu}(hg)\phi(\xi_\ell, \xi_\ell) = \omega_{\psi^{-1}, \mu}({}^{\beta_{\ell, n}} (z_{21} z_{22}) hg)\phi(\xi_\ell, \xi_\ell)$
whence we can collapse the $d z_{21}$-integration and $d z_{22}$-integration into the $dh$-integration.
Define 
\[
Y_\ell^\circ =  \left\{ (0, (y_1, y_2); 0) : y_i \in F^{\ell} \text{ whose last $n_2$ coordinates of $y_i$ are zero} \right\}.
\]
Then we obtain
\begin{multline*}
\mathcal{L}(W, f_{\varphi_s}, \phi) = 
 \int_{U_{\mathrm{GL}_{2\ell}} \backslash \mathrm{GL}_{2\ell}} \int_{Y_\ell^\circ \backslash H_\ell} \int_{R_{\ell, n}} \int
W(g)  \omega_{\psi^{-1}, \mu}(hg)\phi(\xi_\ell, \xi_\ell)  
\\
\varphi_s(\omega_{n_1, n_2} \beta_{\ell,n} r ({}^{ \beta_{\ell,n}}(z_{11}z_{12} z_{31}z_{32})) h g, I_{n_1}, I_{n_1}, I_{2n_2}, I_{n_2}, I_{n_2}) 
\\
 \psi^{-1}(z_{11}+z_{12})\, dz_{11} \, dz_{12} \, dz_{31} \, dz_{32}
\, dr \, dh \, dg. 
\end{multline*}
Let $R_{\ell, n}^+$ denote the subgroup of elements $r ({}^{ \beta_{\ell,n}}z_{31} {}^{ \beta_{\ell,n}}z_{32})$.
Additionally, we have $\psi^{-1}(z_{11}+z_{12}) W(g) = W(z g)$ and 
\[
\omega_{\psi^{-1}, \mu}(hg)\phi(\xi_\ell, \xi_\ell) = \omega_{\psi^{-1}, \mu}( {}^{ \beta_{\ell,n}}z_{11} {}^{ \beta_{\ell,n}}z_{12} hg)\phi(\xi_\ell, \xi_\ell).
\]
Changing a variables, we can change $({}^{ \beta_{\ell,n}}z_{11} {}^{ \beta_{\ell,n}}z_{12}) h \mapsto h({}^{ \beta_{\ell,n}}z_{11} {}^{ \beta_{\ell,n}}z_{12})$.
Since
\[
{}^{ \beta_{\ell,n}}z_{11} {}^{ \beta_{\ell,n}}z_{12} = \begin{pmatrix} I_{n-\ell}&&&&&\\ &I_{n_2}&z_{11}&&&\\ &&I_{\ell-n_2}&&&\\  &&&I_{\ell-n_2}&z_{12}&\\  
&&&&I_{n_2}&\\  &&&&& I_{n-\ell}\\ \end{pmatrix},
\]
collapsing $dz_{11}$-integration and $dz_{12}$-integration into $dg$-integration, we get
\begin{multline*}
 \int_{(U_{\mathrm{GL}_{n_2}})^2 (U_{\mathrm{GL}_{\ell-n_2}})^2 U_{\ell} \backslash \mathrm{GL}_{2\ell}} \int_{Y_\ell^\circ \backslash H_\ell} \int_{R_{\ell, n}^+}
W(g)  \omega_{\psi^{-1}, \mu}(hg)\phi(\xi_\ell, \xi_\ell)  
\\
\varphi_s(\omega_{n_1, n_2} \beta_{\ell,n}  rh g, I_{n_1}, I_{n_1}, I_{2n_2}, I_{n_2}, I_{n_2}) 
\, dr \, dh \, dg. 
\end{multline*}
Here, $U_{\mathrm{GL}_{n_2}} \times U_{\mathrm{GL}_{\ell-n_2}} \times U_{\mathrm{GL}_{\ell-n_2}} \times U_{\mathrm{GL}_{n_2}}$ is embedded diagonally.
Then change a variable $g \mapsto {}^{\omega_{\ell-n_2, n_2}} g$ and we get
\begin{multline*}
 \int_{{}^{\omega_{n_2, \ell-n_2}}((U_{\mathrm{GL}_{n_2}})^2 (U_{\mathrm{GL}_{\ell -n_2}})^2 U_\ell) \backslash \mathrm{GL}_{2\ell}}
  \int_{Y_\ell^\circ \backslash H_\ell} \int_{R_{\ell, n}^+}
W({}^{\omega_{\ell-n_2, n_2}} g)  
\\
\omega_{\psi^{-1}, \mu}(h({}^{\omega_{\ell-n_2, n_2}} g))\phi(\xi_\ell, \xi_\ell)  
\varphi_s(\omega_{n_1, n_2} \beta_{\ell,n}  rh ({}^{\omega_{\ell-n_2, n_2}} g), I_{n_1}, I_{n_1}, I_{2n_2}, I_{n_2}, I_{n_2}) 
\, dr \, dh \, dg. 
\end{multline*}
Observe that $R^{\ell, n_2} X_{n_2}({}^{\omega_{n_2, \ell-n_2}}((U_{\mathrm{GL}_{n_2}})^2 (U_{\mathrm{GL}_{\ell -n_2}})^2 U_\ell) ) = U_{\mathrm{GL}_{2\ell}}$.
We may factor the integral
\[
 \int_{{}^{\omega_{n_2, \ell-n_2}}((U_{\mathrm{GL}_{n_2}})^2 (U_{\mathrm{GL}_{\ell -n_2}})^2 U_\ell) \backslash \mathrm{GL}_{2\ell}}
 = \int_{U_{\mathrm{GL}_{2\ell}} \backslash \mathrm{GL}_{2\ell}} \int_{R^{\ell, n_2}} \int_{X_{n_2}}.
\]
Further, we have $U_{\mathrm{GL}_{2\ell}} = ((U_{\mathrm{GL}_{\ell-n_2}})^2 \ltimes U_{\ell-n_2}) \rtimes U_{\mathrm{GL}_{2n_2}}$.
Then we may factor the integral
\[
 \int_{U_{\mathrm{GL}_{2\ell}} \backslash \mathrm{GL}_{2\ell}}
 = \int_{((U_{\mathrm{GL}_{\ell-n_2}})^2 \ltimes U_{\ell-n_2}) \mathrm{GL}_{2n_2} \backslash \mathrm{GL}_{2\ell} } \int_{U_{\mathrm{GL}_{2n_2}} \backslash \mathrm{GL}_{2n_2} }.
\]
Hence, our integral is
\begin{multline*}
 \int_{((U_{\mathrm{GL}_{\ell-n_2}})^2 \ltimes U_{\ell-n_2}) \mathrm{GL}_{2n_2} \backslash \mathrm{GL}_{2\ell} } \int_{U_{\mathrm{GL}_{2n_2}} \backslash \mathrm{GL}_{2n_2} }
 \int_{R^{\ell, n_2}} \int_{X_{n_2}}
  \int_{Y_\ell^\circ \backslash H_\ell} \int_{R_{\ell, n}^+}
  \\
W(({}^{\omega_{\ell-n_2, n_2}} (r^\prime x^\prime g^\prime))({}^{\omega_{\ell-n_2, n_2}} g))  
\omega_{\psi^{-1}, \mu}(h({}^{\omega_{\ell-n_2, n_2}} (r^\prime x^\prime g^\prime)({}^{\omega_{\ell-n_2, n_2}} g))\phi(\xi_\ell, \xi_\ell)  
\\
\varphi_s(\omega_{n_1, n_2} \beta_{\ell,n}  rh ({}^{\omega_{\ell-n_2, n_2}}( r^\prime x^\prime g^\prime)({}^{\omega_{\ell-n_2, n_2}} g), I_{n_1}, I_{n_1}, I_{2n_2}, I_{n_2}, I_{n_2}) 
\, dr \, dh \, \, dx^\prime \, dr^\prime \, dg^\prime \, dg. 
\end{multline*}
Since the conjugate by ${}^{\omega_{\ell-n_2, n_2}} (r^\prime x^\prime)$ normalize $Y_\ell^\circ$ and $H_\ell$,
we may change the variable $h ({}^{\omega_{\ell-n_2, n_2}} (r^\prime x^\prime)) \mapsto  ({}^{\omega_{\ell-n_2, n_2}} (r^\prime x^\prime))h$.
Note that we may write the $dg^\prime$-integration as the integration over $\overline{Q_{n_2}}$
and that $\overline{Q_{n_2}}$ normalizes $Y_\ell^\circ$ and $H_\ell$.
Then we can change the variable $h ({}^{\omega_{\ell-n_2, n_2}} g^\prime) \mapsto ({}^{\omega_{\ell-n_2, n_2}} g^\prime) h$.
Therefore, we get
\begin{multline*}
 \int_{((U_{\mathrm{GL}_{\ell-n_2}})^2 \ltimes U_{\ell-n_2}) \mathrm{GL}_{2n_2} \backslash \mathrm{GL}_{2\ell} } \int_{U_{\mathrm{GL}_{2n_2}} \backslash \mathrm{GL}_{2n_2} }
 \int_{R^{\ell, n_2}} \int_{X_{n_2}}
  \int_{Y_\ell^\circ \backslash H_\ell} \int_{R_{\ell, n}^+}
  \\
W(({}^{\omega_{\ell-n_2, n_2}} (r^\prime x^\prime g^\prime))({}^{\omega_{\ell-n_2, n_2}} g))  
\omega_{\psi^{-1}, \mu}(({}^{\omega_{\ell-n_2, n_2}} (x^\prime g^\prime) h ({}^{\omega_{\ell-n_2, n_2}} g))\phi(\xi_\ell, \xi_\ell)  
\\
\varphi_s(\omega_{n_1, n_2} \beta_{\ell,n}  r ({}^{\omega_{\ell-n_2, n_2}}( r^\prime x^\prime g^\prime) h ({}^{\omega_{\ell-n_2, n_2}} g), I_{n_1}, I_{n_1}, I_{2n_2}, I_{n_2}, I_{n_2}) 
\, dr \, dh \, \, dx^\prime \, dr^\prime \, dg^\prime \, dg. 
\end{multline*}
Here, we note that from the definition, it is easy to check that
\[
\omega_{\psi^{-1}, \mu}({}^{\omega_{\ell-n_2, n_2}} r^\prime)\phi^\prime(\xi_\ell, \xi_\ell)  = \phi^\prime(\xi_\ell, \xi_\ell)
\]
for any $\phi^\prime \in \mathcal{S}(F^\ell \times F^\ell)$, and we used this identity in the above computation.

For any $\phi^\prime \in \mathcal{S}(F^\ell \times F^\ell)$, let $\phi_{n_2}^\prime \in \mathcal{S}(F^{n_2} \times F^{n_2})$ be the function defined by 
\begin{equation}
\label{phi prime def}
\phi_{n_2}^\prime(\xi_1, \xi_2) = \phi^\prime((\xi_1, 0, \dots, 0, 1), (\xi_2, 0, \dots, 0, 1)).
\end{equation}
As before, we write an element $x^\prime$ of $X_{n_2}$ as $x^\prime = ((x_1, x_2), 0)$ with $x_i \in F^{n_2}$.
Then we have $\omega_{\psi^{-1}, \mu}({}^{\omega_{\ell-n_2, n_2}} (x^\prime ) )\phi(\xi_\ell, \xi_\ell)  = \phi_{n_2}^\prime(x_1, x_2)$,
and the last integral becomes
\begin{multline*}
 \int_{((U_{\mathrm{GL}_{\ell-n_2}})^2 \ltimes U_{\ell-n_2}) \mathrm{GL}_{2n_2} \backslash \mathrm{GL}_{2\ell} }
  \int_{Y_\ell^\circ \backslash H_\ell}  \int_{U_{\mathrm{GL}_{2n_2}} \backslash \mathrm{GL}_{2n_2} }
 \int_{R^{\ell, n_2}} \int_{X_{n_2}}
 \int_{R_{\ell, n}^+}
  \\
W(({}^{\omega_{\ell-n_2, n_2}} (r^\prime x^\prime g^\prime))({}^{\omega_{\ell-n_2, n_2}} g))  
\omega_{\psi^{-1}, \mu}(g^\prime) (\omega_{\psi^{-1}, \mu}(h ({}^{\omega_{\ell-n_2, n_2}} g))\phi)_{n_2}(x_1, x_2)  
\\
\varphi_s(\omega_{n_1, n_2} \beta_{\ell,n}  r ({}^{\omega_{\ell-n_2, n_2}}( r^\prime x^\prime g^\prime) h ({}^{\omega_{\ell-n_2, n_2}} g), I_{n_1}, I_{n_1}, I_{2n_2}, I_{n_2}, I_{n_2}) 
\, dr \, dh \, \, dx^\prime \, dr^\prime \, dg^\prime \, dg. 
\end{multline*}
Set $w_{r^\prime, x^\prime} = {}^{\omega_{\ell-n_2, n_2}}( r^\prime x^\prime)$.
By a direct computation, we can find $b_{r^\prime, x^\prime} \in \mathrm{GL}_{2n}$ and $r_{r^\prime, x^\prime} \in R_{\ell, n}^+$
such that $w^{r^\prime, x^\prime} r =b_{r^\prime, x^\prime} r_{r^\prime, x^\prime}$ and 
${}^{(\omega_{n_1, n_2} \beta_{\ell, n})^{-1}} b_{r^\prime, x^\prime} \in U_{\mathrm{GL}_n} \times U_{\mathrm{GL}_n}$.
Further, $r \mapsto r^\prime$ gives an isomorphism of $R_{\ell, n}^+$.
Then, we obtain
\begin{multline*}
\varphi_s({}^{(\omega_{n_1, n_2} \beta_{\ell, n})^{-1}} ( w_{r^\prime, x^\prime} b_{r^\prime, x^\prime}), I_{n_1}, I_{n_1}, I_{2n_2}, I_{n_2}, I_{n_2} )
\\
=\varphi_s( I_{2n}, I_{n_1}, I_{n_1}, I_{2n_2}, I_{n_2}, I_{n_2} ),
\end{multline*}
and changing the variable $r_{r^\prime, x^\prime} \mapsto r$, we get
\begin{multline*}
 \int_{((U_{\mathrm{GL}_{\ell-n_2}})^2 \ltimes U_{\ell-n_2}) \mathrm{GL}_{2n_2} \backslash \mathrm{GL}_{2\ell} }
  \int_{Y_\ell^\circ \backslash H_\ell}  \int_{U_{\mathrm{GL}_{2n_2}} \backslash \mathrm{GL}_{2n_2} }  \int_{R_{\ell, n}^+} 
 \int_{R^{\ell, n_2}} \int_{X_{n_2}}
  \\
W(({}^{\omega_{\ell-n_2, n_2}} (r^\prime x^\prime g^\prime))({}^{\omega_{\ell-n_2, n_2}} g))  
\omega_{\psi^{-1}, \mu}(g^\prime)( \omega_{\psi^{-1}, \mu}( h ({}^{\omega_{\ell-n_2, n_2}} g))\phi)_{n_2} (x_1, x_2)  
\\
\varphi_s(\omega_{n_1, n_2} \beta_{\ell,n}  r ({}^{\omega_{\ell-n_2, n_2}}(g^\prime) h ({}^{\omega_{\ell-n_2, n_2}} g), I_{n_1}, I_{n_1}, I_{2n_2}, I_{n_2}, I_{n_2}) 
\, dx^\prime \, dr^\prime  \, dr \, dg^\prime \, dh  \, dg. 
\end{multline*}
Also, one sees that ${}^{\omega_{\ell-n_2, n_2}}g^\prime$ normalize $R_{\ell, n}^{+}$ without changing $dr$.
Further, since 
\[
\omega_{n_1, n_2} \beta_{\ell, n} \omega_{n_2, \ell-n_2} =
\begin{pmatrix} &I_{\ell-n_2}&&&&\\ &&&&&I_{n-\ell}\\ &&I_{n_2}&&&\\ &&&I_{n_2}&&\\ -I_{n-\ell}&&&&&\\ &&&&I_{\ell-n_2}&\end{pmatrix},
\]
conjugating by $\omega_{n_1, n_2} \beta_{\ell, n}$ we get
\begin{multline*}
 \int_{((U_{\mathrm{GL}_{\ell-n_2}})^2 \ltimes U_{\ell-n_2}) \mathrm{GL}_{2n_2} \backslash \mathrm{GL}_{2\ell} }
  \int_{Y_\ell^\circ \backslash H_\ell}   \int_{R_{\ell, n}^+} \int_{U_{\mathrm{GL}_{2n_2}} \backslash \mathrm{GL}_{2n_2}  } 
 \int_{R^{\ell, n_2}} \int_{X_{n_2}}
  \\
W(({}^{\omega_{\ell-n_2, n_2}} (r^\prime x^\prime g^\prime))({}^{\omega_{\ell-n_2, n_2}} g))  
\omega_{\psi^{-1}, \mu}(g^\prime) (\omega_{\psi^{-1}, \mu}( h ({}^{\omega_{\ell-n_2, n_2}} g))\phi)_{n_2}(x_1, x_2)  
\\
\varphi_s(\omega_{n_1, n_2} \beta_{\ell,n}  r h ({}^{\omega_{\ell-n_2, n_2}} g), I_{n_1}, I_{n_1}, g^\prime, I_{n_2}, I_{n_2}) 
\, dx^\prime \, dr^\prime  \, dg^\prime  \, dr \, dh  \, dg.  
\end{multline*}
Now we apply the functional equation for $\pi \times (\tau_{21} \otimes \tau_{22})$
to the inner $dx^\prime \, dr^\prime  \, dg^\prime$-integration, and the required identity readily follows.
\subsubsection{$n_2 = \ell < n$}
The proof is similar to the case $n_2 < \ell < n$, we use the same notation and describe the necessary modifications.
Here, we decompose $z = z_{21} z_{22} z_{31} z_{32}$ with
\[
z_{21} = \begin{pmatrix} I_{\ell}&z_{21}&&\\ &1&&\\ &&I_{n-\ell-1}&\\ &&&I_n \end{pmatrix},
z_{31} = \begin{pmatrix} I_{\ell}&&z_{31}&\\ &1&&\\ &&I_{n-\ell-1}&\\ &&&I_n \end{pmatrix},
\]
\[
z_{22} = \begin{pmatrix} I_n&&&\\ &I_{n-\ell-1}&&\\ &&1&z_{22}\\ &&&I_\ell \end{pmatrix},
z_{32} = \begin{pmatrix} I_n&&&\\ &I_{n-\ell-1}&&z_{32}\\ &&1&\\ &&&I_\ell \end{pmatrix}.
\]
Regard $z_{21}$ (resp. $z_{22}$) as a column (resp. row) vector and denote 
$z_{21} = {}^{t}(a_{11}, \dots, a_{1\ell})$ (resp. $z_{22} = (a_{21}, \dots, a_{2\ell})$).
Then 
\[
{}^{\beta_{\ell, n}} z_{21} = ((0, 0), (-a_{21}, \dots, -a_{2\ell}, a_{1 \ell}, \dots a_{11}); 0)
\]
whence 
\[
\omega_{\psi^{-1}, \mu}(hg)(\xi_\ell, \xi_\ell) = \psi(a_{1\ell}-a_{21})\omega_{\psi^{-1}, \mu}( ({}^{\beta_{\ell, n}}(z_{21}z_{22}))hg)(\xi_\ell, \xi_\ell).
\]
Also, $\psi^{-1}(z) = \psi^{-1}(a_{1\ell}-a_{21})$, and thus we can collapse the $dz_{21} dz_{22}$-integration into $dh$-integration.
Then since $\mathrm{GL}_{2\ell}$ normalize $H_\ell$, changing a variable $h \mapsto ghg^{-1}$, we get
\begin{multline*}
\int_{H_\ell} \int_{R_{\ell, n}^+} \int_{U_{ \mathrm{GL}_{2\ell}}\backslash \mathrm{GL}_{2\ell}} W(g) \omega_{\psi^{-1}, \mu}(gh)\phi(\xi_\ell, \xi_\ell)
\\
\varphi_s(\omega_{n_1, n_2} \beta_{\ell, n}th, I_{n_1}, I_{n_1}, g, I_{n_2}, I_{n_2}) \,dg \, dr \, dh.
\end{multline*}
Now we apply the functional equation for $\pi \times (\tau_{21}, \tau_{22})$.
\subsubsection{$\ell >n$}
Note that for $z \in Z_{n_2, n_1}$ and $x_i \in X_n$, we have $\psi^{-1}(z) W(g) = W({}^{\omega_{\ell-n,n}} z g) $ and 
\[
 \omega_{\psi^{-1}, \mu}(g) \phi(x_1, x_2) =  \omega_{\psi^{-1}, \mu}( (x_1,x_2) g) \phi(0, 0)
=  \omega_{\psi^{-1}, \mu}(z(x_1,x_2)g) \phi(0,0)
\]
Since $z$ normalizes $R^{\ell,n}$ and $X_n^2$, after change of variables in $x_i$,
we can collapse the $dz$-integration into $dg$.
Write $U_{\mathrm{GL}_{2n}} = ((N_{n_1} \cdot N_{n_2}) \ltimes Z_{n_2, n_1}) \ltimes U_n$,
where
\[
N_{n_1} = \left\{\mathrm{diag}(I_{n_2}, z_1, z_2, I_{n_2}) : z_i \in U_{\mathrm{GL}_{n_1}} \right\}
\]
and
\[
N_{n_2} = \left\{\mathrm{diag}(z_1,I_{n_1}, I_{n_1}, z_2) : z_i \in U_{\mathrm{GL}_{n_2}} \right\}.
\]
Then our local zeta integral is equal to
\begin{multline*}
 \int_{(N_{n_1} N_{n_2} U_n) \backslash \mathrm{GL}_{2n}}  \int_{R^{\ell, n}}\int_{X_n^2}
W({}^{\omega_{\ell-n, n}}(r(x_1, x_2)g)) 
\\
\varphi_s(\omega_{n_1, n_2}g, I_n, I_n)  \omega_{\psi^{-1}, \mu}(g) \phi(x_1, x_2) \, dx_1\, dx_2 \, dr \, dg .
\end{multline*}
Decompose $R^{\ell, n} = R_1^{\ell, n}  R_2^{\ell, n}$, with
\[
R_1^{\ell, n} = \left\{ \begin{pmatrix}I_{\ell-n-1}&&&r_1&&&&\\ &1&&&&&&\\ &&I_{n_2}&&&&&\\ &&&I_{n_1}&&&& \\
&&&&I_{n_1}&&&r_2\\ &&&&&I_{n_2}&&\\ &&&&&&1&\\ &&&&&&&I_{\ell-n-1}\end{pmatrix}\right\},
\]
\[
R_2^{\ell, n} =  \left\{ \begin{pmatrix}I_{\ell-n-1}&&r_1&&&&&\\ &1&&&&&&\\ &&I_{n_2}&&&&&\\ &&&I_{n_1}&&&& \\
&&&&I_{n_1}&&&\\ &&&&&I_{n_2}&&r_2\\ &&&&&&1&\\ &&&&&&&I_{\ell-n-1}\end{pmatrix}\right\}.
\]
Factor the $dg$-integration through $\overline{Z_{n_2, n_1}}$
and put $N_{n_1, n_2} = \overline{Z_{n_2, n_1}} N_{n_1} N_{n_2}$.
We note that 
${}^{\omega_{\ell-n,n}} (\overline{Z_{n_2, n_1}} \cdot R_2^{\ell, n}) ={}^{\omega_{\ell-n_2, n_2}} R^{\ell, n_2}$.
For $z \in \overline{Z_{n_2, n_1}}$ and $r_1 \in R_1^{\ell, n}$, ${}^z r_1 = r_{2,z} r_1$
with $r_{2,z} \in R_2^{\ell, n}$.
Therefore the integral becomes
\begin{multline*}
 \int_{(N_{n_1, n_2} U_n) \backslash \mathrm{GL}_{2n}}  \int_{R_1^{\ell, n}} \int_{R^{\ell, n_2}} \int_{X_n^2}
W({}^{\omega_{\ell-n_2, n_2}}r {}^{\omega_{\ell-n, n}}(r_1(x_1, x_2)g)) 
\varphi_s(\omega_{n_1, n_2}g, I_n, I_n)  
\\
\omega_{\psi^{-1}, \mu}(g) \phi(x_1, x_2) \, dx_1\, dx_2 \, dr \, \, dr_1 \, dg. 
\end{multline*}
For $a \in X_n$,  write $a= (a_2, a_1)$ with $a_i \in F^{n_i}$
and let $X_n^\circ$ be the subgroups of the elements $(0, a_1)$.
Then for $((a_2, 0), (b_2, 0))$ (regarded as an element of $\mathrm{GL}_{2\ell}$),
${}^{\omega_{n_2, \ell-n_2}} ({}^{\omega_{\ell-n,n}}((a_2, 0), (b_2, 0))) = (a_2, b_2)$,
where on the right-hand side we regard $(a_2, b_2)$ as an element of 
$X_{n_2}^2 < U_{\mathrm{GL}_{2(n_2+1)}} < Q_{\ell-n_2-1} < \mathrm{GL}_{2\ell}$.
Moreover, $(a_2, b_2)$ commutes with any element of $R_1^{\ell, n}$.
Hence, we get
\begin{multline*}
 \int_{(N_{n_1, n_2} U_n) \backslash \mathrm{GL}_{2n}}  \int_{R_1^{\ell, n}} \int_{R^{\ell, n_2}}
\int_{(X_n^\circ)^2}  \int_{X_{n_2}^2}
W({}^{\omega_{\ell-n_2, n_2}}(r(a_2, b_2) {}^{\omega_{\ell-n, n}}(r_1(x_1, x_2)g)) 
\\
\varphi_s(\omega_{n_1, n_2}g, I_n, I_n)  \omega_{\psi^{-1}, \mu}((a_2, b_2)(x_1, x_2)g) \phi(0, 0) \, dx_1\, dx_2 \, da_2 \, db_2
 \, dr \, \, dr_1\, dg .
\end{multline*}
Now factor the $dg$-integration through $G_{n_2}^\prime:= {}^{\omega_{n_1, n_2}} \mathrm{GL}_{2n_2}$,
where $\mathrm{GL}_{2n_2}$ is regarded as a subgroup of $\mathrm{GL}_{2\ell}$.
If $g^\prime \in G_{n_2}^\prime$, $g^\prime$ commutes with $(x_1, x_2)$ and $r_1$ for all $x_i \in X_n^\circ$ and $r_1 \in R^{\ell, n}$
and also ${}^{\omega_{\ell-n, n}} g^\prime = {}^{\omega_{\ell-n_2, n_2}} g^\prime$.
Thus our integral equals
\begin{multline*}
 \int_{G_{n_2}^\prime(N_{n_1, n_2} U_n) \backslash \mathrm{GL}_{2n}}  \int_{R_1^{\ell, n}} 
\int_{(X_n^\circ)^2} 
\int_{U_{G_{n_2}}^\prime \backslash G_{n_2}^\prime}\int_{R^{\ell, n_2}} \int_{X_{n_2}^2}
\\
W({}^{\omega_{\ell-n_2, n_2}}(r(a_2, b_2)g^\prime ({}^{\omega_{\ell-n, n}}(r_1(x_1, x_2)g))) 
\varphi_s(\omega_{n_1, n_2}g, I_{n_1}, g^\prime, I_n)  
\\
\omega_{\psi^{-1}, \mu}(g^\prime)  (\omega_{\psi^{-1}, \mu}((x_1, x_2)g) \phi)(a_2, 0, b_2, 0)
da_2 \, db_2 \, dr \, dg^\prime  \, dx_1\, dx_2 \, 
 \, dr_1 \, dg .
\end{multline*}
Applying the functional equation for $\pi \times (\tau_{2,1}, \tau_{2,2})$ to the inner $da_2 \, db_2 \, dr \, dg^\prime$-integration,
we obtain the required identity.
\subsubsection{$\ell =n$}
In this case, our claim is proved in a similar way as the case of $\ell > n$.
Factoring the $dg$-integration through $\overline{Z_{n_2, n_1}}$ yields that 
\begin{multline*}
\int_{N_{n_1, n_2} U_{\mathrm{GL}_{2\ell}} \backslash \mathrm{GL}_{2\ell}} \int_{\overline{Z_{n_2, n_1}} }W(zg) \omega_{\psi^{-1}, \mu}(zg)\phi(\xi_\ell, \xi_\ell) 
\\
\varphi_s(\omega_{n_1, n_2}g,I_{n_1}, I_{n_1} I_{2n_2}, I_{n_2}, I_{n_2}) \, dz \, dg.
\end{multline*}
Decompose 
\begin{multline*}
z= \begin{pmatrix} I_{n_2}&&&&&\\ r_1&I_{n_1-1}&&&&\\ &&1&&&\\ &&&1&&\\ &&&&I_{n_1-1}&\\ &&&&r_2&I_{n_2} \end{pmatrix}
 \begin{pmatrix} I_{n_2}&&&&&\\ &I_{n_1-1}&&&&\\ a_2&&1&&&\\ &&&1&&\\ &&&&I_{n_1-1}&\\ &&&b_2&&I_{n_2} \end{pmatrix}
 \\
= (r_1, r_2)(a_2, b_2)
\end{multline*}
and set $dz = da_2 db_2 dr_1 dr_2$. Using 
\begin{multline*}
\omega_{\psi^{-1}, \mu}((r_1, r_2)(a_2, b_2)g) \phi(\xi_\ell, \xi_\ell)= \omega_{\psi^{-1}, \mu}((a_2, b_2)g) \phi(\xi_\ell, \xi_\ell)
\\
= \omega_{\psi^{-1}, \mu}(g) \phi((a_2, 0, \dots, 0, 1), (-{}^{t}b_2, 0, \dots, 0, 1))
\end{multline*}
where $(a_2, 0, \dots, 0, 1) , (-{}^{t}b_2, 0, \dots, 0, 1)\in F^\ell$, as in the previous section, we see that this integral is equal to
\begin{multline*}
\int_{N_{n_1, n_2} U_{\mathrm{GL}_{2\ell}} \backslash \mathrm{GL}_{2\ell}}
\int_{R^{\ell, n_2} } \int_{X_{n_2}^2}  
W({}^{\omega_{\ell-n_2, n_2}}((r_1, r_2)(a_2, b_2))g)  \omega_{\psi^{-1}, \mu}(g) 
\\
\phi((a_2, 0, \dots, 0, 1), (-{}^{t}b_2, 0, \dots, 0, 1))
\varphi_s(\omega_{n_1, n_2}g,I_{n_1}, I_{n_1} I_{2n_2}, I_{n_2}, I_{n_2})  
\\
\, da_2 \, db_2 \, dr_1 \, dr_2 \, dg.
\end{multline*}

Now factor the $dg$-integration through $G_{n_2}^{\prime \prime} := {}^{\omega_{\ell-n_2, n_2}} \mathrm{GL}_{2n_2}$,
where $\mathrm{GL}_{2n_2}$ is regarded as a subgroup of $\mathrm{GL}_{2\ell}$.
%We readily apply the functional equation for $\pi \times (\tau_{2,1}, \tau_{2,2})$ to the inner $da_2 \, db_2 \, dr \, dg^\prime$-integration.
For any $\phi^\prime \in \mathcal{S}(F^\ell)$, denote by  $\phi^\prime_{n_2} \in \mathcal{S}(F^{n_2})$ the function defined in \eqref{phi prime def}.
Then we obtain
\begin{multline*}
\int_{G_{n_2}^{\prime \prime} N_{n_1, n_2} U_{\mathrm{GL}_{2\ell}} \backslash \mathrm{GL}_{2\ell}} \int_{U^{\prime \prime}_{G_{n_2}} \backslash G_{n_2}^{\prime \prime}}
\int_{R^{\ell, n_2} } \int_{X_{n_2}^2}  
W({}^{\omega_{\ell-n_2, n_2}}((r_1, r_2)(a_2, b_2)g^\prime)g) 
\\
 (\omega_{\psi^{-1}, \mu}(g) \phi)_{n_2}(a_2, -{}^{t}b_2)
\varphi_s(\omega_{n_1, n_2}g,I_{n_1}, I_{n_1} g^\prime, I_{n_2}, I_{n_2})  \, da_2 \, db_2 \, dg^\prime \, dr_1 \, dr_2 \, dg.
\end{multline*}

Required identity follows by applying the functional equation for $\pi \times (\tau_{2,1}, \tau_{2,2})$ to the inner $da_2 \, db_2 \, dg^\prime \, dr_1 \, dr_2$-integration.
This completes the proof of Proposition~\ref{M1} in the split case.
%This implies that there is a constant $C > 0$ such that for any $\zeta$ with $\mathrm{Re}(\zeta) > C$,
%our integral converges. Hence, we may fix such $\zeta$ for our computation.
%%%%%%%%%%%%%%%%%%%%%%%%%%%%%%%%%%%%%%%%%%%%%%%
%
%
%
%
%
%
%
%
%
%
%
%
%
%
%%%%%%%%%%%%%%%%%%%%%%%%%%%%%%%%%%%%%%%%%%%%%%%
\subsection{Proof of Proposition~\ref{M2} in the split case}
In this section, we shall show  Proposition~\ref{M2} in the split case.
We will study the case $r \leq \ell < n$ in Section~\ref{sect 4.1}, the case $r=\ell > n$ in Section~\ref{9.2}
and the case $r =\ell=n$ in Section~\ref{s:r=l=n}. The missing case $n, r<\ell$ follows from the first case and Proposition~\ref{M1}.
Indeed, we consider $\mathrm{Ind}_{P_{n, m}}^{\mathrm{GL}_{n+m}} (\tau_1\otimes \tau_0) \otimes \mathrm{Ind}_{P_{n, m}}^{\mathrm{GL}_{n+m}} (\tau_1\otimes \tau_0^\prime)$ instead of $\tau = \tau_1 \otimes \tau_2$
where $\tau_0$ and $\tau_0^\prime$ are irreducible representations of $\mathrm{GL}_m$ with $m > \ell$.
Then we may apply the multiplicativity of the first case. Moreover, applying Proposition~\ref{M1}, the required multiplicativity follows.

Let $\sigma = (\sigma_1, \sigma_2)$ and $\pi^\prime$ be representations of $\mathrm{GL}_{r} \times \mathrm{GL}_r$
and $\mathrm{GL}_{2(\ell-r)}$, which are irreducible in the Archimedean case or quotients of representations 
parabolically induced from supercuspidal representations in the non-archimedean case.
Let $\pi =\mathrm{Ind}^{\mathrm{GL}_{2\ell}}_{Q_r}(\sigma_1 \otimes \pi^\prime \otimes \sigma_2^\ast)$.
%Suppose that $\sigma_i$ and $\pi^\prime$ are of Whittaker type.
In this section, we prove
\begin{multline}
\label{mult II}
\gamma(s, \pi \times \tau, \mu, \psi)
\\
= \omega_{\sigma_1}(-1)^n \omega_{\sigma_2}(-1)^n \omega_{\tau_1}(-1)^r \omega_{\tau_2}(-1)^r
\gamma(s, (\sigma_1 \otimes \mu) \times \tau_1, \psi) \gamma(s, (\sigma_1 \otimes \mu)^\ast \times \tau_2, \psi)
\\
\gamma(s, (\sigma_2 \otimes \mu) \times \tau_1, \psi) \gamma(s, (\sigma_2 \otimes \mu)^\ast \times \tau_2, \psi)
\gamma(s, \pi^\prime \times \tau, \mu, \psi).
\end{multline}
For $\zeta \in \mathbb{C}$, 
%we put $\varepsilon_i = \sigma_i |\cdot|^{\zeta}$ and 
we consider 
\[
\pi_\zeta =\mathrm{Ind}^{\mathrm{GL}_{2\ell}}_{Q_r}(\sigma_1 |\cdot|^{\zeta} \otimes \pi^\prime \otimes \sigma_2^\ast |\cdot|^{-\zeta}).
\]
As in the previous section, we realize $\pi_\zeta$ in the space of functions on $\mathrm{GL}_{2\ell} \times \mathrm{GL}_r \times \mathrm{GL}_{2(\ell-r)} \times \mathrm{GL}_r$. 
Then for each $\varphi \in \pi_\zeta$, we realize the Whittaker function of $\mathcal{W}(\pi_\zeta, \psi^{-1})$ by 
%\[
%\mathrm{Ind}_{\overline{Q_r}}^{\mathrm{GL}_{2\ell}}\left(\mathcal{W}(\varepsilon_1, \psi^{-1}) \otimes \mathcal{W}(\pi^\prime, \psi^{-1}) \otimes \mathcal{W}(\varepsilon_2^\ast, \psi) \right)
%\]
%defines a Whittaker function $W_{\varphi} \in \mathcal{W}(\pi_\zeta, \psi^{-1})$ via
\begin{equation}
\label{W_phi def}
W_{\varphi}(g) = \int_{U_r} \varphi(ug, I_{r}, I_{2(\ell-r)}, I_r) \psi(u) \, du,
\end{equation}
with
\[
\psi(u) = 
\left\{
\begin{array}{ll}
\psi(u_{r, r+1} - u_{2\ell-r, 2\ell-r+1}), &r<\ell \\
\\
\psi^{-1}(u_{\ell, \ell+1})& r=\ell,
\end{array}
\right.
\]
which converges absolutely when $\mathrm{Re}(\zeta) \gg 0$. 
Here, we recall that we take a character $\psi_2$ 
as explained in the beginning of Section~\ref{proof:s}.
%Here, $\overline{Q_r}$ denotes the opposite parabolic subgroup of $Q_r$
%with the Levi decomposition $\overline{Q_r} = M_r \overline{U_r}$.
Then as in the previous section, we need additional parameter $\zeta$ in order to justify our argument.
As we noted above, we only consider a formal computation, and thus hereafter we omit the parameter $\zeta$.

For our proof, let us recall the functional equation by Jacquet--Piatetski-Shapiro--Shalika~\cite{JPSS}
in the non-archimedean case and Jacquet--Shalika~\cite{JS2} in the archimedean case (see also \cite[p.70]{S1}).
Let $\eta$ and $\eta^\prime$ be representations of $\mathrm{GL}_k$ and $\mathrm{GL}_{k^\prime}$
with $k <k^\prime$.
Then for $0 \leq j \leq k^\prime - k-1$, $W \in \mathcal{W}(\eta,\psi^{-1})$ and $W^\prime \in \mathcal{W}(\eta^\prime, \psi)$,
the $\gamma$-factor $\gamma(s, \eta^\prime \times \eta, \psi)$ is defined by
\begin{multline}
\label{9.1}
\omega_\eta(-1)^{k^\prime-1} \gamma(s, \eta^\prime \times \eta, \psi)
\times \int_{U_{\mathrm{GL}_k} \backslash \mathrm{GL}_k} \int_{M_{j \times k}} W(a) 
\\
W^\prime \left( \begin{pmatrix} a&&\\ r&I_j&\\ &&I_{k^\prime-k-j}\end{pmatrix}\right) |\det a|^{s-\frac{1}{2}(k^\prime-k)}\, dr \, da
\\
=
\int_{U_{\mathrm{GL}_k} \backslash \mathrm{GL}_k} \int_{M_{k \times k^\prime-k-j-1}} W(a) 
W^\prime \left(\Omega_{k^\prime-k,k} \begin{pmatrix} I_k&&m\\ &I_{j+1}&\\ &&I_{k^\prime-k-j-1}\end{pmatrix}\right)
\\
 |\det a|^{s-\frac{1}{2}(k^\prime-k)+j}\, dm \, da.
\end{multline}
Moreover, we know the following functional equation (e.g. see \cite[p.209]{S2}))
\begin{equation}
\label{fq gamma}
\gamma(s, \eta^\prime \times \eta, \psi)\gamma(1-s, (\eta^\prime)^\vee \times \eta^\vee, \psi) = \omega_{\eta^\prime}(-1)^k \omega_{\eta}(-1)^{k^\prime}.
\end{equation}
%%%%%%%%%
\subsubsection{$r \leq \ell < n$}
\label{sect 4.1}
Putting \eqref{W_phi def} into our local integral, we get
\begin{multline*}
\mathcal{L}(W_{\varphi}, f_s, \phi)
= \int_{U_{ \mathrm{GL}_{2\ell}} \backslash \mathrm{GL}_{2\ell}} \int_{Y_\ell \backslash H_\ell} \int_{R_{\ell, n}}\int_{U_r} 
\varphi(ug, I_{r}, I_{2(\ell-r)}, I_r)  \psi(u)  
\\
f_s(\beta_{\ell, n} rh g, I_n, I_n)  \omega_{\psi^{-1}, \mu}(h g) \phi(\xi_\ell, \xi_\ell) \, du \, dr \, dh \, dg. 
\end{multline*}
%Let $B^\circ$ denote the standard parabolic subgroup of $\mathrm{GL}_{2\ell}$ whose Levi part is $\mathrm{GL}_r \times B_{2\ell-r} \times \mathrm{GL}_r$
%and its unipotent radical is $\overline{U_r}$.
%Then we have $\mathrm{GL}_{2\ell} = B^\circ K_{\mathrm{GL}_{2\ell}}$ and 
%\[
%\int_{(U_{\mathrm{GL}_r})^2 \rtimes U_{\mathrm{GL}_{2(\ell-r)}} \backslash \mathrm{GL}_{2\ell}}
%=\int_{K_{\mathrm{GL}_{2\ell}}} \int_{(U_{\mathrm{GL}_r})^2 \rtimes U_{\mathrm{GL}_{2(\ell-r)}} \backslash B^\circ}
%= \int_{K_{\mathrm{GL}_{2\ell}}} \int_{T_{\mathrm{GL}_{2(\ell-r)}}} \int_{\overline{U_r} \cap \overline{U_\ell}} \int_{(U_{\mathrm{GL}_r} \backslash \mathrm{GL}_r)^2}
%\int_{\overline{Z_{r, \ell-r}}}.
%\]
Using $U_{ \mathrm{GL}_{2\ell}} = ((U_{\mathrm{GL}_r})^2 \ltimes U_r) \rtimes U_{ \mathrm{GL}_{2(\ell-r)}}$,
we may collapse the $du$-integration into $dg$, and thus the above integral becomes
\begin{multline*}
 \int_{K_{\mathrm{GL}_{2\ell}}} \int_{T_{2\ell-2r}} \int_{\overline{U_r} \cap \overline{U_\ell}} \int_{(U_{\mathrm{GL}_r} \backslash \mathrm{GL}_r)^2}
\int_{\overline{Z_{r, \ell-r}}}  \int_{Y_\ell \backslash H_\ell} \int_{R_{\ell, n}}
\\
\varphi(\mathrm{diag}(a_1, I_{2(\ell-r)}, a_2)mvtk, I_{r}, I_{2(\ell-r)}, I_r) 
f_s(\beta_{\ell, n} rh \mathrm{diag}(a_1, I_{2(\ell-r)}, a_2) mvtk, I_n, I_n)  
\\
\omega_{\psi^{-1}, \mu}(h \mathrm{diag}(a_1, I_{2(\ell-r)}, a_2)mvtk) \phi(\xi_\ell, \xi_\ell) 
\, dr \, dh \, dm \, da_1 \, da_2 \, dv \, dt \, dk.
\end{multline*}
Here, $T_{2\ell-2r}$ is the group of diagonal matrices of $\mathrm{GL}_{2\ell-2r}$, which is regarded as a subgroup of $\mathrm{GL}_{2\ell}$,
and $\overline{Z_{r, \ell-r}} = {}^{t}Z_{r, \ell-r}$.
As in the previous section, for $\phi^\prime \in \mathcal{S}(F^\ell \times F^\ell)$ and $\xi_1, \xi_2 \in F^r$, we set 
\[
\phi^\prime_r(\xi_1, \xi_2) =
\left\{
\begin{array}{ll}
\phi^\prime((\xi_1,0, \dots, 0,1),(\xi_2,0, \dots, 0,1) ) & r<\ell
\\
&\\
\phi^\prime(\xi_1,\xi_2) &r=\ell.
\end{array}
\right.
\]
Shifting $\mathrm{diag}(a_1, I_{2(\ell-r)}, a_2)$ to the left gives
\begin{multline}
\label{p.49}
 \int_{K_{\mathrm{GL}_{2\ell}}} \int_{T_{2(\ell-r)}} \int_{\overline{U_r} \cap \overline{U_\ell}} \int_{(U_{\mathrm{GL}_r} \backslash \mathrm{GL}_r)^2}
\int_{\overline{Z_{r, \ell-r}}}  \int_{Y_\ell \backslash H_\ell} \int_{R_{\ell, n}}
\varphi(vtk, a_1, I_{2(\ell-r)}, a_2) 
\\
 |\det a_1 a_2^{-1}|^{\frac{1}{2}(r-n)+s}
f_s(\beta_{\ell,n} rhvtk, \left( \begin{smallmatrix}a_1&&\\z_1&I_{\ell-r}&\\ &&I_{n-\ell} \end{smallmatrix}\right), 
\left( \begin{smallmatrix}I_{n-\ell}&&\\&I_{\ell-r}&\\ &a_2 z_2&a_2 \end{smallmatrix}\right)) 
\\
\left( \omega_{\psi^{-1}, \mu}(hvtk)  \phi \right)_r
\left(\xi_\ell \begin{pmatrix} a_1\\ z_1 \end{pmatrix}, \xi_\ell \begin{pmatrix} a_2^\ast \\ -w({}^{t}z_2)w \end{pmatrix} \right) 
\\
\mu(\det a_1 a_2) \, dr \, dh \, dm \, da_1 \, da_2 \, dv \, dt \, dk
\end{multline}
where we write 
\[
m = \begin{pmatrix} I_r&&&\\ z_1&I_{\ell-r}&&\\ &&I_{\ell-r}&\\ &&z_2&I_r\end{pmatrix}.
\]
For fixed $g \in \mathrm{GL}_{2\ell}$, $\phi^\prime \in \mathcal{S}(F^r \times F^r)$ and $b_1, b_2 \in \mathrm{GL}_\ell$,
put
\begin{multline*}
f_{s}^{\phi^\prime}(g, b_1, b_2)
\\
= \int_{F^r} \int_{F^r} 
f_s\left( g, b_1 \begin{pmatrix}I_r&0&w_r{}^{t}x_1&0\\ &I_{\ell-r}&&\\ &&1&\\ &&&I_{n-\ell-1} \end{pmatrix}, 
b_2 \begin{pmatrix}I_{n-\ell-1}&&&\\ &1&&x_2\\ &&I_{\ell-r}&\\ &&&I_{r} \end{pmatrix} \right)
\\
\widehat{\phi^\prime}(x_1, x_2) \, dx_1 \, dx_2.
\end{multline*}
Then we find that 
\begin{multline*}
f_{s}^{\phi^\prime} \left(g, \left( \begin{smallmatrix}a_1&&\\z_1&I_{\ell-r}&\\ &&I_{n-\ell} \end{smallmatrix}\right), \left( \begin{smallmatrix}I_{n-\ell}&&\\&I_{\ell-r}&\\ &a_2 z_2&a_2 \end{smallmatrix}\right) \right)
\\
=
 \int_{F^r} \int_{F^r} \psi \left(- (z_{1})_{\ell-r} w_r {}^{t}x_1 + x_2w_r {}^{t}(z_2)^1 \right)
f_s\left( g,  \left( \begin{smallmatrix}a_1&&\\z_1&I_{\ell-r}&\\ &&I_{n-\ell} \end{smallmatrix}\right), 
\left( \begin{smallmatrix}I_{n-\ell}&&\\&I_{\ell-r}&\\ &a_2 z_2&a_2 \end{smallmatrix}\right)  \right)
\\
\widehat{\phi^\prime}(x_1, x_2) \, dx_1 \, dx_2
\end{multline*}
where $(z_1)_{\ell-r}$ and $(z_2)^1$ denote the $(\ell-r)$-th row of $z_1$ and the first column of $z_2$.
By the Fourier inversion formula, this equals to
\begin{multline*}
f_s\left( g,  \left( \begin{smallmatrix}a_1&&\\z_1&I_{\ell-r}&\\ &&I_{n-\ell} \end{smallmatrix}\right), 
\left( \begin{smallmatrix}I_{n-\ell}&&\\&I_{\ell-r}&\\ &a_2 z_2&a_2 \end{smallmatrix}\right)  \right)
\phi^\prime((z_{1})_{\ell-r} w_{\ell-r}, (z_2)_\ell) 
\\
=
f_s\left( g,  \left( \begin{smallmatrix}a_1&&\\z_1&I_{\ell-r}&\\ &&I_{n-\ell} \end{smallmatrix}\right), 
\left( \begin{smallmatrix}I_{n-\ell}&&\\&I_{\ell-r}&\\ &a_2 z_2&a_2 \end{smallmatrix}\right)  \right)
\phi^\prime\left(\xi_\ell \begin{pmatrix} a_1\\ z_1 \end{pmatrix}, \xi_\ell \begin{pmatrix} a_2^\ast \\ -w({}^{t}z_2)w \end{pmatrix} \right) .
\end{multline*}
Hence, \eqref{p.49} equals to
\begin{multline*}
 \int_{K_{\mathrm{GL}_{2\ell}}} \int_{T_{2(\ell-r)}} \int_{\overline{U_r} \cap \overline{U_\ell}}
 \int_{Y_\ell \backslash H_\ell} \int_{R_{\ell, n}}  \int_{(U_{\mathrm{GL}_r} \backslash \mathrm{GL}_r)^2}
\int_{\overline{Z_{r, \ell-r}}} 
\varphi(vtk, a_1, I_{2(\ell-r)}, a_2)  
\\
|\det a_1 a_2^{-1}|^{\frac{1}{2}(r-n)+s}
 f_s^{\left( \omega_{\psi^{-1}, \mu}(hvtk)  \phi \right)_r}(\beta_{\ell,n} rhvtk, \left( \begin{smallmatrix}a_1&&\\z_1&I_{\ell-r}&\\ &&I_{n-\ell} \end{smallmatrix}\right), 
\left( \begin{smallmatrix}I_{n-\ell}&&\\&I_{\ell-r}&\\ &a_2 z_2&a_2 \end{smallmatrix}\right)) 
\\
 \mu(\det a_1 a_2) \, dz_1 \, dz_2 \, da_1 \, da_2\, dr \, dh \, dv \, dt \, dk
 \\
 =
  \int_{K_{\mathrm{GL}_{2\ell}}} \int_{T_{2(\ell-r)}} \int_{\overline{U_r} \cap \overline{U_\ell}}
 \int_{Y_\ell \backslash H_\ell} \int_{R_{\ell, n}}  \int_{(U_{\mathrm{GL}_r} \backslash \mathrm{GL}_r)^2}
\int_{(M_{\ell-r \times r})^2} 
\varphi^\ast(vtk, a_1, I_{2(\ell-r)}, a_2)  
\\
|\det a_1 a_2|^{\frac{1}{2}(r-n)+s}
 f_{\tau, s, \ast}^{\left( \omega_{\psi^{-1}, \mu}(hvtk)  \phi \right)_r}(\beta_{\ell,n} rhvtk, \left( \begin{smallmatrix}a_1&&\\z_1&I_{\ell-r}&\\ &&I_{n-\ell} \end{smallmatrix}\right), 
\left( \begin{smallmatrix}a_2&&\\z_2&I_{\ell-r}&\\ &&I_{n-\ell} \end{smallmatrix}\right)) 
\\ \mu(\det a_1 a_2^{-1})
 \, dz_1 \, dz_2 \, da_1 \, da_2\, dr \, dh \, dv \, dt \, dk
\end{multline*}
where
\[
f_{\tau, s, \ast}(g, h_1, h_2) = f_s(g, h_1, h_2^\ast)
\quad
\text{and}
\quad
\varphi^\ast(x, y_1, z, y_2) = \varphi(x, y_1, z, y_2^\ast). 
\]
The inner $dz_1 \, da_1$-integration and $dz_2 \, da_2$-integration comprise for an integral 
for $\mathrm{GL}_\ell \times \mathrm{GL}_n$
and $(\sigma_1 \otimes \mu) \times \tau_1$ and $(\sigma_2 \otimes \mu)^\ast \times \tau_2$.
Now we apply \eqref{9.1} with $j=\ell-r$.
Multiplying the last integral by $\mu(-1)^\ell \omega_{\sigma_1}(-1)^{n-1} \gamma(s, \sigma_1 \otimes \mu \times \tau_1, \psi)$
and $\mu(-1)^\ell \omega_{\sigma_2}(-1)^{n-1} \gamma(s, (\sigma_2 \otimes \mu)^\ast \times \tau_2, \psi)$
leads to
\begin{multline*}
  \int_{K_{\mathrm{GL}_{2\ell}}} \int_{T_{2(\ell-r)}} \int_{\overline{U_r} \cap \overline{U_\ell}}
 \int_{Y_\ell \backslash H_\ell} \int_{R_{\ell, n}}  \int_{(U_{\mathrm{GL}_r} \backslash \mathrm{GL}_r)^2}
\int_{(M_{r \times n-\ell-1})^2 } 
\varphi^\ast(vtk, a_1, I_{2(\ell-r)}, a_2) 
\\
  f_{\tau, s, \ast}^{\left( \omega_{\psi^{-1}, \mu}(hvtk)  \phi \right)_r} \left(\beta_{\ell,n} rhvtk,\Omega_{n-r, r} \left( \begin{smallmatrix}a_1&&z_1\\&I_{\ell-r+1}&\\ &&I_{n-\ell-1} \end{smallmatrix}\right), 
\Omega_{n-r, r}\left( \begin{smallmatrix}a_2&&z_2\\&I_{\ell-r+1}&\\ &&I_{n-\ell-1} \end{smallmatrix}\right) \right) 
\\
 |\det a_1 a_2|^{-\frac{1}{2}(r+n)+s+\ell}
\mu(\det a_1 a_2^{-1}) \, dz_1 \, dz_2 \, da_1 \, da_2\, dr \, dh \, dv \, dt \, dk
\end{multline*}
where we recall
\[
\Omega_{n_1, n_2} = \begin{pmatrix} &I_{n_1}\\ I_{n_2}& \end{pmatrix}.
\]
Shifting $a_1$ and $a_2$ back yields
\begin{multline*}
  \int_{K_{\mathrm{GL}_{2\ell}}} \int_{T_{2(\ell-r)}} \int_{\overline{U_r} \cap \overline{U_\ell}}
 \int_{Y_\ell \backslash H_\ell} \int_{(U_{\mathrm{GL}_r} \backslash \mathrm{GL}_r)^2} \int_{R_{\ell, n}} 
\int_{(M_{r \times n-\ell-1})^2 } \int_{(F^r)^2} \int_{(F^r)^2}
\\
 \psi \left(- x_1 w_r {}^{t}y_1 + y_2w_r {}^{t}x_2 \right) |\det a_1 a_2|^{\ell-r} 
\varphi(\mathrm{diag}(a_1, I_{2(\ell-r)}, a_2^\ast)vtk, I_r, I_{2(\ell-r)}, I_r) 
\\
f_s \left(\beta_{\ell,n} rh\mathrm{diag}(a_1, I_{2(\ell-r)}, a_2^\ast)vtk,\Omega_{n-r, r} 
\left( \begin{smallmatrix}I_{\ell+1}&\begin{matrix}z_1\\ 0\end{matrix}\\ &I_{n-\ell-1}\end{smallmatrix}\right)y_1, 
\Omega_{r, n-r}\left(\begin{smallmatrix}I_{n-\ell-1}&\begin{matrix}0&{}^{t}z_2 \end{matrix}\\ &I_{\ell+1} \end{smallmatrix} \right)y_2 \right) 
\\
\omega_{\psi^{-1}, \mu}((x_1, x_2)h\mathrm{diag}(a_1, I_{2(\ell-r)}, a_2^\ast)vtk)  \phi (\xi_\ell^\prime, \xi_\ell^\prime)
\\
dx_1 \, dx_2 \, dy_1 \, dy_2 \, dz_1 \, dz_2 \, da_1 \, da_2\, dr \, dh \, dv \, dt \, dk
\end{multline*}
where for $x_i = (x_{i1}, \dots, x_{ir}) \in F^r$, $(x_1, x_2)$ denotes 
$(((x_1,0, \dots, 0), (x_2, 0, \dots, 0)), 0; 0) \in X_\ell$,
\[
\xi_\ell^\prime
=
\left\{
\begin{array}{ll}
\xi_\ell, & r<\ell \\
0,& r=\ell
\end{array}
\right.
\quad
y_1= \begin{pmatrix}I_r&0&w_r{}^{t}y_1&0\\ &I_{\ell-r}&&\\ &&1&\\ &&&I_{n-\ell-1} \end{pmatrix}, \,
y_2 = \begin{pmatrix}I_{n-\ell-1}&&&\\ &1&&y_2\\ &&I_{\ell-r}&\\ &&&I_{r} \end{pmatrix}.
\]
For simplicity, we write 
\[
z_1^\sharp = \left( \begin{smallmatrix}I_{\ell+1}&\begin{matrix}z_1\\ 0\end{matrix}\\ &I_{n-\ell-1}\end{smallmatrix}\right)
\quad
\text{and}
\quad
z_2^\flat =\left(\begin{smallmatrix}I_{n-\ell-1}&\begin{matrix}0&{}^{t}z_2 \end{matrix}\\ &I_{\ell+1} \end{smallmatrix} \right).
\]
Since the function 
\begin{multline*}
g \mapsto \int_{R_{\ell, n}} 
\int_{(M_{r \times n-\ell-1})^2 } \int_{(F^r)^2} \int_{(F^r)^2}
\varphi(g, I_r, I_{2(\ell-r)}, I_r) 
\omega_{\psi^{-1}, \mu}((x_1, x_2)hg)  \phi (\xi_\ell^\prime, \xi_\ell^\prime)
\\
f_s \left(\beta_{\ell,n} rhg,\Omega_{n-r, r} 
z_1^\sharp y_1, 
\Omega_{r, n-r} z_2^\flat y_2 \right) 
 \psi \left(- x_1 w_r {}^{t}y_1 + y_2w_r {}^{t}x_2 \right)
dx_1 \, dx_2 \, dy_1 \, dy_2 \, dz_1 \, dz_2 \, dr 
\end{multline*}
is $\overline{Z_{r, \ell-r}} U_{\mathrm{GL}_r} U_{\mathrm{GL}_{2(\ell-r)}}U_{\mathrm{GL}_r }$-invariant on the left, 
we obtain
\begin{multline*}
  \int_{(Z_{r, \ell-r} U_{\mathrm{GL}_r} U_{\mathrm{GL}_{2(\ell-r)}}U_{\mathrm{GL}_r }) \backslash \mathrm{GL}_{2\ell}}
 \int_{Y_\ell \backslash H_\ell} \int_{R_{\ell, n}} 
\int_{(M_{r \times n-\ell-1})^2 } \int_{(F^r)^2} \int_{(F^r)^2}
\varphi(g, I_r, I_{2(\ell-r)}, I_r) 
\\
\omega_{\psi^{-1}, \mu}((x_1, x_2)hg)  \phi (\xi_\ell^\prime, \xi_\ell^\prime)
f_s \left(\beta_{\ell,n} rhg,\Omega_{n-r, r} 
z_1^\sharp y_1, 
\Omega_{r, n-r} z_2^\flat y_2 \right) 
\\
 \psi \left(- x_1 w_r {}^{t}y_1 + y_2w_r {}^{t}x_2 \right) 
dx_1 \, dx_2 \, dy_1 \, dy_2 \, dz_1 \, dz_2 \,  dr \, dh \, dg.
\end{multline*}
After factoring the integral through $\overline{U_r} \cap \overline{U_\ell}$,
we obtain
\begin{multline*}
  \int_{ (U_{\mathrm{GL}_r})^2 \overline{U_r} \backslash \mathrm{GL}_{2\ell}} \int_{\overline{U_r} \cap \overline{U_\ell}}
   \int_{Y_\ell \backslash H_\ell} \int_{R_{\ell, n}} 
\int_{(M_{r \times n-\ell-1})^2 } \int_{(F^r)^2} \int_{(F^r)^2}
\\
\varphi(g, I_r, I_{2(\ell-r)}, I_r) 
\omega_{\psi^{-1}, \mu}((x_1, x_2)((0,0), (0_{\ell-r}, 2y_2w_r,0_{\ell-r}, 2y_1))hug)  \phi (\xi_\ell^\prime, \xi_\ell^\prime)
\\
f_s \left(\beta_{\ell,n} rhug,\Omega_{n-r, r}  z_1^\sharp y_1, 
\Omega_{r, n-r} z_2^\flat y_2 \right) 
dx_1 \, dx_2 \, dy_1 \, dy_2 \, dz_1 \, dz_2 \, dr \, dh \, du \, dg.
\end{multline*}
We see that ${}^{\beta_{\ell, n}}( \mathrm{diag}(z_1^\sharp y_1, z_2^\flat y_2))$ normalizes $R_{\ell, n}$.
Also, if $h = ((c_{11}, c_{12}, c_{21}, c_{22}), (0,0);c)$ where $c_{i1} \in F^r$
$c_{i2} \in F^{\ell-r}$ and $c \in F$, we have
\[
({}^{\beta_{\ell, n}}( \mathrm{diag}(z_1^\sharp y_1, z_2^\flat y_2)))^{-1} h \, 
{}^{\beta_{\ell, n}}( \mathrm{diag}(z_1^\sharp y_1, z_2^\flat y_2))= r^\prime h (0_{2\ell}, 0_{2\ell}; -2c_{11} w_r{}^{t}y_1+2y_2{}^{t}c_{21})
\]
with $r^\prime \in R_{\ell, n}$.
Therefore, we can change the order 
${}^{\beta_{\ell, n}}( \mathrm{diag}(z_1^\sharp y_1, z_2^\flat y_2))h \mapsto h {}^{\beta_{\ell, n}}( \mathrm{diag}(z_1^\sharp y_1, z_2^\flat y_2))$, if we change variables in $r$ and $h$.
Additionally, if 
\[
u= \begin{pmatrix}I_r&&&\\ &I_{\ell-r}&&\\ u_1&&I_{\ell-r}&\\ u_3&u_2&&I_r \end{pmatrix},
\]
which is regarded as an element of $\mathrm{GL}_{2n}$,
we have
\[
({}^{\beta_{\ell, n}}( \mathrm{diag}(z_1^\sharp y_1, z_2^\flat y_2)))^{-1} u \, 
{}^{\beta_{\ell, n}}( \mathrm{diag}(z_1^\sharp y_1, z_2^\flat y_2)) 
= r_{u, m} uA(u_1, u_2, u_3; y_1, y_2)
\]
with $r_{u, m} \in R_{\ell, n}$ and 
\begin{multline*}
A(u_1, u_2, u_3; y_1, y_2) 
\\
= (-2y_2 w_r u_3, -2y_2 w_r u_2, 2y_1w_r {}^{t}u_3 w_r, 2y_1 w_r {}^{t}u_1 w_r, 0 \dots, 0, 4y_2w_ru_3w_r{}^{t}y_1) \in H_\ell.
\end{multline*}
Since $h$ normalizes $R_{\ell, n}$, we can change again variables in $r$ to remove $r_{u, m}$.
Also note that $u$ commute with $Y_\ell \backslash H_\ell$. Thus we have
\begin{multline*}
  \int_{ (U_{\mathrm{GL}_r})^2 \overline{U_r} \backslash \mathrm{GL}_{2\ell}} \int_{\overline{U_r} \cap \overline{U_\ell}}
\int_F \int_{(F^r)^2}   \int_{(F^{\ell-r})^2} \int_{R_{\ell, n}} 
\int_{(M_{r \times n-\ell-1})^2 } \int_{(F^r)^2} \int_{(F^r)^2}
\\
\varphi(g, I_r, I_{2(\ell-r)}, I_r) 
\omega_{\psi^{-1}, \mu}((x_1, x_2)((0,0), (0_{\ell-r}, 2y_2w_r,0_{\ell-r}, 2y_1))
\\
u( (c_{11}, c_{12}, c_{21}, c_{22}), (0, 0); c+2c_{11} w_r{}^{t}y_1-2y_2{}^{t}c_{21}  )g)  \phi (\xi_\ell^\prime, \xi_\ell^\prime)
\\
f_s \left(\beta_{\ell,n} ru( (c_{11}, c_{12}, c_{21}, c_{22}), (0, 0); c) A(u_1, u_2, u_3; y_1, y_2) \right.
\\
\,  \left. {}^{\beta_{\ell, n}}( \mathrm{diag}(z_1^\sharp y_1, z_2^\flat y_2)) g,\Omega_{n-r, r} , 
\Omega_{r, n-r}\right) 
\\
dx_1 \, dx_2 \, dy_1 \, dy_2 \, dz_1 \, dz_2 \, dr \, dc_{11} \, dc_{21} \, dc_{12} \, dc_{22} \, dc \, du \, dg.
\end{multline*}
Further, we note that 
\begin{multline*}
((0,0), (0_{\ell-r}, 2y_2w_r,0_{\ell-r}, 2y_1))u 
= u((0,0), (0_{\ell-r}, 2y_2w_r,0_{\ell-r}, 2y_1))
\\
 \cdot (-2y_2 w_r u_3, -2y_2 w_r u_2, 2y_1w_r {}^{t}u_3 w_r, 2y_1 w_r {}^{t}u_1 w_r, 0 \dots, 0, -4y_2w_ru_3w_r{}^{t}y_1).
\end{multline*}
On the other hand, we have
\begin{multline*}
((0,0), (0_{\ell-r}, 2y_2w_r,0_{\ell-r}, 2y_1))
(c_{11}, 0, c_{21}, 0, 0_{2\ell}; -2c_{11} w_r{}^{t}y_1+2y_2{}^{t}c_{21})
\\
\cdot ((0,0), (0_{\ell-r}, -2y_2w_r,0_{\ell-r}, -2y_1))
=
(c_{11}, 0, c_{21}, 0, 0_{2\ell}; 0).
\end{multline*}
Then changing variables, the above integral is equal to
\begin{multline*}
  \int_{ (U_{\mathrm{GL}_r})^2 \overline{U_r} \backslash \mathrm{GL}_{2\ell}} \int_{\overline{U_r} \cap \overline{U_\ell}}
\int_F \int_{(F^r)^2}   \int_{(F^{\ell-r})^2}  \int_{R_{\ell, n}} 
\int_{(M_{r \times n-\ell-1})^2 } \int_{(F^r)^2} \int_{(F^r)^2}\varphi(g, I_r, I_{2(\ell-r)}, I_r) 
\\
\omega_{\psi^{-1}, \mu}( u (x_1, x_2)((0,0), (0_{\ell-r}, 2y_2w_r,0_{\ell-r}, 2y_1))( (0, 0, c_{21}, c_{22}), (0, 0); c )g)  \phi (\xi_\ell^\prime, \xi_\ell^\prime)
\\
f_s \left(\beta_{\ell,n} ru( (c_{11}, c_{12}, c_{21}, c_{22}), (0, 0); c)  \,  {}^{\beta_{\ell, n}}( \mathrm{diag}(z_1^\sharp y_1, z_2^\flat y_2)) g,\Omega_{n-r, r} , 
\Omega_{r, n-r}\right) 
\\
dx_1 \, dx_2 \, dy_1 \, dy_2 \, dz_1 \, dz_2 \, dr \, dc_{11} \, dc_{21} \, dc_{12} \, dc_{22} \, dc \, du \, dg.
\end{multline*}
In the same argument as \cite[p.441]{Ka}, we can prove 
\begin{multline*}
\int_{(F^r)^2}
\omega_{\psi^{-1}, \mu}(u((x_1, x_2), (0, 0); 0))\phi(\xi_\ell^\prime, \xi_\ell^\prime) \, dx_1 \, dx_2
\\
=
\int_{(F^r)^2}
\omega_{\psi^{-1}, \mu}((x_1, x_2), (0, 0); 0))\phi(\xi_\ell^\prime, \xi_\ell^\prime) \, dx_1 \, dx_2.
\end{multline*}
Applying this identity to our local integral, we get
\begin{multline*}
  \int_{ (U_{\mathrm{GL}_r})^2 \overline{U_r} \backslash \mathrm{GL}_{2\ell}} \int_{\overline{U_r} \cap \overline{U_\ell}}
\int_F \int_{(F^r)^2}   \int_{(F^{\ell-r})^2}  \int_{R_{\ell, n}} 
\int_{(M_{r \times n-\ell-1})^2 } \int_{(F^r)^2} \int_{(F^r)^2}
\\
\omega_{\psi^{-1}, \mu}(  (x_1, x_2)((0,0), (0_{\ell-r}, 2y_2w_r,0_{\ell-r}, 2y_1))( (0, 0, c_{21}, c_{22}), (0, 0); c )g)  \phi (\xi_\ell^\prime, \xi_\ell^\prime)
\\
f_s \left(\beta_{\ell,n} ru( (c_{11}, c_{12}, c_{21}, c_{22}), (0, 0); c)  \,  {}^{\beta_{\ell, n}}( \mathrm{diag}(z_1^\sharp y_1, z_2^\flat y_2)) g,\Omega_{n-r, r} , 
\Omega_{r, n-r}\right) 
\\
\varphi(g, I_r, I_{2(\ell-r)}, I_r) \, dx_1 \, dx_2 \, dy_1 \, dy_2 \, dz_1 \, dz_2 \, dr \,dc_{11} \, dc_{21} \, dc_{12} \, dc_{22} \, dc \, du \, dg.
\end{multline*}
Let $w = \beta_{\ell-r, n}^{-1} \omega_{n-r, r} \beta_{\ell, n}$.
Let $X_\ell^\circ$ be the subgroup of $X_\ell$ consisting of elements $\{ ( (c_{11}, 0_{\ell-r}, c_{21}, 0_{\ell-r}), (0, 0); 0) \} $.
Then $X_\ell = X_\ell^\circ \times X_{\ell-r}$ where $X_{\ell-r}$ is the subgroup of elements $\{ ( (0_r, c_{12}, 0_r, c_{22}), (0, 0); 0)\}$,
and we have $Y_\ell \backslash H_\ell \simeq X_\ell^\circ \times X_{\ell-r} \times F$, where $F = \{ (0_\ell, 0_\ell, (0, 0); c) \}$.
Also $Y_{\ell-r} \backslash H_{\ell-r} \simeq X_{\ell-r} \times F$
and when $X_{\ell-r} \times F$ is regarded as a subgroup of $\mathrm{GL}_{2\ell+2}$, ${}^{w^{-1}} (X_{\ell-r} \times F)$
is equal to the image of $Y_{\ell-r} \backslash H_{\ell-r}$ as a subgroup of $\mathrm{GL}_{2(\ell-r+1)}$.
Moreover, ${}^{w^{-1}} (R_{\ell, n} (\overline{U_r} \cap \overline{U_\ell}) X_\ell^\circ) = R_{\ell-r, n}$.
Factoring $dg$-integration through $U_{\mathrm{GL}_{2(\ell-r)}} \backslash \mathrm{GL}_{2(\ell-r)}$, the above integral is equal to
\begin{multline}
\label{9.3}
  \int_{ (U_{\mathrm{GL}_r})^2 \overline{U_r} \mathrm{GL}_{2(\ell-r)} \backslash \mathrm{GL}_{2\ell}}  \int_{U_{\mathrm{GL}_{2(\ell-r)}} \backslash \mathrm{GL}_{2(\ell-r)}}
 \int_{Y_{\ell-r} \backslash H_{\ell-r}}  \int_{R_{\ell-r, n}} 
\int_{(M_{r \times n-\ell-1})^2 } \int_{(F^r)^2} \int_{(F^r)^2}
\\
\omega_{\psi^{-1}, \mu}( hg^\prime (x_1, x_2)((0,0), (0_{\ell-r}, 2y_2w_r,0_{\ell-r}, 2y_1))g)  \phi (\xi_\ell^\prime, \xi_\ell^\prime)
\\
f_s \left(\beta_{\ell-r,n} rhg^\prime w  \,  {}^{\beta_{\ell, n}}( \mathrm{diag}(z_1^\sharp y_1, z_2^\flat y_2)) g, 1_n , 
1_n\right) 
\\
\varphi(g, I_r, g^\prime, I_r)  \, dx_1 \, dx_2 \, dy_1 \, dy_2 \, dz_1 \, dz_2 \, dr \,dh \, dg^\prime  \, dg.
\end{multline}
Assume $r < \ell$. Then we can apply the functional equation for $\pi^\prime \times \tau$ to $dr \,dh \,dg^\prime$-integration and get that 
\eqref{9.3} multiplied by $\gamma(s, \pi^\prime \times \tau, \psi)$ equals
\begin{multline}
\label{9.4}
  \int_{ (U_{\mathrm{GL}_r})^2 \overline{U_r} \mathrm{GL}_{2(\ell-r)} \backslash \mathrm{GL}_{2\ell}}  \int_{U_{\mathrm{GL}_{2(\ell-r)}} \backslash \mathrm{GL}_{2(\ell-r)}}
 \int_{Y_{\ell-r} \backslash H_{\ell-r}}  \int_{R_{\ell-r, n}} 
\int_{(M_{r \times n-\ell-1})^2 } \int_{(F^r)^2} \int_{(F^r)^2}
\\
\varphi(g, I_r, g^\prime, I_r) \omega_{\psi^{-1}, \mu}( hg^\prime (x_1, x_2)((0,0), (0_{\ell-r}, 2y_2w_r,0_{\ell-r}, 2y_1))g)  \phi (\xi_\ell^\prime, \xi_\ell^\prime)
\\
M^\ast(\tau, s)f_s \left(\beta_{\ell-r,n} rhg^\prime w  \,  {}^{\beta_{\ell, n}}( \mathrm{diag}(z_1^\sharp y_1, z_2^\flat y_2)) g, 1_n , 
1_n\right) 
\\
dx_1 \, dx_2 \, dy_1 \, dy_2 \, dz_1 \, dz_2 \, dr \,dh \, dg^\prime  \, dg.
\end{multline}
Thus, $\mathcal{L}(W_\varphi, f_s, \phi)$ multiplied by 
\begin{multline*}
\mu(-1)^{n-1} \omega_{\sigma_1}(-1)^{n-1} \mu(-1)^{n-1} \omega_{\sigma_2}(-1)^{n-1}
\gamma(s, \sigma_1 \otimes \mu \times \tau_1, \psi)  
\\
\gamma(s, (\sigma_2 \otimes \mu)^\ast \times \tau_2, \psi) \gamma(s, \pi^\prime \times \tau, \psi)
\end{multline*}
equals \eqref{9.4}.
Reversing the arguments above with $(\tau, s, f_s)$ replaced by $(\tau^\ast, 1-s, M^\ast(\tau, s)f_s)$, we see that 
\eqref{9.4} multiplied by $\mu(-1)^{n-1} \omega_{\sigma_1}(-1)^{n-1}\gamma(1-s, \sigma_1 \otimes \mu \times \tau_2^\ast, \psi)^{-1} \mu(-1)^{n-1} \omega_{\sigma_2}(-1)^{n-1}
\gamma(1-s, (\sigma_2 \otimes \mu)^\ast  \times \tau_1^\ast, \psi)^{-1}$
equals $\mathcal{L}^\ast(W_\varphi, f_s, \phi)$. 
Then by \eqref{fq gamma}, the required identity readily follows
%\[
%\gamma(s, \pi \times \tau, \psi, \mu) = \gamma(s, (\sigma_1 \otimes \mu)^\ast \times \tau_2, \psi)
%\gamma(s, \varepsilon_2 \otimes \mu  \times \tau_1, \psi)
%\gamma(s, \varepsilon_1 \otimes \mu \times \tau_1, \psi)  \gamma(s, (\varepsilon_2 \otimes \mu)^\ast \times \tau_2, \psi) \gamma(s, \pi^\prime \times \tau, \psi, \mu)
%\]
%with $\tau = (\tau_1, \tau_2)$.
%

If $r = \ell$, the integral \eqref{9.3} is
\begin{multline*}
  \int_{ (U_{\mathrm{GL}_r})^2 \overline{U_r} \mathrm{GL}_{2(\ell-r)} \backslash \mathrm{GL}_{2\ell}} 
\int_{(M_{r \times n-\ell-1})^2 } \int_{(F^r)^2} \int_{(F^r)^2}
\\
\varphi(g, I_\ell, I_\ell)
\omega_{\psi^{-1}, \mu}(  (x_1, x_2)((0,0), (0_{\ell-r}, 2y_2w_r,0_{\ell-r}, 2y_1))g)  \phi (\xi_\ell^\prime, \xi_\ell^\prime)
\\
\left( \int_{U_n} f_s \left(\beta_{\ell-r,n} uw  \,  {}^{\beta_{\ell, n}}( \mathrm{diag}(z_1^\sharp y_1, z_2^\flat y_2)) g, 1_n , 
1_n\right) \psi^{-1} (u_{n, n+1})\, du \right)
\\
dx_1 \, dx_2 \, dy_1 \, dy_2 \, dz_1 \, dz_2  \, dg.
\end{multline*}
By \eqref{3.7 sp}, this integral is equal to \eqref{9.4}
and we can prove the required multiplicativity (without $\gamma(s, \pi^\prime \times \tau, \mu,  \psi)$).
\begin{Remark}
\label{conv rem}
When $r=\ell$, the above integral is equal to our local zeta integral in the case of $\mathrm{GL}_{2n} \times (\mathrm{GL}_\ell \times \mathrm{GL}_\ell)$. 
Hence, the meromorphic continuation of the local zeta integrals for the case of $\ell < n$ is proved as in \cite[Section~5, Section~6]{S2}.
This is an analogue of symplectic or metaplectic case remarked in Kaplan~\cite[Remark~9.1]{Ka}.
\end{Remark}
\begin{Remark}
\label{qs rem1}
We would like to mention the non-split case. 
For example, when $r<\ell$, in a similar argument as a proof of \cite[(9.4)]{Ka}, we may show that $\mathcal{L}(W, f_{\varphi_s}, \phi) $
is equal to the following analogue of \eqref{9.4} (up to suitable local $\gamma$-factors):
\begin{multline*}
\int_{\overline{U}_r U_{\mathrm{GL}_r} G_{\ell-r} \backslash G_\ell} \int_{M_{r \times (n-\ell-1)}} \int_{E^r} \int_{E^r}
 \int_{U_{G_{\ell-r}} \backslash G_{\ell-r}} \int_{Y_{\ell-r} \backslash H_{\ell-r}} \int_{R_{\ell-r, n}}
 \\
 \varphi(g, I_r, g^\prime) \omega_{\psi^{-1}, \Upsilon}(hg^\prime(x, 0\dots, 0)(0,\cdots, 0, z, 0))\phi(\xi_\ell^\prime)
 \\
 M^\ast(\tau, s) f_s(\beta_{\ell-r, n} rhg^\prime w ({}^{\beta_{\ell, n}} (mz) ) g, I_n) \,dr \, dh \, dg^\prime \, dx \, dz \, dm \, dg.
\end{multline*}
\end{Remark}
\subsubsection{$r=\ell > n$}
\label{9.2}
Recall that our local integral is 
\[
 \int_{U_{\mathrm{GL}_{2n}} \backslash \mathrm{GL}_{2n}}  \int_{R^{\ell, n}}\int_{X_n}
W_\varphi({}^{\omega_{\ell-n, n}}(rxg)) f_s(g, I_n, I_n)  \omega_{\psi^{-1}, \mu}(g) \phi(x) \, dx \, dr \, dg .
\]
Write the $dg$-integration in this integral over $\overline{U_n} \times \left( U_{\mathrm{GL}_n} \backslash \mathrm{GL}_n \right)^2$.
For $a_1, a_2 \in \mathrm{GL}_n$, $\mathrm{diag}(a_1, a_2)$ normalizes $X_n$ and $R^{\ell, n}$. Then we see that it changes 
$dx_1 dx_2 \rightarrow |\det a_1|^{-1} |\det a_2| dx_1 dx_2$ and $dr \rightarrow |\det a_1|^{-\ell+n+1} |\det a_2|^{\ell-n-1} dr$.
Further, ${}^{\omega_{\ell-n ,n}}\mathrm{diag}(a_1, a_2)$ normalizes $U_\ell$ and fixes $\psi(u)$.
Then it changes $du \mapsto |\det a_1|^{\ell} |\det a_2|^{-\ell} du$, and the integration formula introduces a factor of $|\det a_1|^{-n} |\det a_2|^n$.
Additionally, 
\begin{align*}
&f_s(\mathrm{diag}(a_1, a_2)g, I_n, I_n) = |\det a_1 a_2^{-1}|^{\frac{1}{2}n+s-\frac{1}{2}}, \\
&\varphi(\mathrm{diag}(a_1, a_2)g, I_\ell, I_\ell) = |\det a_1 a_2^{-1}|^{-\frac{1}{2}\ell} \varphi(g, a),\\
&\omega_{\psi^{-1}, \mu}(axv)\phi(0) = \mu(\det(a_1 a_2)) |\det (a_1 a_2^{-1})|^{\frac{1}{2}}\omega_{\psi^{-1}, \mu}(v)\phi(x)
\end{align*}
and ${}^{\omega_{\ell-n, n}} r$ normalizes $U_\ell$ and $\psi$.
Thus, the integral becomes
\begin{multline*}
 \int_{\overline{U_n}}  \int_{X_n} \int_{U_\ell} \int_{(U_{\mathrm{GL}_n} \backslash \mathrm{GL}_n)^2} \int_{R^{\ell, n}}
\varphi \left( u\, {}^{\omega_{\ell-n, n}}(xv), \left(\begin{smallmatrix}a_1&&\\ r_1&I_{\ell-n-1}&\\ &&1 \end{smallmatrix} \right), 
 \left(\begin{smallmatrix}1&&\\ &I_{\ell-n-1}&\\ &a_2 r_2&a_2 \end{smallmatrix} \right) \right) \psi(u)\\
\mu(\det a_1 a_2) |\det a_1 a_2^{-1}|^{\frac{1}{2}(n-\ell)+s} f_s(v, a_1, a_2)  \omega_{\psi^{-1}, \mu}(v) \phi(x) \, dx \, dr \, dg .
\end{multline*}
As in the previous case, multiplying the last integral by 
$\omega_{\tau_1}(-1)^{\ell-1} \gamma(s, \sigma_1 \otimes \mu \times \tau_1, \psi^{-1}) \omega_{\tau_2}(-1)^{\ell-1} 
\gamma(s, (\sigma_2 \otimes \mu)^\ast \times \tau_1, \psi^{-1})$ leads to
\begin{multline*}
 \int_{\overline{U_n}}  \int_{X_n} \int_{U_\ell} \int_{(U_{\mathrm{GL}_n} \backslash \mathrm{GL}_n)^2} 
\varphi \left( u\, {}^{\omega_{\ell-n, n}}(xv), \Omega_{\ell-n,n} \left(\begin{smallmatrix}a_1&\\ &I_{\ell-n}\end{smallmatrix} \right), \Omega_{n, n-\ell} \begin{pmatrix}I_{\ell-n}&\\ &a_2 \end{pmatrix} \right) \\
 \psi(u) \mu(\det a_1 a_2) |\det a_1 a_2^{-1}|^{-\frac{1}{2}(n-\ell)+s-1} f_s(v, a_1, a_2)  \omega_{\psi^{-1}, \mu}(v) \phi(x) \, dx \, dr \, dg .
\end{multline*}
Reversing the previous steps, one arrives at
\begin{multline*}
 \int_{U_{\mathrm{GL}_{2n}} \backslash \mathrm{GL}_{2n}}  \int_{X_n} \int_{U_\ell} 
\varphi \left( u\, {}^{\omega_{\ell-n, n}}(xg), \Omega_{\ell-n,n}, \Omega_{n, n-\ell} \right) 
\\
\psi_{U_\ell}(u) f_s(g, I_n, I_n)  \omega_{\psi^{-1}, \mu}(g) \phi(x) \, dx \, dr \, dg .
\end{multline*}
Decompose $U_\ell = U_\ell^\prime \cdot {}^{\omega_{\ell-n, n}} U_n$ where 
$U_n$ is considered as the unipotent subgroup of $\mathrm{GL}_{2\ell}$ and 
\begin{equation}
\label{9.5}
U_\ell^\prime = \left\{ \begin{pmatrix} I_n&&v_1&\\ &I_{\ell-n}&v_2&v_3\\ &&I_{\ell-n}&\\ &&&I_n \end{pmatrix} \right\}.
\end{equation}
For $u \in U_n$, ${}^{u^{-1}} x = b_{u,x} x$ where $b_{u, x} \in H_n$ is such that ${}^{\omega_{\ell-n, n}} b_{u, x} \in U_\ell^\prime$ and 
\[
\omega_{\psi^{-1}, \mu}(xg)\phi(0)
=\omega_{\psi^{-1}, \mu}(uxg)\phi(0)
= \psi_{U_{\ell}}^{-1}({}^{\omega_{\ell-n,n}} b_{u,x})\omega_{\psi^{-1}, \mu}(xug)\phi(0).
\]
Further, we note that $f_s(ug, I_n) = f_s(g, I_n)$.
Hence, collapsing the integration over ${}^{\omega_{\ell-n, n}} U_n$ into the $dg$-integration, we get
\begin{multline*}
 \int_{Z_n \backslash \mathrm{GL}_{2n}}  \int_{X_n} \int_{U_\ell^\prime} 
\varphi \left( u\, {}^{\omega_{\ell-n, n}}(xg), \Omega_{\ell-n,n}, \Omega_{n, n-\ell} \right) 
\\
\psi_{U_\ell}(u) f_s(g, I_n, I_n)  \omega_{\psi^{-1}, \mu}(g) \phi(x) \, dx \, dr \, dg. 
\end{multline*}
Let us denote $u \in U_\ell^\prime$ in the form \eqref{9.5}.
We also write by $y_1$ the last row of $v_3$ and $y_2$ the first column of $v_1$
 and set $z =(v_2)_{\ell-n, 1}$.
 Then ${}^{\omega_{n, \ell-n}}u = {}^{\omega_n}r ((0, 0), y_1, {}^{t}y_2 w, z)$,
 where $r$ is a general element of $R_{n, \ell}$, $((0, 0), y_1, {}^{t}y_2 w, z) \in H_n$
 and recall $\omega_n =\left( \begin{smallmatrix}&I_n\\ -I_n \end{smallmatrix} \right)$
 is regarded as an element of $\mathrm{GL}_{2\ell}$ by the embedding 
 $g \mapsto \left( \begin{smallmatrix} 1&&\\ &g&\\ &&1\end{smallmatrix}\right)$. We obtain
 \begin{multline*}
 \int_{Z_n \backslash \mathrm{GL}_{2n}}  \int_{X_n} \int_{Y_n} \int_F \int_{R_{n,\ell}} 
\varphi \left( {}^{\omega_n} rxg \omega_{\ell-n,n}, I_\ell, I_\ell \right)  f_s(g, I_n, I_n)  
\\
\omega_{\psi^{-1}, \mu}((0, y, z)xg) \phi(0) \, dx \, dr \, dg .
\end{multline*}
Because $x \in X_n$ normalizes ${}^{\omega_n} R_{n, \ell}$, the integral becomes
 \begin{multline*}
 \int_{Z_n \backslash \mathrm{GL}_{2n}}  \int_{Y_n} \int_F \int_{R_{n,\ell}} 
\varphi \left( {}^{\omega_n} rg 
\omega_{\ell-n,n}, I_\ell, I_\ell \right) f_s(g, I_n, I_n) 
\\
( \omega_{\psi^{-1}, \mu}((0, y, z)xg))^\wedge \phi(0) \, dx \, dr \, dg .
\end{multline*}
Next factor the $dg$-integration through $\overline{U_n}$.
Let $\overline{u} =\left( \begin{smallmatrix} I_n&\\ u&I_n\end{smallmatrix} \right)$.
Observe that 
\begin{multline*}
{}^{\bar{u}}((0, 0), (y_1, y_2), z) = ((y_1u, y_2{}^{t}(\omega_n u \omega_n)), (y_1, y_2), z)
\\
=  ((y_1u, y_2{}^{t}(\omega_n u \omega_n)), (0, 0), 0)  ((0, 0), (y_1, y_2), z-(y_1 u)w{}^{t}y_2 )
\end{multline*}
and 
\begin{multline*}
( \omega_{\psi^{-1}, \mu}((0, y, z)\bar{u}g))^\wedge \phi(0)
\\
=  \omega_{\psi^{-1}, \mu}(\omega_n \bar{u} ((y_1u, y_2{}^{t}(\omega_n u \omega_n)), (y_1, y_2), z) g) \phi(0)
\\
=  \omega_{\psi^{-1}, \mu}( \left( \begin{smallmatrix}I_n&-u\\ &I_n \end{smallmatrix} \right)\omega_n ((y_1u, y_2{}^{t}(\omega_n u \omega_n)), (y_1, y_2), z) g) \phi(0)
\\
= \psi\left(u_{n,1} \right)\psi\left(\xi_n(u J_n{}^t y_2)+ \xi_n J_n {}^{t}(y_1u) \right)
\\
\omega_{\psi^{-1}, \mu}(((y_1, y_2),(0, 0), z-(y_1 u)w{}^{t}y_2 ) \omega_ng).
\end{multline*}
Also, $\varphi(\bar{u}, I_n, I_n) = \varphi(I_{2n}, I_n, I_n)$, and
for $r \in R_{n, \ell}$,  ${}^{\bar{u}} ({}^{\omega_n}r) = b_{\bar{u}, r} {}^{\omega_n}r$
with $\varphi(b_{\bar{u}, r}, I_n, I_n) = \varphi(I_{2n}, I_n, I_n)$.
Moreover, $((y_1u, y_2{}^{t}(\omega_n u \omega_n)), (0, 0), 0) $ normalizes ${}^{\omega_n} R_{n,\ell}$ and contributes $\psi^{-1}((y_1 u )_1-(u\omega_n {}^{t}y_2)_n)$
because of $\varphi$.
Combining these observations and changing $z \mapsto z+(y_1 u)w{}^{t}y_2 ) \omega_ng$ yield
 \begin{multline*}
 \int_{\overline{U_n}Z_n \backslash \mathrm{GL}_{2n}}  \int_{Y_n \backslash H_n} \int_{R_{n,\ell}} 
\varphi \left( {}^{\omega_n} rhg \omega_{\ell-n,n}, I_\ell, I_\ell \right) ( \omega_{\psi^{-1}, \mu}( h\omega_ng))\phi(\xi_n, \xi_n)
\\
 \left( \int_{U_n} f_s(\omega_n u \omega_n^{-1}g, I_n, I_n) \psi^{-1}\left(u_{n,1} \right) \, du \right) \, dr \, dh \, dg .
\end{multline*}
The inner $du$-integration resembles the left-hand side of \eqref{3.7 split} and we see that this integral is equal to
 \begin{multline}
 \label{9.7}
 \int_{\overline{U_n}Z_n \backslash \mathrm{GL}_{2n}}  \int_{Y_n \backslash H_n} \int_{R_{n,\ell}} 
\varphi \left( {}^{\omega_n} rhg \omega_{\ell-n,n}, I_\ell, I_\ell \right) ( \omega_{\psi^{-1}, \mu}( h\omega_ng))\phi(\xi_n, \xi_n)
\\
 \left( \int_{U_n} M^\ast(\tau, s)f_s(\omega_n u \omega_n^{-1}g, I_n, I_n) \psi^{-1}\left(u_{n,1} \right) \, du \right) \, dr \, dh \, dg .
\end{multline}
Summing up, we see that $\mathcal{L}(W_\varphi, f_s, \phi)$
multiplied by $\omega_{\tau_1}(-1)^{n-1} \gamma(s, \sigma_1 \otimes \mu \times \tau_1, \psi^{-1}) \omega_{\tau_2}(-1)^{n-1} 
\gamma(s, (\sigma_2 \otimes \mu)^\ast \times \tau_2, \psi^{-1})$ equals \eqref{9.7}.
Reversing the arguments above, \eqref{9.7} multiplied by 
\[
\omega_{\tau_1}(-1)^{n-1}\gamma(1-s, \sigma_1 \otimes \mu \times \tau_2^\ast, \psi)^{-1} \omega_{\tau_2}(-1)^{n-1}
\gamma(1-s, (\sigma_2 \otimes \mu)^\ast  \times \tau_1^\ast, \psi)^{-1}
\]
equals $\mathcal{L}^\ast(W_\varphi, f_s, \phi)$. The required identity follows by \eqref{fq gamma}.
\begin{Remark}
 \label{qs rem2}
In the non-split case, in a similar argument as the proof of \cite[(9.7)]{Ka}, we may show that our local zeta integral is equal to the following 
analogue of \eqref{9.7} (up to suitable local $\gamma$-factors):
 \begin{multline*}
 \int_{\overline{U_n}Z_n \backslash G_n}  \int_{Y_n \backslash H_n} \int_{R_{n,\ell}} 
\varphi \left( {}^{\omega_n} rhg \omega_{\ell-n,n}, I_\ell, I_\ell \right) ( \omega_{\psi^{-1}, \Upsilon}( h\omega_ng))\phi(\xi_n, \xi_n)
\\
 \left( \int_{U_n} M^\ast(\tau, s)f_s(\omega_n u \omega_n^{-1}g, I_n, I_n) \psi^{-1}\left(\frac{1}{2}u_{n,1} \right) \, du \right) \, dr \, dh \, dg .
\end{multline*}
As we remarked in the beginning of Section~\ref{proof:s}, in order to simplify the notation, we replaced $\psi$ by $\psi_2$.
Hence, in the formula~\eqref{9.7}, $\psi^{-1}\left(u_{n,1} \right)$ appear. On the other hand, 
from the definition of our intertwining operator \eqref{int op}, $\psi^{-1}\left(\frac{1}{2}u_{n,1} \right)$ appears in the above formula.
\end{Remark}
 %%%%%%%%%%%%%%%%%%%
%
%
%
%
%
%
%
%
%
%
%%%%%%%%%%%%%%%%%%%%
\subsubsection{$r=\ell=n$}
\label{s:r=l=n}
As in Section~\ref{sect 4.1}, collapsing the $du$-integration, using the Iwasawa decomposition 
$\mathrm{GL}_{2\ell} = K_{\mathrm{GL}_{2\ell}} ((\mathrm{GL}_\ell \times \mathrm{GL}_\ell) \rtimes \overline{U_\ell})$
and shifting $a$ to the left, $\mathcal{L}(W_\varphi, f_s, \phi)$ is equal to 
\begin{multline*}
\int_{K_{\mathrm{GL}_{2\ell}}} \int_{\overline{U_\ell}} \int_{ \left( U_{\mathrm{GL}_\ell}\backslash \mathrm{GL}_\ell \right)^2}
\varphi(k, a_1, a_2) \omega_{\psi^{-1}, \mu}(\bar{u}k) \phi(\xi_\ell a_1, \xi_\ell a_2^\ast) f_s(\bar{u}k, a_1, a_2)
\\
\mu(\det a_1 a_2)|\det(a_1 a_2^{-1})|^s \, da_1 \, da_2 \, d\bar{u} \, dk.
\end{multline*}
The inner $da_1$ and $da_2$-integrations form integrals for $\mathrm{GL}_\ell \times \mathrm{GL}_\ell$
and $\tau_1 \times \sigma_1$ and $\tau_2 \times \sigma_2^\ast$.
According to the functional equation \eqref{9.1},
multiplying the last intergral by $\omega_{\tau_1}(-1)^{\ell-1} \gamma(s, (\sigma_1 \otimes \mu) \times \tau_1, \psi)
\omega_{\tau_2}(-1)^{\ell-1} \gamma(s, (\sigma_2 \otimes \mu)^\ast \times \tau_1, \psi)$ leads to
\begin{multline*}
\int_{K_{\mathrm{GL}_{2\ell}}} \int_{\overline{U_\ell}} \int_{ \left( U_{\mathrm{GL}_\ell}\backslash \mathrm{GL}_\ell \right)^2}
\varphi(k, a_1, a_2)( \omega_{\psi^{-1}, \mu}(\bar{u}k) \phi)^\wedge(\xi_\ell a_2, \xi_\ell a_1^\ast)
\\
 f_s(\bar{u}k, a_1, a_2)
\mu(\det a_1 a_2)|\det(a_1 a_2^{-1})|^{s-1} \, da_1 \, da_2 \, d\bar{u} \, dk.
\end{multline*}
Shifting $a_1$ and $a_2$ back, we obtain
\begin{multline*}
\int_{\overline{U_\ell} U_{\mathrm{GL}_\ell} \backslash \mathrm{GL}_{2\ell}}
\varphi(g, I_\ell, I_\ell) \omega_{\psi^{-1}, \mu}(w_\ell g) \phi (\xi_\ell, \xi_\ell ) 
\\
\left( \int_{U_\ell} f_s(w_\ell u w_\ell^{-1}, I_\ell, I_\ell) \psi^{-1} \left(u_{\ell, \ell+1} \right) \, du \right)
\, dg.
\end{multline*}
Now proceed as in the previous section~\ref{9.2},
and this completes the proof of Proposition~\ref{M2} in the split case.
%%%%%%%%%%%%%%%%%%%
%
%
%
%
%
%
%
%
%
%
%%%%%%%%%%%%%%%%%%%%
%%%%%%%%%%%%%%%%%%%%%%%%%%%%%%%%%%%%%%%%%%%%%%%
%
%
%
%
%
%
%
%
%
%
%
%
%
%
%%%%%%%%%%%%%%%%%%%%%%%%%%%%%%%%%%%%%%%%%%%%%%%
%%%%%%%%%%%%%%%%%%%%%%%%%%%%%%%%%%%%%%%%%%%%%%%
%
%
%
%
%
%
%
%
%
%
%
%
%
%
%%%%%%%%%%%%%%%%%%%%%%%%%%%%%%%%%%%%%%%%%%%%%%%
\appendix
\section{Uniqueness of Fourier-Jacobi models for certain representations of Whittaker type}
\label{section:app}
Suppose that $F$ is non-archimedean that $\ell \geq n$.
For a unipotent subgroup $C$ and a character $\chi$ of $C$, 
$J_{C, \chi}$ denote the Jacquet functor with respect to $C$ and $\chi$.
When $E \slash F$ is a quadratic extension of $F$, for a representation $\sigma$ of $G_\ell$, we define
\[
FJ_{\psi, \Upsilon}(\sigma) = 
J_{\mathcal{H}_{n}} \left(J_{U_\ell, \psi_\ell}(\sigma) \otimes \omega_{\psi^{-1}, \Upsilon^{-1}} \right).
\]
This is the Fourier-Jacobi module of $\sigma$ with respect to $\Upsilon$ and $\psi_\ell$ defined in \eqref{def Ul char}
and this gives a representation of $G_n$.

Similarly, in the split case, for a representation $\sigma$ of $\mathrm{GL}_{2\ell}(F)$ and 
for a character $\mu$ of $F^\times$, we define Fourier-Jacobi module 
$\sigma$ with respect to $(\mu, \mu^{-1})$ and $\psi_{\ell, GL}$ defined in \eqref{def Ul char GL} by
\[
FJ_{\psi, \mu}(\sigma) = J_{\mathcal{H}_n} \left( J_{U_{\ell, GL}, \psi_{\ell, GL}}(\sigma)
\otimes \omega_{\psi^{-1}, \mu^{-1}, \mu} \right)
\] 
This is the Fourier-Jacobi module of $\sigma$ and this gives a representation of $\mathrm{GL}_{2n}(F)$.
In this appendix, we prove the following uniqueness. Here, $q$ denotes the cardinality of the residue field of $F$.
\begin{theoremA}
\label{unique}
Let $\pi$ (resp. $\tau$) be a representation of $G_\ell$ (resp. $\mathrm{GL}_n(E)$) of Whittaker type.
Then except a finite number of values of $q^{-s}$, we have 
\[
\dim_\mathbb{C} \mathrm{Hom}_{G_n}(FJ_{\psi, \Upsilon}(\pi), \rho_{\tau, s} ) \leq 1.
\]
Here, we set $\rho_{\tau, s} = \mathrm{Ind}_{Q_n}^{G_n}(\tau|\cdot|^{s-\frac{1}{2}})$.
\end{theoremA}
\begin{theoremA}
\label{GL unique}
Let $\pi$ (resp. $\tau_1, \tau_2$) be a representation of $\mathrm{GL}_{2\ell}(F)$ 
(resp. $\mathrm{GL}_n(F)$) of Whittaker type.
Then except a finite number of values of $q^{-s}$, we have 
\[
\dim_\mathbb{C} \mathrm{Hom}_{\mathrm{GL}_{2n}(F)}(FJ_{\psi, \mu}(\pi), \rho_{\tau_1, \tau_2, s} ) \leq 1.
\]
Here, we set $\rho_{\tau_1, \tau_2, s} = \mathrm{Ind}_{Q_n}^{\mathrm{GL}_{2n}(F)}(\tau_1|\cdot|^{s-\frac{1}{2}} \otimes \tau_2^\ast|\cdot|^{-s+\frac{1}{2}}  )$.
\end{theoremA}
\begin{RemarkA}
Since Theorem~\ref{GL unique} is proved in a similar way as Theorem~\ref{unique},
we only give a proof of Theorem~\ref{unique}.
\end{RemarkA}
\begin{RemarkA}
As remarked in Remark~\ref{verify conv}, in the non-archimedean case, we need these uniqueness in order to define $\gamma$-factors for certain reducible representations
(see also the statements of Proposition~\ref{M1} and Proposition~\ref{M2}).
On the other hand, in the archimedean case, we need these uniqueness only for irreducible representations.
We may expect that similar results as Theorem~\ref{unique}, \ref{GL unique} hold for archimedean fields, but we do not study them in this paper
for our purpose.
\end{RemarkA}
\begin{proof}
For simplicity, we write $\pi_{(\ell)} = J_{U_\ell, \psi_\ell}(\pi)$.
Then our space is written as
\[
\mathrm{Hom}_{G_n \mathcal{H}_n }( \pi_{(\ell)} \otimes \omega_{\psi^{-1}, \Upsilon^{-1}}, \rho_{\tau, s} )
\]
where we extend $\rho_{\tau, s} $ to $G_n \mathcal{H}_n $ so that it is trivial on $\mathcal{H}_n$.
Moreover, by the Frobenius reciprocity, the above space is isomorphic to
\begin{equation}
\label{e:app above sp}
\mathrm{Hom}_{Q_n \mathcal{H}_n}( \pi_{(\ell)} \otimes (\omega_{\psi^{-1}, \Upsilon^{-1}})|_{Q_n \mathcal{H}_n}, \tau \otimes | \cdot |^{s-\frac{1}{2}} )
\end{equation}
Here, recall that 
\[
Q_n \mathcal{H}_n = \left\{ \begin{pmatrix} 1&a&b&c\\ &g&u&b^\prime\\ &&g^\ast&a^\prime\\ &&&1\end{pmatrix} \in G_{2n}\right\}
\]
and $\tau \otimes | \cdot |^{s-\frac{1}{2}}$ is regarded as a representation of $Q_n \mathcal{H}_n$ 
so that it is trivial on $U_n \mathcal{H}_n$.
Moreover, we note that 
\[
(\omega_{\psi^{-1}, \Upsilon^{-1}})|_{Q_n \mathcal{H}_n} = 
\mathrm{ind}_{Q_n \mathcal{H}^\prime_n}^{Q_n \mathcal{H}_n}  (\omega_{\psi^{-1}, \Upsilon^{-1}} |_{Q_n \mathcal{H}_n^\prime}),
\]
where $\mathcal{H}_n^\prime = \left\{(0, y; 0) \in \mathcal{H}_n \right\}$ and $\mathrm{ind}$ denotes the compact induction.
Then by the Frobenius reciprocity, the space \eqref{e:app above sp} is isomorphic to
\[
\mathrm{Hom}_{Q_n \mathcal{H}_n^\prime}( \pi_{(\ell)} \otimes (\omega_{\psi^{-1}, \Upsilon^{-1}})|_{Q_n \mathcal{H}_n^\prime}, 
\tau \otimes | \cdot |^{s-\frac{1}{2}} ).
\]
Let us consider the exact sequence 
\[
0 \rightarrow \mathcal{S}(E^n - \{0\}) \rightarrow \mathcal{S}(E^n) \overset{p}{\longrightarrow} \mathbb{C} \rightarrow 0
\]
where the surjective linear functional $p$ on $\mathcal{S}(E^n)$ is given by $p(\phi) = \phi(0)$.
Then in the same argument as the proof of \cite[Theorem~5.3]{BAS}, the above space is isomorphic to
\[
\mathrm{Hom}_{Q_n \mathcal{H}_n^\prime}
( \pi_{(\ell)} \otimes ind_{Q_n^1 \mathcal{H}_n^\prime}^{Q_n \mathcal{H}_n^\prime}(\Upsilon^{-1} c_\psi), 
\tau \otimes | \cdot |^{s-\frac{1}{2}} )
\]
where $\mathcal{P}_n$ denotes the mirabolic subgroup of $\mathrm{GL}_n(E)$,
$Q_n^1$ denotes the subgroup of $Q_n$ consisting of elements whose Levi part is in $\mathcal{P}_n$,
and we define 
\[
\Upsilon^{-1} c_\psi \left( m \left(\begin{pmatrix}b&e\\&1 \end{pmatrix} \right) 
\begin{pmatrix}1_{n}&v\\ &1_{n} \end{pmatrix} (0, y; 0)\right)
= \Upsilon^{-1}(\det g)\psi^{-1}(v_{n, 1}) \psi^{-1}(y_1).
\]
Hence, by the Frobenius reciprocity, we should consider the space 
\[
\mathrm{Hom}_{Q_n^1 \mathcal{H}_n^\prime}
( \pi_{(\ell)} \otimes \Upsilon^{-1} c_\psi, 
\tau \otimes | \cdot |^{s-\frac{1}{2}} )
\simeq 
\mathrm{Hom}_{\mathcal{P}_n}(
J_{S_n, \psi_{S_n}}(\pi_{(\ell)}), 
\tau \otimes | \cdot |^{s-\frac{1}{2}} ).
\]
Here, $S_n$ denotes the unipotent subgroup 
\[
S_n =\left\{ \begin{pmatrix} 1&&b&\\ &1_n&u&b^\prime\\ &&1_n&\\ &&&1\end{pmatrix} \in G_{2n}\right\}
\]
and 
\[
\psi_{S_n} \left( \begin{pmatrix} 1&&b&\\ &1_n&u&b^\prime\\ &&1_n&\\ &&&1\end{pmatrix}\right)
=\psi^{-1}(u_{n,1}) \psi^{-1}(b_1).
\]
Now, the Jacquet module $J_{S_n, \psi_{S_n}}(\pi_{(\ell)})$ is a $\mathcal{P}_n$-module,
and thus we may consider the Bernstein-Zelevinsky derivative.
Then by the standard argument (for example, see the proof of \cite[Theorem~5.3]{BAS}), 
we see that except a finite number of values of $q^{-s}$,  this space is isomorphic to
\[
\mathrm{Hom}_{\mathcal{P}_n} \left(ind_{U_{\mathrm{GL}_n}(E)}^{\mathcal{P}_n} \psi_n, \tau \right)
\simeq 
\mathrm{Hom}_{U_{\mathrm{GL}_n}(E)} \left(\psi_n, \tau \right).
\]
The last space is the space of Whittaker functionals on $\tau^\ast$, and thus 
our result follows from the uniqueness of Whittaker model since $\tau$ is of Whittaker type.
\end{proof}
%%%%%%%%%%%%%%%%%%%%%%%%%%%%%%%%%%%%%%%%%%%%%%%
%
%
%
%
%
%
%
%
%
%
%
%
%
%
%%%%%%%%%%%%%%%%%%%%%%%%%%%%%%%%%%%%%%%%%%%%%%%

%%%%%%%%%%%%%%%%%%%%%%%%%%%%%%%%%%%%%%%%%%%%%%%
%
%
%
%
%
%
%
%
%
%
%
%
%
%
%%%%%%%%%%arithmetic of automorphic forms%%%%%%%%%%%%%%%%%%%%%%
\end{document}